\documentclass[10.5pt]{article}
\usepackage[leqno]{amsmath}
\usepackage{amsfonts}
\usepackage{graphicx}

\usepackage{amsmath}
\usepackage{amssymb}
\usepackage{latexsym}
\usepackage{amsmath, amsfonts,amssymb, amsthm, euscript,makeidx,color,mathrsfs}

\def\sqr#1#2{{\vcenter{\vbox{\hrule height.#2pt
              \hbox{\vrule width.#2pt height#1pt \kern#1pt \vrule width.#2pt}
              \hrule height.#2pt}}}}
\def\signed #1{{\unskip\nobreak\hfil\penalty50
              \hskip2em\hbox{}\nobreak\hfil#1
              \parfillskip=0pt \finalhyphendemerits=0 \par}}
\def\endpf{\signed {$\sqr69$}}

\def\3n{\negthinspace \negthinspace \negthinspace }
\def\2n{\negthinspace \negthinspace }
\def\1n{\negthinspace }

\def\dbE{\mathbb{E}}
\def\dbF{\mathbb{F}}

\def\dbI{\mathbb{I}}

\def\dbP{\mathbb{P}}
\def\dbQ{\mathbb{Q}}
\def\dbR{\mathbb{R}}
\def\dbS{\mathbb{S}}

\def\sD{\mathscr{D}}

\def\sU{\mathscr{U}}

\def\sX{\mathscr{X}}

\def\={\buildrel \triangle \over =}

\def\ds{\displaystyle}

\def\ns{\noalign{\ss}}
\def\b{\beta}
\def\g{\gamma}
\def\d{\delta}
\def\e{\varepsilon}

\def\l{\lambda}

\def\si{\sigma}
\def\t{\tau}
\def\f{\varphi}
\def\th{\theta}
\def\o{\omega}

\def\G{\Gamma}
\def\D{\Delta}
\def\Th{\Theta}
\def\L{\Lambda}

\def\F{\Phi}
\def\O{\Omega}

\def\cA{{\cal A}}
\def\cB{{\cal B}}
\def\cC{{\cal C}}
\def\cD{{\cal D}}

\def\cF{{\cal F}}
\def\cG{{\cal G}}

\def\cQ{{\cal Q}}
\def\cR{{\cal R}}
\def\cS{{\cal S}}

\def\cU{{\cal U}}

\def\no{\noindent}

\def\ss{\smallskip}
\def\ms{\medskip}
\def\bs{\bigskip}
\def\q{\quad}
\def\qq{\qquad}
\def\hb{\hbox}

\def\liminf{\mathop{\underline{\rm lim}}}

\def\da{\mathop{\downarrow}}

\def\lan{\mathop{\langle}}
\def\ran{\mathop{\rangle}}

\def\esssup{\mathop{\rm esssup}}

\def\h{\widehat}
\def\wt{\widetilde}

\def\cd{\cdot}
\def\cds{\cdots}

\def\as{\hbox{\rm a.s.{ }}}

\def\var{\hbox{\rm var$\,$}}

\def\deq{\mathop{\buildrel\D\over=}}

\def\({\Big (}
\def\){\Big )}
\def\[{\Big[}
\def\]{\Big]}
\def\bde{\begin{definition}}
\def\ede{\end{definition}}
\def\be{\begin{equation}}
\def\bel{\begin{equation}\label}
\def\ee{\end{equation}}
\def\bt{\begin{theorem}}
\def\et{\end{theorem}}
\def\bc{\begin{corollary}}
\def\ec{\end{corollary}}
\def\bl{\begin{lemma}}
\def\el{\end{lemma}}
\def\bp{\begin{proposition}}
\def\ep{\end{proposition}}
\def\bas{\begin{assumption}}
\def\eas{\end{assumption}}
\def\br{\begin{remark}}
\def\er{\end{remark}}
\def\ba{\begin{array}}
\def\ea{\end{array}}
\def\ed{\end{document}}

\def\square#1{\vbox{\hrule\hbox{\vrule height#1%
     \kern#1\vrule}\hrule}}
\def\rectangle#1#2{\vbox{\hrule\hbox{\vrule height#1%
     \kern#2\vrule}\hrule}}

\font\tenbb=msbm10 \font\sevenbb=msbm7 \font\fivebb=msbm5

\newfam\bbfam
\scriptscriptfont\bbfam=\fivebb \textfont\bbfam=\tenbb
\scriptfont\bbfam=\sevenbb

\newtheorem{lemma}{Lemma}[section]
\newtheorem{remark}{Remark}[section]

\newtheorem{theorem}{Theorem}[section]
\newtheorem{corollary}{Corollary}[section]

\newtheorem{definition}{Definition}[section]
\newtheorem{proposition}{Proposition}[section]
\newtheorem{assumption}{Assumption}[section]

\makeatletter
   
   \@addtoreset{equation}{section}
\makeatother

\begin{document}
\title{\bf Linear-Quadratic
Optimal Control Problems\\ for Mean-Field Stochastic Differential
Equations\\ --- Time-Consistent Solutions\footnote{This work is
supported in part by NSF Grant DMS-1007514.}}

\author{Jiongmin Yong\\ Department of
Mathematics, University of Central Florida, Orlando, FL 32816, USA.}

\maketitle

\begin{abstract} Linear-quadratic optimal control problems are considered for mean-field stochastic
differential equations with deterministic coefficients.
Time-inconsistency feature of the problems is carefully
investigated. Both open-loop and closed-loop equilibrium solutions
are presented for such kind of problems. Open-loop solutions are
presented by means of variational method with decoupling of
forward-backward stochastic differential equations, which leads to a
Riccati equation system lack of symmetry. Closed-loop solutions are
presented by means of multi-person differential games, the limit of
which leads to a Riccati equation system with a symmetric structure.

\end{abstract}


\ms


\section{Introduction.}

Let $(\O,\cF,\dbP,\dbF)$ be a complete filtered probability space,
on which a one-dimensional standard Brownian motion $W(\cd)$ is
defined with $\dbF\equiv\{\cF_t\}_{t\ge0}$ being its natural
filtration augmented by all the $\dbP$-null sets. Consider the
following controlled linear stochastic differential equation (SDE,
for short):
\bel{1.1}\left\{\2n\ba{ll}
\ns\ds dX(s)=\[A(s)X(s)+B(s)u(s)\]ds+\[C(s)X(s)+D(s)u(s)\]dW(s),\q
s\in[t,T],\\
\ns\ds X(t)=x,\ea\right.\ee
with a quadratic cost functional
\bel{1.2}J(t,x;u(\cd))=\dbE_t\Big\{\int_t^T\[\lan
Q(s)X(s),X(s)\ran+\lan R(s)u(s),u(s)\ran\]ds+\lan
GX(T),X(T)\ran\Big\}.\ee
Here, $A(\cd),B(\cd),C(\cd),D(\cd),Q(\cd),R(\cd)$ are suitable
matrix-valued (deterministic) functions, $G$ is a matrix, and
$\dbE_t=\dbE[\,\cd\,|\cF_t]$ stands for the conditional expectation
operator. In the above, $X(\cd)$, valued in $\dbR^n$, is called the
{\it state process}, $u(\cd)$, valued in $\dbR^m$, is called the
{\it control process}, and $(t,x)\in\sD$ is called the {\it initial
pair} where
$$\sD=\Big\{(t,x)\bigm|t\in[0,T],~x\hb{ is $\cF_t$-measurable, }\dbE|x|^2<\infty\Big\}.$$
Denote
$$\ba{ll}\sU[t,T]\equiv L^2_\dbF(t,T;\dbR^m)=\Big\{u:[t,T]\times\O\to\dbR^m\bigm|u(\cd)
\hb{ is $\dbF$-progressively measurable,}\\
\ns\ds\qq\qq\qq\qq\qq\qq\qq\qq\qq\qq\qq
\dbE\int_t^T|u(s)|^2ds<\infty\Big\},\ea$$
which is the set of all control processes. Under some mild
conditions on the coefficients, for any initial pair $(t,x)\in\sD$
and a control $u(\cd)\in\sU[t,T]$, the state equation (\ref{1.1})
admits a unique solution $X(\cd)=X(\cd\,;t,x,u(\cd))$, and the cost
functional $J(t,x;u(\cd))$ is well-defined. A standard
linear-quadratic (LQ, for short) optimal control problem can be
stated as follows.

\ms

\bf Problem (LQ). \rm For any given initial pair $(t,x)\in\sD$, find
a $u^*(\cd)\in\sU[t,T]$ such that
\bel{1.3}J(t,x;u^*(\cd))=\inf_{u(\cd)\in\sU[t,T]}J(t,x;u(\cd))\deq
V(t,x).\ee
Any $u^*(\cd)\in\sU[t,T]$ satisfying (\ref{1.3}) is called an {\it
optimal control} for the given initial pair $(t,x)$, the
corresponding state process $X^*(\cd)$ is called an {\it optimal
state process} for $(t,x)$, $(X^*(\cd),u^*(\cd))$ is called an {\it
optimal pair} for $(t,x)$, and $V(\cd\,,\cd)$ is called the {\it
value function} of Problem (LQ).

\ms

It is well-known that (\cite{Yong-Zhou 1999}) under proper
conditions, the above Problem (LQ) admits a unique optimal pair
$(X^*(\cd),u^*(\cd))$. Moreover, the following {\it Riccati
equation}
\bel{}\left\{\2n\ba{ll}
\ns\ds\dot P(s)\1n+\1n P(s)A(s)\1n+\1n A(s)^TP(s)\1n+\1n
C(s)^TP(s)C(s)\1n+\1n Q(s)\\
\ns\ds\q-\big[P(s)B(s)\1n+\1n C(s)^T\1n P(s)D(s)\big]\1n
\big[R(s)\1n+\1n D(s)^T\1n P(s)D(s)\big]\1n^{-\1n1}\big[B(s)^T\1n
P(s)\1n+\1n D(s)^T\1n P(s)C(s)\big]\1n=\1n0,\\
\ns\ds\qq\qq\qq\qq\qq\qq\qq\qq\qq\qq\qq\qq\qq\qq\qq\qq s\in[0,T],\\
\ns\ds P(T)=G,\\
\ns\ds\det\big[R(s)+D(s)^TP(s)D(s)\big]>0,\qq
s\in[0,T],\ea\right.\ee
has a unique solution $P(\cd)$ which is positive definite, and the
optimal control $u^*(\cd)$ admits the following state feedback
representation:
\bel{}u^*(s)=-\Th(s)X^*(s),\qq\qq s\in[t,T],\ee
with
\bel{}\Th(s)=\big[R(s)+D(s)^TP(s)D(s)\big]^{-1}\big[B(s)^TP(s)+D(s)^TP(s)C(s)
\big],\qq s\in[0,T],\ee
and $X^*(\cd)$ being the solution to the following closed-loop
system:
\bel{closed-loop00}\left\{\2n\ba{ll}
\ns\ds
dX^*(s)=\big[A(s)-B(s)\Th(s)\big]X^*(s)ds+\big[C(s)-D(s)\Th(s)\big]X^*(s)dW(s),
\q s\in[t,T],\\
\ns\ds X^*(t)=x,\ea\right.\ee
Further, the value function $V(\cd\,,\cd)$ is explicitly given by
\bel{V}V(t,x)=\lan P(t)x,x\ran,\qq\forall(t,x)\in\sD.\ee
Next, if $\F(\cd\,;\cd)$ is the {\it fundamental matrix} of
(\ref{closed-loop00}), i.e.,
\bel{}\left\{\2n\ba{ll}
\ns\ds
d\F(s;t)=\big[A(s)-B(s)\Th(s)\big]\F(s;t)ds+\big[C(s)-D(s)\Th(s)\big]\F(s;t)dW(s),
\q s\in[t,T],\\
\ns\ds\F(t;t)=I,\ea\right.\ee
then
\bel{}X^*(s)=\F(s;t)x,\qq\forall s\in[t,T],~(t,x)\in\sD,\ee
and
\bel{}u^*(s)\equiv u^*(s;t,x)=-\Th(s)\F(s;t)x,\qq
s\in[t,T],~(t,x)\in\sD.\ee
Note that $\Th(\cd)$ is independent of the initial pair $(t,x)$.
According to the above representation of the optimal pair, if we
start from the initial pair $(\t,X^*(\t))\in\sD$ with $\t\in(t,T]$
being any $\dbF$-stopping time, then we have
$$\ba{ll}
\ns\ds u^*(s;\t,X^*(\t))=-\Th(s)\F(s;\t)X^*(\t)=-\Th(s)\F(s;\t)\F(\t;t)x\\
\ns\ds\qq\qq\qq~=-\Th(s)\F(s;t)x=u^*(s;t,x),\qq s\in[\t,T].\ea$$
This means that the restriction $u^*(\cd\,;t,x)\big|_{[\t,T]}$ of
the optimal control $u^*(\cd\,;t,x)$ for the initial pair
$(t,x)\in\sD$ on a later time interval $[\t,T]$ coincides with the
optimal control $u^*(\cd\,;\t,X^*(\t))$ for the initial pair
$(\t,X^*(\t))\in\sD$. Such a phenomenon is called the {\it
time-consistency} of Problem (LQ).

\ms

We now look at the situation for the following cost functional:
\bel{1.9}\ba{ll}
\ns\ds J^\l(t,x;u(\cd))=\dbE_t\Big\{\int_t^Te^{-\l(s-t)}\[\lan
Q(s)X(s),X(s)\ran+\lan R(s)u(s),u(s)\ran\]ds\\
\ns\ds\qq\qq\qq\qq\qq+e^{-\l(T-t)}\lan GX(T),X(T)\ran\Big\},\ea\ee
for some constant $\l>0$ which is referred to as a {\it discount
rate}. We also call the function $s\mapsto e^{-\l s}$ an {\it
exponential discounting}. If we denote
$$\wt Q(s)=e^{-\l s}Q(s),\q\wt R(s)=e^{-\l s}R(s),\q\wt G=e^{-\l
T}G,$$
then
\bel{1.10}\ba{ll}
\ns\ds\wt J(t,x;u(\cd))\deq\dbE_t\int_t^T\[\lan\wt
Q(s)X(s),X(s)\ran+\lan
\wt R(s)u(s),u(s)\ran\]ds+\lan\wt GX(T),X(T)\ran\]\\
\ns\ds\qq\qq\q=e^{-\l t}J^\l(t,x;u(\cd)).\ea\ee
Therefore, the optimal controls for the problems with resepective
cost functionals (\ref{1.9}) and (\ref{1.10}) coincide. Since the
problem with cost functional (\ref{1.10}) is time-consistent, so is
the problem with cost functional (\ref{1.9}).

\ms

Note that the significant difference between the cost functionals
$J^\l(t,x;u(\cd))$ defined by (\ref{1.9}) and $J(t,x;u(\cd))$
defined by (\ref{1.2}) is that the running cost rate and the
terminal cost of $J^\l(t,x;u(\cd))$ explicitly depend on the initial
time moment $t$. However, from the above, we see that such a special
dependence on $t$ does not change the time-consistency of the problem.
The reason is that exponential discounting exhibits a time-consistent memory
effect.

\ms

It then will be interesting to know what will happen if we replace
the exponential discounting by some general discounting function. To
get some feeling, let us look at the following simple example.

\ms

\bf Example 1.1. \rm Consider a one-dimensional controlled linear
SDE:
\bel{ex2.1a}\left\{\2n\ba{ll}
\ns\ds dX(s)=u(s)ds+X(s)dW(s),\qq s\in[t,T],\\
\ns\ds X(t)=x,\ea\right.\ee
with cost functional
\bel{ex2.1b}J(t,x;u(\cd))=\dbE_t\[\int_t^T\rho(s,t)|u(s)|^2ds+g(t)|X(T)|^2\],\ee
where $\rho(s,t)$ and $g(t)$ are deterministic non-constant,
continuous and positive functions. Let $P(\cd\,,t)$ be the unique
solution to the following Riccati differential equation:
\bel{P(s,t)}\left\{\2n\ba{ll}
\ns\ds P_s(s,t)+P(s,t)-{P(s,t)^2\over\rho(s,t)}=0,\qq s\in[t,T],\\
\ns\ds P(T,t)=g(t).\ea\right.\ee
Observe the following:
$$\dbE_t\[g(t)|X(T)|^2\]-P(t,t)x^2=\dbE_t\int_t^T\(P_s(s,t)X(s)^2
+2P(s,t)X(s)u(s)+P(s,t)X(s)^2\)ds.$$
Hence,
$$\ba{ll}
\ns\ds J(t,x;u(\cd))=\dbE_t\[\int_t^T\rho(s,t)|u(s)|^2ds+g(t)|X(T)|^2\]\\
\ns\ds\qq\qq\q=P(t,t)x^2+\dbE_t\int_t^T\(\rho(s,t)|u(s)|^2+2P(s,t)X(s)u(s)
+\big[P_s(s,t)+P(s,t)\big]X(s)^2\)ds\\
\ns\ds\qq\qq\q=P(t,t)x^2+\dbE_t\int_t^T\rho(s,t)\Big|u(s)+{P(s,t)\over
\rho(s,t)}X(s)\Big|^2ds.\ea$$
Let
\bel{u(s,t,x)}u^*(s;t,x)=-{P(s,t)\over\rho(s,t)}X^*(s;t,x),\qq
s\in[t,T],\ee
with $X^*(\cd\,;t,x)$ being the solution of the following
closed-loop system:
\bel{closed-loop1}\left\{\2n\ba{ll}
\ns\ds
dX^*(s;t,x)=-{P(s,t)\over\rho(s,t)}X^*(s;t,x)ds+X^*(s;t,x)dW(s),\qq
s\in[t,T],\\
\ns\ds X^*(t;t,x)=x.\ea\right.\ee
Then
$$\inf_{u(\cd)\in\sU[t,T]}J(t,x;u(\cd))=J(t,x;u^*(\cd\,;t,x))=P(t,t)x^2.$$
Hence, for any $\t\in(t,T)$, if we let $P(\cd\,,\t)$ be the solution
of (\ref{P(s,t)}) with $t$ replaced by $\t$, then
$$\ba{ll}
\ns\ds
J(\t,X^*(\t;t,x);u(\cd))=\dbE_\t\[\int_\t^T\rho(s,\t)|u(s)|^2ds
+g(\t)|\wt X(T)|^2\]\\
\ns\ds\qq\qq\qq\qq\q\,=P(\t,\t)X^*(\t;t,x)^2
+\dbE_\t\int_\t^T\rho(s,\t)\big|u(s)+{P(s,\t)\over\rho(s,\t)}\wt
X(s)\big|^2ds,\ea$$
where
$$\wt X(\cd)\equiv X(\cd\,;\t,X^*(\t;t,x),u(\cd)),\qq u(\cd)\in\cU[\t,T].$$
Consequently, by denoting
$$\wt X^*(\cd)=X(\cd,;\t,X^*(\t;t,x),u^*(\cd)\big|_{[\t,T]}),$$
we have
\bel{1.17}\ba{ll}
\ns\ds
J(\t,X^*(\t;t;x);u^*(\cd)\big|_{[\t,T]})=P(\t,\t)X^*(\t;t,x)^2
+\dbE_\t\int_\t^T\rho(s,\t)|u^*(s;t,x)+{P(s,\t)\over\rho(s,\t)}\wt X^*(s)|^2ds\\
\ns\ds\ge P(\t,\t)X^*(\t;t,x)^2=\inf_{u(\cd)\in\sU[\t,T]}
J(\t,X^*(\t;t,x),u(\cd))=J(\t,X^*(\t;t,x),u^*(\cd\,;\t,X^*(\t;t,x))),\ea\ee
and the inequality actually becomes an equality if and only if
\bel{1.19}\left\{\2n\ba{ll}
\ns\ds u^*(s;t,x)=u^*(s;\t,X^*(\t;t,x)),\\
\ns\ds X^*(s;t,x)=X^*(s;\t,X^*(\t;t,x)),\ea\right.\qq s\in[\t,T].\ee
Similar to (\ref{u(s,t,x)}), we should also have
\bel{u(s,t,x)0}u^*(s;\t,X^*(\t;t,x))=-{P(s,\t)\over\rho(s,\t)}
X^*(s;\t,X^*(\t;t,x)),\qq s\in[\t,T],\ee
Combining (\ref{u(s,t,x)}) and (\ref{1.19})--(\ref{u(s,t,x)0}), we
obtain
\bel{}{P(s,\t)\over\rho(s,\t)}={P(s,t)\over\rho(s,t)},\qq\forall
t\le\t\le s\le T.\ee
In particular, evaluating at $s=T$, one has
$${g(\t)\over\rho(T,\t)}={g(t)\over\rho(T,t)}.$$
This means that if $t\mapsto{g(t)\over\rho(T,t)}$ is not a constant
function, then a strict inequality in (\ref{1.17}) will hold for
some $\t\in(t,T]$, which means that the restriction
$u^*(\cd)\big|_{[\t,T]}$ of the optimal control $u^*(\cd)$ (for the
initial pair $(t,x)$) is not optimal on $[\t,T]$ with initial pair
$(\t,X^*(\t;t,x))\in\sD$. Such a feature is called {\it
time-inconsistency} of the problem.

\ms

Next example exhibits a different situation.

\ms

\bf Example 1.2. \rm Consider state equation (\ref{ex2.1a}) with
cost functional
\bel{1.18}J(t,x;u(\cd))=\dbE_t\[\int_t^T|u(s)|^2ds+|\dbE_t[X(T)]|^2\].\ee
Note that in the above, $\dbE_t[X(T)]$ is nonlinearly appeared,
which is non-classical. To solve the corresponding LQ problem, let
us introduce the following Riccati equation:
$$\left\{\2n\ba{ll}
\ns\ds\dot{\h P}(s)-\h P(s)^2=0,\qq s\in[0,T],\\
\ns\ds\h P(T)=1,\ea\right.$$
whose solution is given by
$$\h P(s)={1\over T-s+1},\qq s\in[0,T].$$
For any given $(t,x)\in\sD$ and $u(\cd)\in\sU[t,T]$, let
$X(\cd)=X(\cd\,;t,x,u(\cd))$ be the corresponding state process.
Then
$$\left\{\2n\ba{ll}
\ns\ds d\dbE_t[X(s)]=\dbE_t[u(s)]ds,\qq s\in[t,T],\\
\ns\ds\dbE_t[X(t)]=x.\ea\right.$$
Observe the following:
$$\ba{ll}
\ns
J(t,x;u(\cd))\3n\2n&\ds=\dbE_t\[\int_t^T|u(s)|^2ds+\h P(T)|\dbE_t[X(T)]|^2\]\\
\ns&\ds=\h P(0)|x|^2+\dbE_t\int_t^T\(|u(s)|^2+\h
P(s)^2\big|\dbE_t[X(s)]\big|^2
+2\h P(s)\dbE_t[X(s)]u(s)\)ds\\
\ns&\ds=\h P(0)|x|^2+\dbE_t\int_t^T\big|u(s)+\h
P(s)\dbE_t[X(s)]\big|^2ds\ge\h P(0)|x|^2,\ea$$
with the equality holds when
$$u^*(s)=-\h P(s)\dbE_t[X^*(s)],\qq s\in[t,T],$$
where $X^*(\cd)$ is the solution to the following closed-loop
system:
$$\left\{\2n\ba{ll}
\ns\ds dX^*(s)=-\h P(s)\dbE_t[X^*(s)]ds+X^*(s)dW(s),\qq s\in[t,T],\\
\ns\ds X^*(t)=x.\ea\right.$$
Such an equation is called a mean-field SDE. The above yields
$$\left\{\2n\ba{ll}
\ns\ds d\big(\dbE_t[X^*(s)]\big)=-\h P(s)\dbE_t[X(s)]ds,\qq
s\in[t,T],\\
\ns\ds\dbE_t[X^*(t)]=x.\ea\right.$$
Consequently,
$$\dbE_t[X^*(s)]=e^{-\int_t^s\h P(r)dr}x=e^{-\int_t^s{dr\over
T-r+1}}x={T-s+1\over T-t+1}x,\qq s\in[t,T].$$
Then
$$u^*(s)\equiv u^*(s;t,x)=-\h P(s)\dbE_t[X^*(s)]=-{x\over T-t+1}\,,\qq s\in[t,T].$$
Note that $u^*(\cd)$ stays as a constant on $[t,T]$, with $(t,x)$ as
a parameter. With such a control, the closed-loop system reads
$$\left\{\2n\ba{ll}
\ns\ds dX^*(s)=-{x\over T-t+1}\,ds+X^*(s)dW(s),\qq s\in[t,T],\\
\ns\ds X^*(t)=x,\ea\right.$$
whose solution is given explicitly by the following:
$$X^*(s)=e^{-{1\over2}(s-t)+W(s)-W(t)}x+{x\over
T-t+1}\int_t^se^{-{1\over2}(s-r)+W(s)-W(r)}dr,\qq s\in[t,T].$$
Consequently, for any $\t\in(t,T]$, the optimal control starting
from $(\t,X^*(\t))$ should be
$$u^*(s;\t,X^*(\t))=-{X^*(\t)\over T-\t+1},\qq s\in[\t,T],$$
which is $\cF_\t$-measurable, and is not $\cF_t$-measurable for any
$t<\t$, whereas $u^*(s;t,x)$ is merely $\cF_t$-measurable, with
$t<\t$. Hence,
$$u^*(s;\t,X^*(\t))=u^*(s;t,x),\qq s\in[\t,T],~\as,$$
cannot be true. This implies that the problem associated with
(\ref{ex2.1a}) and (\ref{1.18}) is time-inconsistent.

\ms

The above two examples present two different situations for which
the associated problems are time-inconsistent. In the first example,
the time-inconsistency is due to the fact that the discounting is
non-exponential, and in the second example, the time-inconsistency
is due to the appearance of the conditional expectation of the state
in a nonlinear way. The non-exponential discounting situation has
been discussed in \cite{Ekeland-Lazrak 2010,Ekeland-Mbodji-Pirvu
2012,Ekeland-Privu 2008a,Yong 2011,Yong 2012a,Yong
2012b,Marin-Solano-Navas 2010}. The main motivation of that is
trying to catch people's subjective preferences on the discounting
(see \cite{Strotz,Pollak 1968,Goldman 1980,Miller-Salmon 1985,Tesfatsion 1986,
L1997,Karp-Lee 2003,Caplin-Leahy 2006,Herings-Rohde,Grenadier-Wang 2007,
Marin-Solano-Shevkoplyas}). The case with nonlinear appearance of conditional expectation in the
terminal cost of the cost functional has been studied in
\cite{Bjork-Murgoci 2008,Bjork-Murgoci-Zhou,Hu-Jin-Zhou 2011}. A
main motivation of that is a kind of conditional/dynamic mean-variance
problem (see \cite{Basak-Chabakarui 2008,Borkar-Kumar 2010}) See also
\cite{Buckdahn-Li-Peng 2009,Lorenz 2010,Carmona-Delarue 2012}
for relevant results. On the other
hand, in \cite{Yong 2012c}, an LQ problem for
mean-field SDEs in a fixed time-duration was investigated (without
mentioning the time-inconsistency), for which a motivation is from
optimal control theory point of view: People might like to have the
optimal control, as well as the optimal state, to be not too
``random''. To achieve that, one could include conditional variation
$\var_t[X(s)]$ of the state process and/or conditional variation
$\var_t[u(s)]$ of the control process in the cost functional, where
$$\ba{ll}
\ns\ds\var_t[X(s)]=\dbE_t\big|X(s)-\dbE_t[X(s)]\big|^2=\dbE_t|X(s)|^2
-\big|\dbE_t[X(s)]\big|^2,\\
\ns\ds\var_t[u(s)]=\dbE_t\big|u(s)-\dbE_t[u(s)]\big|^2=\dbE_t|u(s)|^2
-\big|\dbE_t[u(s)]\big|^2.\ea$$
Then $\dbE_t[X(s)]$ and $\dbE_t[u(s)]$ naturally appear
quadratically.

\ms

The purpose of this paper is to study a general linear quadratic
optimal control problem for mean-field stochastic differential
equations (MF-SDEs, for short), which include the cases of
non-exponential discounting and conditional mean-variance problems.
We point out that the problem can be stated for MF-SDEs with random
coefficients. However, at the moment, we only have a relatively
complete theory for deterministic coefficient case. Therefore, for
the definiteness of the presentation, we restrict ourselves to the
deterministic coefficients. A theory for general random coefficient
case will be reported in a future publication.

\ms

The rest of the paper is organized as follows. Section 2 will
collect some basic preliminary results. In Section 3, we formulate
our linear-quadratic optimal control problem for mean-field SDEs
(MF-LQ problem for short), and the so-called pre-commitment solution
to such a problem (which is a result from \cite{Yong
2012c}). In Section 4, we introduce open-loop and closed-loop
time-consistent solutions to our MF-LQ problem. We derive the
Riccati equation system for open-loop equilibrium solution by means
of variational method and decoupling a forward-backward stochastic
differential equation (FBSDE, for short). We show that as long as
the Riccati equation system has a solution, the open-loop
equilibrium control can be constructed. It turns out that the Riccati
equation system for open-loop equilibrium solution does not have a
desired symmetric structure. Therefore, it seems to be difficult to
establish the well-posedness for such a system in general. Further,
in this section, we state the main result of the paper: the existence of
closed-loop equilibrium solution via a Riccati equation system which has a
symmetric structure. Some interesting special cases are indicated.
In Section 5, a multi-person differential game is introduced from
which we find a kind of approximate time-consistent strategy.
Convergence of the approximate time-consistent strategy is established
in Section 6, which leads to the Riccati equation system and its
well-posedness.

\section{Preliminaries.}

In this section, we present some preliminary results which will be
useful in the following sections. Let us first introduce some
notations. For any $0\le t<T$, let
$$\ba{ll}
\ns\ds\sX_t=\big\{\xi:\O\to\dbR^n\bigm|\xi\hb{ is
$\cF_t$-measurable, }
\dbE_t|\xi|^2<\infty\big\},\\
\ns\ds\sX[t,T]=\Big\{X:[t,T]\times\O\to\dbR^n\bigm|X(\cd)\hb{ is
$\dbF$-adapted, $t\mapsto X(t)$ is continuous,
}\dbE_t\[\sup_{s\in[t,T]}|X(s)|^2\]<\infty\Big\}.\ea$$
For any Euclidean space $H=\dbR^n,\dbR^{n\times m},\dbS^n$ (where
$\dbS^n$ is the set of all $(n\times n)$ symmetric matrices), we let
$$\ba{ll}
\ns\ds C\big([t,T];H\big)=\big\{\f:[t,T]\to H\bigm|\f(\cd)\hb{ is
continuous}\big\},\\
\ns\ds C^1\big([t,T];H\big)=\big\{\f:[t,T]\to
H\bigm|\f(\cd),\dot\f(\cd)\hb{ are
continuous}\big\},\\
\ns\ds L^\infty(t,T;H)=\big\{\f:[t,T]\to
H\bigm|\|\f(\cd)\|_\infty\equiv\esssup_{s\in[t,T]}|\f(s)|<\infty\big\}.\ea$$

\ms

Now, we consider the following linear mean-field SDE (MF-SDE, for
short):
\bel{MF-SDE1}\left\{\2n\ba{ll}
\ns\ds dX(s)\1n=\2n\[\cA(s)X(s)\1n+\1n\bar\cA(s)\dbE_t[X(s)]\1n+\1n
b(s)\]ds\1n+\1n\[\cC(s)X(s)\1n+\1n\bar\cC(s)
\dbE_t[X(s)]\1n+\1n\si(s)\]dW(s),\q s\1n\in\1n[t,T],\\
\ns\ds X(t)=x,\ea\right.\ee
with suitable coefficients
$\cA(\cd),\bar\cA(\cd),\cC(\cd),\bar\cC(\cd)$, and nonhomogeneous
terms $b(\cd),\si(\cd)$, as well as the initial pair $(t,x)\in\sD$.
We refer to $\{\cA(\cd),\bar\cA(\cd),\cC(\cd),\bar\cC(\cd)\}$ as the
{\it generator} of the linear MF-SDE (\ref{MF-SDE1}).

\ms

The following result is concerned with the well-posedness of
(\ref{MF-SDE1}).

\ms

\bf Proposition 2.1. \sl Let
$\cA(\cd),\bar\cA(\cd),\cC(\cd),\bar\cC(\cd)$ be in
$L^\infty(0,T;\dbR^{n\times n})$. Then for any $(t,x)\in\sD$ and
$b(\cd),\si(\cd)\in L^2_\dbF(0,T;\dbR^n)$, MF-SDE $(\ref{MF-SDE1})$
admits a unique solution $X(\cd)\in\sX[t,T]$. Moreover, the
following variation of constants formula holds:
\bel{X(s)}\ba{ll}
\ns\ds
X(s)\1n=\1n\Psi(s,t)\Big\{\1n\big[I\1n+\2n\int_t^s\1n\Psi(r,t)^{-1}\1n\big(\bar\cA(r)\1n
-\1n\cC(r)\bar\cC(r)\big)\bar\F(r,t)dr
\1n+\2n\int_t^s\1n\Psi(r,t)^{-1}\bar\cC(r)\bar\F(r,t)dW(r)\big]x\\
\ns\ds\qq\qq+\int_t^s\Psi(r,t)^{-1}\[b(r)+\(\bar\cA(r)-\cC(r)\bar\cC(r)\)\int_t^r\bar\F(r,t)
\dbE_t[b(\t)]d\t-\cC(r)\si(r)\]dr\\
\ns\ds\qq\qq+\int_t^s\Psi(r,t)^{-1}\[\bar\cC(r)\int_t^r\bar\F(r,\t)\dbE_t[b(\t)]d\t
+\si(r)\]dW(r)\Big\},\qq s\in[t,T],\ea\ee
where $\bar\F(\cd\,,\cd)$ is the fundamental matrix of
$\cA(\cd)+\bar\cA(\cd)$, and $\Psi(\cd\,,\cd)$ is the fundamental
matrix of $\{\cA(\cd),\cC(\cd)\}$, i.e.,
\bel{}\left\{\2n\ba{ll}
\ns\ds\dot{\bar\F}(s,t)=\big[\cA(s)+\bar\cA(s)\big]\bar\F(s,t),\qq
s\in[t,T],\\
\ns\ds\bar\F(t,t)=I,\ea\right.\ee
and
\bel{}\left\{\2n\ba{ll}
\ns\ds d\Psi(s,t)=\cA(s)\Psi(s,t)ds+\cC(s)\Psi(s,t)dW(s),\qq
s\in[t,T],\\
\ns\ds\Psi(t,t)=I.\ea\right.\ee
%

%
%
%
%
%
%
%
%
%

\ms

\it Proof. \rm First, we introduce the following linear ordinary
differential equation (ODE, for short):
\bel{ODE1}\left\{\2n\ba{ll}
\ns\ds d\bar X(s)=\big\{\big[\cA(s)+\bar\cA(s)\big]\bar
X(s)+\dbE_t[b(s)]\big\}ds,\qq s\in[t,T],\\
\ns\ds\bar X(t)=x,\ea\right.\ee
whose solution is given by
$$\bar X(s)=\bar\F(s,t)x+\int_t^s\bar\F(s,r)\dbE_t[b(r)]dr,\qq s\in[t,T],$$
where $\bar\F(\cd\,,\cd)$ is the fundamental matrix of
$\cA(\cd)+\bar\cA(\cd)$ (which is deterministic). Then we consider
the following linear non-homogeneous SDE:
$$\left\{\2n\ba{ll}
\ns\ds dX(s)\1n=\1n\[\cA(s)X(s)\1n+\1n\bar\cA(s)\bar X(s)\1n+\1n
b(s)\]ds\1n+\1n\[\cC(s)X(s)\1n+\1n\bar\cC(s)
\bar X(s)\1n+\1n\si(s)\]dW(s),\q s\1n\in\1n[t,T],\\
\ns\ds X(t)=x,\ea\right.$$
It is standard that the above admits a unique solution $X(\cd)\equiv
X(\cd\,;t,x,b(\cd),\si(\cd))\in\sX[t,T]$. Clearly, for such a
solution, one has
$$\dbE_t[X(s)]=\bar X(s),\qq s\in[t,T].$$
Hence, $X(\cd)$ is the unique solution to (\ref{MF-SDE1}). The rest
of the conclusion follows easily. \endpf

\ms

It is known that if $X(\cd\,;t,x)$ is the solution to the following
SDE:
\bel{}\left\{\2n\ba{ll}
\ns\ds dX(s)\1n=\1n\[\cA(s)X(s)\1n+\1n
b(s)\]ds\1n+\1n\[\cC(s)X(s)\1n+\1n\si(s)\]dW(s),\q s\1n\in\1n[t,T],\\
\ns\ds X(t)=x,\ea\right.\ee
for any given $(t,x)\in\sD$, then
\bel{X=X}X(s;\t,X(\t;t,x))=X(s;t,x),\qq t\le\t\le s\le T.\ee
In another word, the map $\t\mapsto(\t,X(\t;t,x))$ has the so-called
semigroup property. Whereas, we point out that this is not
necessarily the case for MF-SDEs, in general. Here is a simple
example.

\ms

\bf Example 2.2. \rm Consider
$$\left\{\ba{ll}
\ns\ds dX(s)=\dbE_t[X(s)]ds+\dbE_t[X(s)]dW(s),\qq s\in[t,T],\\
\ns\ds X(t)=x.\ea\right.$$
Then
$$d\dbE_t[X(s)]=\dbE_t[X(s)]ds,\qq s\in[t,T],$$
which leads to
$$\dbE_t[X(s)]=e^{s-t}x,\qq s\in[t,T].$$
Hence, the original equation becomes
$$\left\{\ba{ll}
\ns\ds dX(s)=e^{s-t}xds+e^{s-t}xdW(s),\qq s\in[t,T],\\
\ns\ds X(t)=x.\ea\right.$$
Then one has
$$\ba{ll}
\ns\ds X(s;t,x)=\[1+\int_t^se^{r-t}dr+\int_t^se^{r-t}dW(r)\]x\\
\ns\ds\qq\qq=\[e^{s-t}+\int_t^se^{r-t}dW(r)\]x,\qq s\in[t,T].\ea$$
Now, for any $t<\t<s$, we have
$$\ba{ll}
\ns\ds X(s;\t,X(\t;t,x))=\[e^{s-\t}+\int_\t^se^{r-\t}dW(r)\]X(\t;t,x)\\
\ns\ds\qq\qq\qq\qq=\[e^{s-\t}+\int_\t^se^{r-\t}dW(r)\]\[e^{\t-t}+\int_t^\t
e^{r-t}dW(r)\]x\\
\ns\ds\qq\qq\qq\qq=\[e^{s-t}+\int_\t^s e^{r-t}dW(r)+\int_t^\t
e^{s-\t+r-t}dW(r)+\(\int_\t^se^{r-\t}dW(r)\)\(\int_t^\t
e^{r-t}dW(r)\)\]x.\ea$$
Hence,
$$X(s;\t,X(\t;t,x))-X(s;t,x)=\[e^{s-\t}+\int_\t^se^{r-\t}dW(r)\]\(\int_t^\t
e^{r-t}dW(r)\)x.$$
Consequently, as long as $x\ne0$,
$$\ba{ll}
\ns\ds\dbE_t|X(s;\t,X(\t;t,x))-X(s;t,x)|^2=\dbE_t\Big\{\dbE_\t\[e^{s-\t}+\int_\t^se^{r-\t}dW(r)\]^2\(\int_t^\t
e^{r-t}dW(r)\)^2\Big\}x^2\\
\ns\ds\qq\qq\qq\qq\qq\qq\qq\q=\[e^{2(s-\t)}+\int_\t^se^{2(r-\t)}dr\]\(\int_t^\t
e^{2(r-t)}dr\)x^2>0.\ea$$
This shows that (\ref{X=X}) fails in general.

\ms

The following result will be useful later.

\ms

\bf Proposition 2.3. \sl Let
\bel{}\left\{\2n\ba{ll}
\ns\ds\cA(\cd),\cC(\cd)\in
L^\infty(t,T;\dbR^n),\q\cB(\cd),\cD(\cd),\cS(\cd)^T\in
L^\infty(t,T;\dbR^{n\times m}),\\
\ns\ds\cQ(\cd)\in L^\infty(t,T;\dbS^n),\q\cR(\cd)\in
L^\infty(t,T;\dbS^m),\q\cG\in\dbS^n.\ea\right.\ee
Suppose that for some $\d>0$,
\bel{2.12}\left\{\2n\ba{ll}
\ns\ds\cR(s)\ge\d I,\q\cQ(s)-\cS(s)^T\cR(s)^{-1}\cS(s)\ge0,\qq
s\in[t,T],\\
\ns\ds\cG\ge0.\ea\right.\ee
Then the following Riccati equation admits a unique solution
$P(\cd)\in C^1([t,T];\dbS^n)$:
\bel{Riccati0}\left\{\2n\ba{ll}
\ns\ds\dot P(s)+P(s)\cA(s)+\cA(s)^TP(s)+\cC(s)^TP(s)\cC(s)+\cQ(s)\\
\ns\ds\q-\big[P(s)\cB(s)+\cS(s)^T+\cC(s)^TP(s)\cD(s)\big]\big[\cR(s)+\cD(s)^TP(s)
\cD(s)\big]^{-1}\\
\ns\ds\qq\qq\cd\big[\cB(s)^TP(s)+\cS(s)+\cD(s)^TP(s)\cC(s)\big]=0,\qq
s\in[t,T],\\
\ns\ds P(T)=\cG,\\
\ns\ds\cR(s)+\cD(s)^TP(s)\cD(s)>0.\ea\right.\ee
Further,
\bel{}0\le P(s)\le\Pi(s)\le K_0I,\qq s\in[t,T],\ee
where $\Pi(\cd)\in C^1([0,T];\dbS^n)$ is the unique solution to the
following Lyapunov equation:
\bel{Lyapunov1}\left\{\2n\ba{ll}
\ns\ds\dot\Pi(s)+\Pi(s)\cA(s)+\cA(s)^T\Pi(s)+\cC(s)^T\Pi(s)\cC(s)+\cQ(s)=0,\qq
s\in[t,T],\\
\ns\ds\Pi(T)=\cG,\ea\right.\ee
and
$$K_0=\(|\cG|+T\|\cQ(\cd)\|_\infty\)e^{T\|\cA(\cd)+\cA(\cd)^T
+\cC(\cd)^T\cC(\cd)\|_\infty}.$$

\ms

\it Proof. \rm For any $(r,x)\in[t,T)\times\dbR^n$, we introduce the
following controlled SDE:
\bel{2.14}\left\{\2n\ba{ll}
\ns\ds
dX(s)=\[\cA(s)X(s)+\cB(s)u(s)\]ds+\[\cC(s)X(s)+\cD(s)u(s)\]dW(s),\qq
s\in[r,T],\\
\ns\ds X(r)=x,\ea\right.\ee
with cost functional
$$\ba{ll}
\ns\ds
J(r,x;u(\cd))=\dbE\[\int_r^T\2n\(\1n\lan\cQ(s)X(s),X(s)\ran\1n+\1n
2\lan\cS(s)X(s),u(s)\ran
\1n+\1n\lan\cR(s)u(s),u(s)\ran\1n\)ds+\lan\cG X(T),X(T)\ran\].\ea$$
From \cite{Yong-Zhou 1999}, we know that under our conditions
(\ref{2.12}), the associated LQ problem admits a unique optimal
control $u^*(\cd)$ which has the following state feedback
representation:
$$u^*(s)=-\Th(s)X^*(s),\qq s\in[r,T],$$
with $\Th(\cd)$ given by
\bel{Th}\Th(s)=\big[\cR(s)+\cD(s)^TP(s)\cD(s)\big]^{-1}\big[\cB(s)^TP(s)+\cS(s)
+\cD(s)^TP(s)\cC(s)\big],\q s\in[t,T].\ee
and $X^*(\cd)$ is the solution to the following closed-loop system:
\bel{}\left\{\2n\ba{ll}
\ns\ds
dX^*(s)=\big[\cA(s)\1n-\1n\cB(s)\Th(s)\big]X^*(s)ds\1n+\1n\big[\cC(s)\1n-\1n\cD(s)
\Th(s)\big]X^*(s)dW(s),
\q s\in[r,T],\\
\ns\ds X^*(r)=x.\ea\right.\ee
Here, $P(\cd)$ is the unique solution to the Riccati equation
(\ref{Riccati0}). Moreover,
\bel{2.15}\ba{ll}
\ns\ds0\le\lan
P(r)x,x\ran=\inf_{u(\cd)\in\cU[r,T]}J(r,x;u(\cd))=J(r,x;u^*(\cd))\\
\ns\ds\q\le J(r,x;0)=\dbE\[\int_r^T\lan\cQ(s)X^0(s),X^0(s)\ran
ds+\lan\cG X^0(T),X^0(T)\ran\],\ea\ee
where $X^0(\cd)$ is the solution to (\ref{2.14}) with $u(\cd)=0$,
i.e.,
\bel{2.8}\left\{\2n\ba{ll}
\ns\ds dX^0(s)=\cA(s)X^0(s)ds+\cC(s)X^0(s)dW(s),\qq s\in[r,T],\\
\ns\ds X^0(r)=x.\ea\right.\ee
On the other hand, since $\Pi(\cd)$ is the solution to
(\ref{Lyapunov1}), applying It\^o's formula to
$\lan\Pi(\cd)X^0(\cd),X^0(\cd)\ran$, we have
\bel{2.17}\ba{ll}
\ns\ds\lan\Pi(r)x,x\ran=\dbE\[\lan\Pi(T)X^0(T),X^0(T)\ran\\
\ns\ds\qq\qq\qq-\1n\int_r^T\3n\lan\big[\dot\Pi(s)\1n+\1n
\Pi(s)\cA(s)\1n+\1n\cA(s)^T\1n\Pi(s)\1n
+\1n\cC(s)^T\1n\Pi(s)\cC(s)\big]X^0(s),X^0(s)\ran ds\]\\
\ns\ds\qq\qq\q=\dbE\[\lan\cG
X^0(T),X^0(T)\ran+\int_r^T\lan\cQ(s)X^0(s),X^0(s)\ran
ds\]\\
\ns\ds\qq\qq\q=J(r,x;0)\ge\lan P(r)x,x\ran.\ea\ee
Next, applying It\^o's formula to $|X^0(\cd)|^2$, one has
$$\ba{ll}
\ns\ds\dbE|X^0(s)|^2=|x|^2+\dbE\int_r^s\lan\big[\cA(\t)+\cA(\t)^T+\cC(\t)^T\cC(\t)
\big]X^0(\t),X^0(\t)\ran d\t\\
\ns\ds\qq\qq~\le|x|^2+\|\cA(\cd)+\cA(\cd)^T+\cC(\cd)^T\cC(\cd)\|_\infty\int_r^s\dbE|X^0(r)|^2dr.\ea$$
Then, by Gronwall's inequality, together with (\ref{2.15}) and
(\ref{2.17}), we obtain
$$0\le P(r)\le\Pi(r)\le K_0I,\qq r\in[t,T].$$
This completes the proof. \endpf

\ms

Note that with the function $\Th(\cd)$ defined by (\ref{Th}),
Riccati equation (\ref{Riccati0}) can be written as
\bel{}\left\{\2n\ba{ll}
\ns\ds\dot
P(s)+P(s)\big[\cA(s)-\cB(s)\Th(s)\big]+\big[\cA(s)-\cB(s)\Th(s)\big]^TP(s)+\cQ(s)\\
\ns\ds\qq+\big[\cC(s)-\cD(s)\Th(s)\big]^TP(s)\big[\cC(s)-\cD(s)\Th(s)\big]
=0,\qq s\in[t,T],\\
\ns\ds P(T)=\cG.\ea\right.\ee
If $\Th(\cd)$ is given, the above can be regarded as a Lyapunov
equation for $P(\cd)$. The following lemma will play an interesting
role in our later investigations.

\ms

\bf Lemma 2.4. \sl Let
\bel{}\left\{\2n\ba{ll}
\ns\ds\cA(\cd),\cC(\cd)\in
L^\infty(t,T;\dbR^n),\q\cB(\cd),\cD(\cd)\in
L^\infty(t,T;\dbR^{n\times m}),\\
\ns\ds\cQ(\cd),\wt\cQ(\cd),\bar\cQ(\cd)\in
L^\infty(t,T;\dbS^n),\q\cG,\bar\cG\in\dbS^n.\ea\right.\ee
Then there are $\wt\G(\cd)$, $\G(\cd)$ and $\bar\G(\cd)$ uniquely
solve the following three Lyapunov equations:
\bel{}\left\{\2n\ba{ll}
\ns\ds\dot{\wt\G}(s)+\wt\G(s)\cA(s)+\cA(s)^T\wt\G(s)+\cC(s)^T\wt\G(s)\cC(s)
+\cQ(s)=0,\qq s\in[t,T],\\
\ns\ds\wt\G(T)=\cG,\ea\right.\ee
\bel{}\left\{\2n\ba{ll}
\ns\ds\dot\G(s)+\G(s)[\cA(s)+\bar\cA(s)]+[\cA(s)+\bar\cA(s)]^T
\G(s)\\
\ns\ds\qq+\big[\cC(s)+\bar\cC(s)\big]^T\wt\G(s)\big[\cC(s)+\bar\cC(s)\big]
+\cQ(s)+\wt\cQ(s)=0,\qq s\in[t,T],\\
\ns\ds\G(T)=\cG,\ea\right.\ee
and
\bel{}\left\{\2n\ba{ll}
\ns\ds\dot{\bar\G}(s)+\bar\G(s)[\cA(s)+\bar\cA(s)]+[\cA(s)+\bar\cA(s)]^T
\bar\G(s)+\bar\cQ(s)=0,\qq
s\in[t,T],\\
\ns\ds\bar\G(T)=\bar\cG.\ea\right.\ee
Let $(t,x)\in\sD$ and $X(\cd)\equiv X(\cd\,;t,x)$ be the solution to
$(\ref{MF-SDE1})$ with $b(\cd)=\si(\cd)=0$. Then for any
$\t\in[0,t]$,
\bel{2.4}\ba{ll}
\ns\ds\dbE_\t\[\int_t^T\2n\(\1n\lan\cQ(s)X(s),X(s)\ran\1n+\1n\lan\wt\cQ(s)
\dbE_t[X(s)],\dbE_t[X(s)]\ran\1n+\1n\lan\bar\cQ(s)\dbE_\t[X(s)],
\dbE_\t[X(s)]\ran\1n\)ds\\
\ns\ds\qq\qq+\lan\cG
X(T),X(T)\ran+\lan\bar\cG\dbE_\t[X(T)],\dbE_\t[X(T)]\ran\]\\
\ns\ds\q=\dbE_\t\big[\lan\G(t)x,x\ran\big]+\lan\bar\G(t)\dbE_\t[x],
\dbE_\t[x]\ran.\ea\ee
Further, if
\bel{2.16}\left\{\2n\ba{ll}
\ns\ds\cQ(s),~\cQ(s)+\wt\cQ(s),~\cQ(s)+\wt\cQ(s)+\bar\cQ(s)\ge0,\qq
s\in[0,T],\\
\ns\ds\cG,~\cG+\bar\cG\ge0,\ea\right.\ee
then
\bel{2.28}\wt\G(s),~\G(s),~\G(s)+\bar\G(s)\ge0,\qq s\in[0,T].\ee

\ms

\rm

Note that in the left hand side of (\ref{2.4}), the conditional
expectation terms $\dbE_t[X(s)]$, $\dbE_\t[X(s)]$ and
$\dbE_\t[X(T)]$ (with $0\le\t\le t$) appear quadratically (not
linearly). It turns out that one can represent them in terms of $x$
and $\dbE_\t[x]$ by using $\G(\cd)$ and $\bar\G(\cd)$ (which also
involves $\wt\G(\cd)$ indirectly).

\ms

\it Proof. \rm First of all, since Lyapunov equations are linear,
under our conditions, $\wt\G(\cd)$, $\G(\cd)$, and $\bar\G(\cd)$
uniquely exist. Note that for $0\le\t\le t<T$, and $x\in\sX_t$, let
$X(\cd)=X(\cd\,;t,x)$ be the solution (\ref{MF-SDE1}) with
$b(\cd)=\si(\cd)=0$. We have
$$\left\{\2n\ba{ll}
\ns\ds
d\big(\dbE_\t[X(s)]\big)=\big[\cA(s)+\bar\cA(s)\big]\dbE_\t[X(s)]ds,\qq s\in[t,T],\\
\ns\ds\dbE_\t[X(t)]=\dbE_\t[x],\ea\right.$$
$$\left\{\2n\ba{ll}
\ns\ds
d\big(\dbE_t[X(s)]-\dbE_\t[X(s)]\big)=\big[\cA(s)+\bar\cA(s)\big]\big(\dbE_t[X(s)]
-\dbE_\t[X(s)]\big)ds,\q s\in[t,T],\\
\ns\ds\dbE_t[X(t)]-\dbE_\t[X(t)]=x-\dbE_\t[x],\ea\right.$$
and
$$\left\{\2n\ba{ll}
\ns\ds
d\big(X(s)-\dbE_t[X(s)]\big)=\cA(s)\big(X(s)-\dbE_t[X(s)]\big)ds+\[\cC(s)X(s)
+\bar\cC(s)\dbE_t[X(s)]\]dW(s),\q s\in[t,T],\\
\ns\ds X(t)-\dbE_t[X(t)]=0.\ea\right.$$
Denote $\h\G(\cd)=\G(\cd)+\bar\G(\cd)$. Then $\h\G(\cd)$ satisfies
the following:
\bel{}\left\{\2n\ba{ll}
\ns\ds\dot{\h\G}(s)+\h\G(s)[\cA(s)+\bar\cA(s)]+[\cA(s)+\bar\cA(s)]^T
\h\G(s)\\
\ns\ds\qq+[\cC(s)+\bar\cC(s)]^T\wt\G(s)[\cC(s)+\cC(s)]+\cQ(s)+\wt\cQ(s)+\bar\cQ(s)=0,\qq
s\in[t,T],\\
\ns\ds\h\G(T)=\cG+\bar\cG,\ea\right.\ee
By It\^o's formula, we have
$$\ba{ll}
\ns\ds\dbE_\t\[\lan\cG
X(T),X(T)\ran+\lan\bar\cG\dbE_\t[X(T)],\dbE_\t[X(T)]\ran\]\\
\ns\ds=\dbE_\t\[\lan\cG
\big(X(T)-\dbE_t[X(T)]\big),X(T)-\dbE_t[X(T)]\ran\\
\ns\ds\qq+\lan\cG
\big(\dbE_t[X(T)]-\dbE_\t[X(T)]\big),\dbE_t[X(T)]-\dbE_\t[X(T)]\ran
+\lan(\cG+\bar\cG)\dbE_\t[X(T)],\dbE_\t[X(T)]\ran\]\\
\ns\ds=\dbE_\t\[\lan\wt\G(T)
\big(X(T)-\dbE_t[X(T)]\big),X(T)-\dbE_t[X(T)]\ran\\
\ns\ds\qq+\lan\G(T)
\big(\dbE_t[X(T)]-\dbE_\t[X(T)]\big),\dbE_t[X(T)]-\dbE_\t[X(T)]\ran
+\lan\h\G(T)\dbE_\t[X(T)],\dbE_\t[X(T)]\ran\]\\
\ns\ds=\dbE_\t\[\lan\G(t)\big(x-\dbE_\t[x]\big),x-\dbE_\t[x]\ran
+\lan\h\G(t)\dbE_\t[x],
\dbE_\t[x]\ran\]\\
\ns\ds\qq+\dbE_\t\int_t^T\1n\[\lan(\dot{\wt\G}\1n+\1n\wt\G\cA\1n+\1n\cA^T\wt\G)\big(X\1n-\1n\dbE_t[X]\big),
X\1n-\1n\dbE_t[X]\ran\1n+\1n\lan\wt\G\big(\cC
X\1n+\1n\bar\cC\dbE_t[X]\big),\cC
X\1n+\1n\bar\cC\dbE_t[X]\ran\\
\ns\ds\qq+\lan\big[\dot\G+\G(\cA+\bar\cA)+(\cA+\bar\cA)^T\G\big]
\big(\dbE_t[X]-\dbE_\t[X]\big),\dbE_t[X]-\dbE_\t[X]\ran\\
\ns\ds\qq+\lan\big\{\dot{\h\G}+\h\G(\cA+\bar\cA)
+(\cA+\bar\cA)^T\h\G\big\}\dbE_\t[X],\dbE_\t[X]\ran\]ds\\
\ns\ds=\1n\dbE_\t\1n\[\1n\lan\G(t)x,x\ran\1n+\1n\lan\bar\G(t)\dbE_\t[x],
\dbE_\t[x]\ran\1n\]\1n+\1n\dbE_\t\3n\int_t^T\3n\[\1n\lan(\dot{\wt\G}\1n+\1n
\wt\G\cA\1n+\1n\cA^T\wt\G\1n+\1n\cC^T\wt\G\cC)
\big(X\1n-\1n\dbE_t[X]\big),\1n X\1n-\1n\dbE_t[X]\ran\\
\ns\ds\qq+\lan\big[\dot\G+\G(\cA+\bar\cA)+(\cA+\bar\cA)^T\G
+(\cC+\bar\cC)^T\wt\G(\cC+\bar\cC)\big]
\big(\dbE_t[X]-\dbE_\t[X]\big),\dbE_t[X]-\dbE_\t[X]\ran\\
\ns\ds\qq+\lan\big\{\dot{\h\G}+\h\G(\cA+\bar\cA)
+(\cA+\bar\cA)^T\h\G+(\cC+\bar\cC)^T\wt\G(\cC+\bar\cC)\big\}\dbE_\t[X],\1n
\dbE_\t[X]\ran\1n\]ds.\ea$$
On the other hand,
$$\ba{ll}
\ns\ds\dbE_\t\int_t^T\[\lan\cQ
X,X\ran+\lan\wt\cQ\dbE_t[X],\dbE_t[X]\ran
+\lan\bar\cQ\dbE_\t[X],\dbE_\t[X]\ran\]ds\\
\ns\ds\q=\dbE_\t\int_t^T\[\lan\cQ\big(X-\dbE_t[X]\big),X-\dbE_t[X]\ran
+\lan(\cQ+\wt\cQ)\big(\dbE_t[X]-\dbE_\t[X]\big),\dbE_t[X]-\dbE_\t[X]\ran\\
\ns\ds\qq\qq+\lan(\cQ+\wt\cQ+\bar\cQ)\dbE_\t[X],\dbE_\t[X]\ran\]ds.\ea$$
Hence,
$$\ba{ll}
\ns\ds\dbE_\t\[\int_t^T\(\lan\cQ X,X\ran+\lan\wt
Q\dbE_t[X],\dbE_t[X]\ran+\lan\bar\cQ\dbE_\t[X],\dbE_\t[X]\ran\)ds\\
\ns\ds\qq\qq+\lan\cG X(T),X(T)\ran+\lan\bar\cG\dbE_\t[X(T)],\dbE_\t[X(T)]\ran\]\\
\ns\ds\q=\dbE_\t\[\lan\G(t)x,x\ran\]+\lan\bar\G(t)\dbE_\t[x],
\dbE_\t[x]\ran\\
\ns\ds\qq+\dbE_\t\int_t^T\[\lan(\dot{\wt\G}+\wt\G\cA+\cA^T\wt\G+\cC^T\wt\G\cC
+\cQ)\big(X-\dbE_\t[X]\big),
X-\dbE_\t[X]\ran\\
\ns\ds\qq+\lan\2n\big[\dot\G\1n+\1n\G(\cA\1n+\1n\bar\cA)\1n+\1n(\cA\1n+\1n\bar\cA)^T\1n\G
\1n+\1n(\cC\1n+\1n\bar\cC)^T\wt\G(\cC\1n+\1n\bar\cC)\1n+\1n\cQ\1n+\1n\wt\cQ\big]
\big(\dbE_t[X]\1n-\1n\dbE_\t[X]\big),\dbE_t[X]\1n-\1n\dbE_\t[X]\ran\\
\ns\ds\qq+\lan\big\{\dot{\h\G}+\h\G(\cA+\bar\cA)
+(\cA+\bar\cA)^T\h\G+(\cC+\bar\cC)^T\wt\G(\cC+\bar\cC)+\cQ+\wt\cQ+\bar\cQ\big\}\dbE_\t[X],\dbE_\t[X]\ran\]ds\\
\ns\ds\q=\dbE_\t\big[\lan\G(t)x,x\ran\big]+\lan\bar\G(t)\dbE_\t[x],
\dbE_\t[x]\ran.\ea$$
This proves (\ref{2.4}). Finally, if (\ref{2.16}) hold, by
Proposition 2.3, we obtain (\ref{2.17}). \endpf

\ms

\section{Linear-Quadratic Optimal Control Problem for MF-SDEs}

We now consider the following controlled linear MF-SDE:
\bel{state}\left\{\2n\ba{ll}
\ns\ds dX(s)=\Big\{A(s)X(s)+\bar A(s)\dbE_t[X(s)]+B(s)u(s)+\bar
B(s)\dbE_t[u(s)]\Big\}ds\\
\ns\ds\qq\qq+\Big\{C(s)X(s)\1n+\1n\bar C(s)\dbE_t[X(s)]\1n+\1n
D(s)u(s)\1n+\1n\bar D(s)
\dbE_t[u(s)]\Big\}dW(s),\q s\in[t,T],\\
\ns\ds X(t)=x\in\sX_t,\ea\right.\ee
Note that for any $(t,x)\in\sD$ and $u(\cd)\in\sU[t,T]$, the
corresponding state process $X(\cd)\1n=\2n X(\cd\,;t,x,u(\cd))$
depends on $(t,x,u(\cd))$. The cost functional is as follows:
\bel{cost}\ba{ll}
\ns\ds J(t,x;u(\cd))=\dbE_t\Big\{\int_t^T\[\lan
Q(s,t)X(s),X(s)\ran+\lan\bar
Q(s,t)\dbE_t[X(s)],\dbE_t[X(s)]\ran\\
\ns\ds\qq\qq\qq\qq\q+\lan R(s,t)u(s),u(s)\ran+\lan\bar
R(s,t)\dbE_t[u(s)],\dbE_t[u(s)]\ran\]ds\\
\ns\ds\qq\qq\qq\qq\q+\lan G(t)X(T),X(T)\ran+\lan\bar
G(t)\dbE_t[X(T)],\dbE_t[X(T)]\ran\Big\}.\ea\ee
Let us introduce the following hypotheses:

\ms

{\bf(H1)} The following hold:
\bel{H1}\left\{\2n\ba{ll}
\ns\ds A(\cd),\bar A(\cd),C(\cd),\bar C(\cd)\in
C([0,T];\dbR^{n\times n}),\\
\ns\ds B(\cd),\bar B(\cd),D(\cd),\bar D(\cd)\in
C([0,T];\dbR^{n\times m}).\ea\right.\ee

\ms

{\bf(H2)} The following hold:
\bel{H2-1}\left\{\2n\ba{ll}
\ns\ds Q(\cd\,,\cd),\bar Q(\cd\,,\cd)\in C([0,T]^2;\dbS^n),\q
R(\cd\,,\cd),\bar R(\cd\,,\cd)\in C([0,T]^2;\dbS^m),\\
\ns\ds G(\cd),\bar G(\cd)\in C([0,T];\dbS^n),\ea\right.\ee
and for some $\d>0$,
\bel{H2-2}\left\{\2n\ba{ll}
\ns\ds Q(s,t),\;Q(s,t)+\bar Q(s,t)\ge0,\q R(s,t),\;R(s,t)+\bar
R(s,t)\ge\d I,\qq0\le t\le s\le T,\\
\ns\ds G(t),\;G(t)+\bar G(t)\ge0,\qq0\le t\le T.\ea\right.\ee

{\bf(H3)} The following monotonicity conditions are satisfied:
\bel{}\left\{\2n\ba{ll}
\ns\ds Q(s,t)\le Q(s,\t),\q Q(s,t)+\bar Q(s,t)\le Q(s,\t)+\bar
Q(s,\t),\\
\ns\ds R(s,t)\le R(s,\t),\q R(s,t)+\bar R(s,t)\le R(s,\t)+\bar
R(s,\t),\\
\ns\ds G(t)\le G(\t),\q G(t)+\bar G(t)\le G(\t)+\bar
G(\t),\ea\right.\qq0\le t\le\t\le s\le T.\ee \ms

It is clear that under (H1)--(H2), for any $(t,x)\in\sD$ and
$u(\cd)\in\sU[t,T]$, state equation (\ref{state}) admits a unique
solution $X(\cd)\equiv X(\cd\,;t,x,u(\cd))$, and the cost functional
$J(t,x;u(\cd))$ is well-defined. Then we can state the following
problem.

\ms

\bf Problem (MF-LQ). \rm For any $(t,x)\in\sD$, find a
$u^*(\cd)\in\sU[t,T]$ such that
\bel{J}J(t,x;u^*(\cd))=\inf_{u(\cd)\in\sU[t,T]}J(t,x;u(\cd))\equiv
V(t,x).\ee
For given $(t,x)\in\sD$, any $u^*(\cd)\in\sU[t,T]$ satisfying the
above is called a {\it pre-commitment optimal control} for Problem
(MF-LQ) at $(t,x)$. The corresponding $X^*(\cd)$ and
$(X^*(\cd),u^*(\cd))$ are called {\it pre-commitment optimal state
process} and {\it pre-commitment optimal pair} of Problem (MF-LQ),
respectively, and $V(\cd\,,\cd)$ is called the {\it pre-commitment
value function}.

\ms

In the case that
\bel{A=0}\bar A(s)=\bar C(s)=0,\qq\bar B(s)=\bar D(s)=0,\qq
s\in[0,T],\ee
and
\bel{Q=0}\bar Q(s,t)=0,\qq\bar R(s,t)=0,\qq\bar G(t)=0,\qq0\le t\le
s\le T,\ee
Problem (MF-LQ) is reduced to the problem studied in \cite{Yong
2011}, which is an LQ problem with purely general (possibly
non-exponential) discounting. In addition, if
\bel{Q=Q}Q(s,t)=Q(s),\q R(s,t)=R(s),\q G(t)=G,\qq0\le t\le s\le
T,\ee
then Problem (MF-LQ) is further reduced to a classical LQ problem
for SDEs with deterministic coefficients (see \cite{Yong-Zhou 1999}
for details). In the case that (\ref{A=0}) and (\ref{Q=Q}) hold, and
\bel{barQ=barQ}\bar Q(s,t)=0,\q\bar R(s,t)=0,\q\bar G(t)=\bar
G,\qq0\le t\le T,\ee
Problem (MF-LQ) is a case of those studied in \cite{Hu-Jin-Zhou
2011} where the problem was investigated by means of variational
method together with forward-backward SDEs and the time-consistent
solution studied is of open-loop type (see the next section for
further details).

\ms

We also point out that if
$$\left\{\2n\ba{ll}
\ns\ds Q(s,t)=Q(s)e^{-\l(s-t)},\q\bar Q(s,t)=e^{-\l(s-t)},\\
\ns\ds R(s,t)=R(s)e^{-\l(s-t)},\q\bar R(s,t)=e^{-\l(s-t)},\\
\ns\ds G(t)=Ge^{-\l(T-t)},\qq\bar G(t)=\bar
Ge^{-\l(T-t)},\ea\right.\qq0\le t\le s\le T,$$
with $\l>0$, and with
$$\left\{\2n\ba{ll}
\ns\ds Q(s),Q(s)+\bar Q(s)\ge0,\q R(s),R(s)+\bar R(s)\ge0,\q
s\in[0,T],\\
\ns\ds G,G+\bar G\ge0,\ea\right.$$
then monotonicity condition (H3) holds. Hypothesis (H3) will be used
later.

\ms

In what follows, we will denote
\bel{hat}\left\{\2n\ba{ll}
\ns\ds\h A(s)=A(s)+\bar A(s),\qq\h B(s)=B(s)+\bar B(s),\\
\ns\ds\h C(s)=C(s)+\bar C(s),\qq\h D(s)=D(s)+\bar
D(s),\\
\ns\ds\h Q(s,t)=Q(s,t)+\bar Q(s,t),\q\h R(s,t)=R(s,t)+\bar
R(s,t),\\
\ns\ds\h G(t)=G(t)+\h G(t).\ea\right.\qq0\le t\le s\le T.\ee

\ms

We now recall a result from \cite{Yong 2012c} for the pre-commitment
solution of Problem (MF-LQ).

\ms

\bf Proposition 3.1. \sl Let {\rm(H1)}--{\rm(H2)} hold. Then for any
fixed $t\in[0,T)$, the following Riccati equation system admits a
unique solution $(P(\cd),\h P(\cd))\in C^1([t,T];\dbS^n)^2$
(suppressing $s$):
\bel{Riccati P}\left\{\2n\ba{ll}
\ns\ds\dot P+PA+A^TP+C^TPC+Q(t)-(PB+C^TPD)\big[R(t)+D^TPD\big]^{-1}
(B^TP+D^TPC)=0,\\
\dot{\h P}+\h P\h A+\h A^T\h P+\h C^TP\h C+\h Q(t)-(\h P\h B+\h C^T
P\h D)\big[\h R(t)+\h D^TP\h D\big]^{-1}(\h B^T\h P+\h D^TP\h C)=0,\\
\ns\ds\qq\qq\qq\qq\qq\qq\qq\qq\qq\qq\qq\qq\qq s\in[t,T],\\
\ns\ds P(T)=G(t),\q\h P(T)=\h G(t).\ea\right.\ee
Further, let $x\in\sX_t$ and $X^*(\cd)\equiv X^*(\cd\,;t,x)$ be the
solution to the following closed-loop system:
\bel{closed-loop}\left\{\2n\ba{ll}
\ns\ds dX^*\1n(s)\1n=\2n\Big\{\1n\big[A(s)\1n-\1n
B(s)\Th(s)\big]X^*\1n(s)\1n+\1n\[\bar A(s)\1n+\1n
B(s)\big[\Th(s)\1n-\1n\h\Th(s)\big]\1n-\1n\bar
B(s)\h\Th(s)\]\dbE_t[X^*\1n(s)]\1n\Big\}ds\\
\ns\ds\;+\Big\{\1n\big[C(s)\2n-\2n
D(s)\Th(s)\big]X^*(s)\1n+\1n\[\bar C(s)\2n+\2n
D(s)\1n\big[\Th(s)\1n-\1n\h\Th(s)\big]\1n-\1n\bar
D(s)\h\Th(s)\]\dbE_t[X^*(s)]\1n\Big\}dW(s),\q s\1n\in\1n[t,\1n T],\\
\ns\ds X^*(t)=x,\ea\right.\ee
with
\bel{Th1}\left\{\2n\ba{ll}
\ns\ds\Th(s)=\big[R(s,t)+D(s)^TP(s)D(s)\big]^{-1}\big[B(s)^TP(s)+D(s)^TP(s)C(s)\big],\\
\ns\ds\h\Th(s)=\big[\h R(s,t)+\h D(s)^TP(s)\h D(s)\big]^{-1}\big[\h
B(s)^T\h P(s)+\h D(s)^TP(s)\h C(s)\big],\ea\right.\qq s\in[t,T],\ee
and define $u^*(\cd)$ as follows:
\bel{feedback}
u^*(s)=-\Th(s)X^*(s)+\big[\Th(s)-\h\Th(s)\big]\dbE_t[X^*(s)],\qq
s\in[t,T].\ee
Then $(X^*(\cd),u^*(\cd))$ is the pre-commitment optimal pair of
Problem {\rm(MF-LQ)} at $(t,x)$, and
\bel{inf}\ba{ll}
\ns\ds V(t,x)=\inf_{u(\cd)\in\cU[t,T]}J(t,x;u(\cd))
=J(t,x;u^*(\cd))=\lan\h P(t)x,x\ran,\qq\forall x\in\sX_t.\ea\ee

\ms

\rm

We point out that system (\ref{Riccati P}) is decoupled. With
$\Th(\cd)$ and $\h\Th(\cd)$ given by (\ref{Th1}), we can write
(\ref{Riccati P}) as follows:
\bel{Riccati P*}\left\{\2n\ba{ll}
\ns\ds\dot
P(s)+P(s)\[A(s)-B(s)\Th(s)\]+\[A(s)-B(s)\Th(s)\]^TP(s)+Q(s,t)\\
\ns\ds\q+\[C(s)-D(s)\Th(s)\]^TP(s)\[C(s)-D(s)\Th(s)\]+\Th(s)^TR(s,t)\Th(s)=0,\\
\ns\ds\dot{\h P}(s)+\h P(s)\[\h A(s)-\h B(s)\h\Th(s)\]+\[\h A(s)-\h
B(s)\h\Th(s)\]^T\h P(s)+\h Q(s)\\
\ns\ds\q+\[\h C(s)-\h D(s)\h\Th(s)\]^TP(s)\[\h C(s)-\h
D(s)\h\Th(s)\]+\h\Th(s)^T\h R(s,t)\h\Th(s)=0,\q s\in[t,T],\\
\ns\ds P(T)=G(t),\qq\h P(T)=\h G(t).\ea\right.\ee
Next, we note that in the case
$$\bar A(\cd)=\bar C(\cd)=0,\qq\bar B(\cd)=\bar D(\cd)=0,$$
the conditional expectation terms are all absent in the state
equation, and
\bel{Si0}\left\{\2n\ba{ll}
\ns\ds\Th(s)=\big[R(s,t)+D(s)^TP(s)D(s)\big]^{-1}\big[B(s)^TP(s)+D(s)^TP(s)C(s)\big],\\
\ns\ds\h\Th(s)=\big[\h R(s,t)+D(s)^TP(s)D(s)\big]^{-1}\big[B(s)^T\h
P(s)+D(s)^TP(s)C(s)\big],\ea\right.\qq s\in[t,T].\ee
In this case, the equation for $\h P(\cd)$ becomes
\bel{Riccati Pi*}\left\{\2n\ba{ll}
\dot{\h P}(s)+\h P(s)A(s)+A(s)^T\h P+C(s)^TP(s)C(s)+\h Q(t,s)\\
\ns\ds\q-\big[\h P(s)B(s)+C(s)^TP(s)D(s)\big]\big[\h
R(s,t)+D(s)^TP(s)D(s)\big]^{-1}\\
\ns\ds\qq\qq\qq\qq\cd\big[B(s)^T\h P(s)+D(s)^TP(s)C(s)\big]=0,\qq s\in[t,T],\\
\ns\ds\h P(T)=\h G(t),\ea\right.\ee
We see that as long as
$$\bar Q(\cd\,,\cd)=0,\qq\bar R(\cd\,,\cd)=0,\qq\bar G(\cd)=0$$
are not true,
$$\Th(s)=\h\Th(s),\qq s\in[t,T]$$
is not true in general. Hence, the term $\dbE_t[X^*(\cd)]$ will
present in the state feedback representation of $u^*(\cd)$ (see
(\ref{feedback})), and the closed-loop system reads
\bel{closed-loop0a}\left\{\2n\ba{ll}
\ns\ds dX^*\1n(s)\1n=\2n\Big\{\1n\big[A(s)\1n-\1n
B(s)\Th(s)\big]X^*\1n(s)\1n+\1n B(s)\big[\Th(s)\1n-\1n\h\Th(s)\big]\dbE_t[X^*\1n(s)]\Big\}ds\\
\ns\ds\qq\qq+\Big\{\big[C(s)\1n-\1n D(s)\Th(s)\big]X^*(s)\1n+\1n
D(s)\big[\Th(s)\1n-\1n\h\Th(s)\big]\dbE_t[X^*(s)]\Big\}dW(s),\q s\1n\in\1n[t,T],\\
\ns\ds X^*(t)=x,\ea\right.\ee
which is a linear MF-SDE. The above shows that even if the original
state equation does not contain the conditional expectation terms,
as long as some conditional expectation terms appear in the cost
functional, the closed-loop system for the pre-committed solution of
the problem will still be an MF-SDE. This is why we start with a
controlled MF-SDE.

\ms

We now present the following result which is relevant to the Riccati
equation system (\ref{Riccati P}).

\ms

\bf Proposition 3.2. \sl Let {\rm(H1)}--{\rm(H2)} hold. For a fixed
$t\in[0,T)$, let $(P(\cd),\h P(\cd))$ be the solution to the Riccati
equation system $(\ref{Riccati P})$, let $\Pi(\cd)$ and $\h\Pi(\cd)$
be the solutions to the following Lyapunov equations:
\bel{Lyapunov P0}\left\{\2n\ba{ll}
\ns\ds\dot
\Pi(s)+\Pi(s)A(s)+A(s)^T\Pi(s)+C(s)^T\Pi(s)C(s)+Q(s,t)=0,\qq s\in[t,T],\\
\ns\ds\Pi(T)=G(t),\ea\right.\ee
and
\bel{Lyapunov Pi0}\left\{\2n\ba{ll}
\ns\ds\dot{\h\Pi}(s)+\h\Pi(s)\h A(s)+\h A(s)^T\h\Pi(s)+\h
C(s)^T\Pi(s)\h C(s)
+\h Q(s,t)=0,\qq s\in[t,T],\\
\ns\ds\h\Pi(T)=\h G(t).\ea\right.\ee
Then
\bel{3.28}0\le P(s)\le\Pi(s),\qq s\in[t,T],\ee
and
\bel{3.29}0\le\h P(s)\le\h\Pi(s),\qq s\in[t,T].\ee

\ms

\it Proof. \rm Fix $t\in[0,T)$. Applying Proposition 2.3 with
$$\left\{\2n\ba{ll}
\ns\ds\cA(s)=A(s),\q\cB(s)=B(s),\q\cC(s)=C(s),\q\cD(s)=D(s),\q t\le s\le T,\\
\ns\ds\cQ(s)=Q(s,t),\q\cS(s)=0,\q\cR(s)=R(s,t),\q\cG=G(t),\q t\le
s\le T,\ea\right.$$
we obtain (\ref{3.28}). Next, we let
$$\left\{\2n\ba{ll}
\ns\ds\cA(s)\1n=\1n\h A(s),\q\cB(s)\1n=\1n\h B(s),\q\cC(s)\1n=\1n \h
C(s),\q\cD(s)\1n=\1n\h D(s),\\
\ns\ds\cQ(s)=\h Q(s,t)+\cC(s)^TP(s)\cC(s),\q\cS(s)=\cD(s)^TP(s)\cC(s),\\
\ns\ds\cR_0(s)=\h R(s,t),\q\cR(s)=\cR_0(s)+\cD(s)^TP(s)\cD(s),\q
\cG=\h G(t).\ea\right.$$
Then the equation for $\h P(\cd)$ can be written as
\bel{}\left\{\2n\ba{ll}
\dot{\h P}(s)+\h P(s)\cA(s)+\cA(s)^T\h P(s)+\cQ(s)\\
\ns\ds\q-\big[\h P(s)\cB(s)+\cS(s)^T\big]\cR(s)^{-1}
\big[\cB(s)\h P(s)\1n+\1n\cS(s)\big]\1n=\1n0,\qq s\in[t,T],\\
\ns\ds\h P(T)=\cG,\ea\right.\ee
which is the Riccati equation for the deterministic LQ problem with
the state equation
$$\left\{\2n\ba{ll}
\ns\ds\dot X(s)=\cA(s)X(s)+\cB(s)u(s),\qq s\in[r,T],\\
\ns\ds X(r)=x,\ea\right.$$
and the cost functional
$$\ba{ll}
\ns\ds
J(r,x;u(\cd))=\int_r^T\(\lan\cQ(s)X(s),X(s)\ran+2\lan\cS(s)X(s),u(s)\ran
+\lan\cR(s)u(s),u(s)\ran\)ds+\lan\cG X(T),X(T)\ran.\ea$$
We now look at the following: (suppressing $s$)
$$\ba{ll}
\ns\ds\cQ-\cS^T\cR^{-1}\cS\ge\cC^TP\cC-\cC^TP\cD(\cR_0+\cD^TP\cD)^{-1}\cD^TP\cC\\
\ns\ds=\cC^TP^{1\over2}
\[I-P^{1\over2}\cD(\cR_0+\cD^TP\cD)^{-1}\cD^TP^{1\over2}\]P^{1\over2}\cC\\
\ns\ds=\cC^TP^{1\over2}\[I\1n-\1n
P^{1\over2}\cD\cR_0^{-{1\over2}}\(I\1n+\1n\cR_0^{-{1\over2}}\cD^T
P^{1\over2}P^{1\over2}\cD\cR_0^{-{1\over2}}\)^{-1}\cR_0^{-{1\over2}}\cD^T
P^{1\over2}\]P^{1\over2}\cC\\
\ns\ds\equiv\cC^T
P^{1\over2}\[I\1n-\1n\L(I\1n+\1n\L^T\L)^{-1}\L^T\]P^{1\over2}\cC
=\cC^T\1n P^{1\over2}\1n(I\1n+\1n\L\L^T\1n)^{-1}\1n
P^{1\over2}\cC\ge0,\ea$$
In the above, we have denoted
$$\L=P^{1\over2}\cD\cR_0^{-{1\over2}},$$
and used the fact
$$I-\L(I+\L^T\L)^{-1}\L^T=(I+\L\L^T)^{-1}.$$
Then by Proposition 2.3 with $\cC(\cd)=0$ and $\cD(\cd)=0$, we have
$$\h P(s)\le\bar P^0(s),\qq s\in[t,T],$$
where $\bar P^0(\cd)$ is the solution to the following Lyapunov
equation:
\bel{}\left\{\2n\ba{ll}
\dot{\bar P}^0(s)+\bar P^0(s)\cA(s)+\cA(s)^T\bar P^0(s)+\cQ(s)=0,\qq s\in[t,T],\\
\ns\ds\bar P^0(T)=\cG,\ea\right.\ee
Next, by the proved (\ref{3.28}), we see that
$$\cQ(s)\equiv\h Q(s,t)+\cC(s)^TP(s)\cC(s)\le\h
Q(s,t)+\cC(s)^T\Pi(s)\cC(s).$$
Thus, by Proposition 2.3 again, we obtain
$$\bar P^0(s)\le P^0(s),$$
with $P^0(\cd)$ being the solution to the following Lyapunov
equation: (note $\cA(\cd)=\h A(\cd)$, $\cC(\cd)=\h C(\cd)$, and
$\cG=\h G(t)$)
$$\left\{\2n\ba{ll}
\dot P^0(s)+P^0(s)\h A(s)+\h A(s)^TP^0(s)+\h C(s)^T\Pi(s)\h C(s)+\h Q(s,t)=0,\qq s\in[t,T],\\
\ns\ds P^0(T)=\h G(t).\ea\right.$$
Comparing with (\ref{Lyapunov Pi0}), by uniqueness, we obtain
$$P^0(\cd)=\h\Pi(\cd).$$
This completes the proof. \endpf

\ms

\section{Time-Consistent Equilibrium Solutions.}

From Examples 1.1 and 1.2, we see that in general, Problem (MF-LQ)
is time-inconsistent. Our goal is to find time-consistent
(equilibrium) solutions to Problem (MF-LQ).

\ms

\subsection{Open-loop equilibrium control.}

In this subsection, we consider the case that
\bel{=0}\bar A(\cd)=\bar C(\cd)=0,\qq\bar B(\cd)=\bar D(\cd)=0,\ee
Thus, the state equation reads as (\ref{1.1}). However, we still
consider cost functional (\ref{cost}) in which all the conditional
expectation terms appear, and all the weighting matrices depend on
$(s,t)$. Note that in \cite{Hu-Jin-Zhou 2011}, the state equation is
similar to ours (with additional nonhomogeneous terms, which do not
bring additional essential difficulties) and the cost functional has
the only conditional expectation appeared at the terminal cost, and
all the weighting matrices are independent of $t$ (the initial
time). Therefore, in some sense, the case of this subsection can be
regarded as an extension of relevant cases appeared in
\cite{Hu-Jin-Zhou 2011}. Similar to \cite{Hu-Jin-Zhou 2011}, we
introduce the following notion.

\ms

\bf Definition 4.1. \rm For given $x\in\dbR^n$, a state-control pair
$(X^*(\cd),u^*(\cd))\in\sX[0,T]\times\sU[0,T]$ is called an {\it
open-loop equilibrium pair} of Problem (MF-LQ) for the initial state
$x$ if
$$X^*(0)=x,$$
and for almost all $t\in[0,T)$, and any $u(\cd)\in\sU[t,T]$,
\bel{local}\liminf_{\e\da0}{J(t,X^*(t);u^\e(\cd))-J(t,X^*(t);u^*(\cd))\over\e}\ge0,\ee
where
\bel{ue}u^\e(\cd)=u(\cd)I_{[t,t+\e)}(\cd)+u^*(\cd)I_{[t+\e,T]}(\cd).\ee
In this case, $X^*(\cd)$ and $u^*(\cd)$ are called an {\it open-loop
equilibrium state process} and an {\it open-loop equilibrium
control}, respectively.

\ms

We refer to (\ref{local}) as a {\it local optimality condition} at
$t\in[0,T)$. One sees that if $(X^*(\cd),u^*(\cd))$ is an open-loop
equilibrium pair of Problem (MF-LQ) for the initial state $x$, then
along the open-loop equilibrium state $X^*(\cd)$, the open-loop
equilibrium control $u^*(\cd)$ stays locally optimal. On the other
hand, since $(X^*(\cd),u^*(\cd))$ is a fixed state-control pair of
(\ref{1.1}) on $[0,T]$, in which the conditional expectation terms
are absent, the above defined open-loop equilibrium pair is
time-consistent. Note that if we consider the general state equation
(\ref{state}) in which some conditional expectation terms appear, we
do not know if one can directly define time-consistent state-control
pairs. This is why for open-loop equilibrium solutions of Problem
(MF-LQ), we only consider (\ref{1.1}). The following result, in some
sense, is an extension of a relevant result found in
\cite{Hu-Jin-Zhou 2011}.

\ms

\bf Proposition 4.2. \sl Let {\rm(H1)--(H2)} and $(\ref{=0})$ hold.
Suppose $(X^*(\cd),u^*(\cd))\in\sX[0,T]\times\sU[0,T]$ is a
state-control pair starting from initial state $x$. For each
$t\in[0,T)$, let $(Y(\cd\,,t),Z(\cd\,,t))$ be the adapted solution
of the following BSDE:
\bel{BSDE(t)}\left\{\2n\ba{ll}
\ns\ds dY(s,t)=-\Big\{A(s)^TY(s,t)+C(s)^TZ(s,t)+Q(s,t)X^*+\bar Q(s,t)
\dbE_t[X^*(s)]\Big\}ds+Z(s,t)dW(s),\\
\ns\ds\qq\qq\qq\qq\qq\qq\qq\qq\qq\qq\qq\qq\qq\qq\qq s\in[t,T],\\
\ns\ds Y(T,t)=G(t)X^*(T)+\bar G(t)\dbE_t[X^*(T)].\ea\right.\ee
Suppose $(s,t)\mapsto(Y(s,t),Z(s,t))$ is continuous on $0\le t\le
s\le T$ and suppose
\bel{4.5}u^*(t)=-\h
R(t,t)^{-1}\big\{B(t)^TY(t,t)+D(t)^TZ(t,t)\big\},\q t\in[0,T].\ee
Then $(X^*(\cd),u^*(\cd))$ is an open-loop equilibrium pair of
Problem (MF-LQ) for initial state $x$.

\ms

\it Proof. \rm For any fixed $t\in[0,T)$ and $u(\cd)\in\sU[t,T]$,
define $u^\e(\cd)$ by (\ref{ue}) and let $X^\e(\cd)\equiv
X(\cd\,;t,X^*(\cd),u^\e(\cd))$ be the solution to the following:
\bel{state(e)}\left\{\2n\ba{ll}
\ns\ds dX^\e(s)=\Big\{A(s)X^\e(s)+B(s)u^\e(s)\Big\}ds+\Big\{C(s)X^\e(s)+D(s)u^\e(s)\Big\}dW(s),\q s\in[t,T],\\
\ns\ds X^\e(t)=X^*(t).\ea\right.\ee
Then (suppressing $s$)
$$\ba{ll}
\ns\ds J(t,X^*(t);u^\e(\cd))-J(t,X^*(t),u^*(\cd))\\
\ns\ds=\dbE_t\Big\{\int_t^T\[\lan
Q(t)(X^\e+X^*),X^\e-X^*\ran+\lan\bar
Q(t)\dbE_t[X^\e+X^*],\dbE_t[X^\e-X^*]\ran\\
\ns\ds\qq\qq\qq+\lan R(t)(u^\e+u^*),u^\e-u^\e\ran+\lan\bar
R(t)\dbE_t[u^\e+u^*],\dbE_t[u^\e-u^*]\ran\]ds\\
\ns\ds\qq\qq\qq+\lan
G(t)\big[X^\e(T)+X^*(s)\big],X^\e(T)-X^*(T)\ran\\
\ns\ds\qq\qq\qq+\lan\bar
G(t)\dbE_t\big[X^\e(T)+X^*(T)\big],\dbE_t\big[X^\e(T)-X^*(T)\big]\ran\Big\}\ea$$
$$\ba{ll}
\ns\ds=\dbE_t\Big\{\int_t^T\lan Q(t)(X^\e+X^*)+\bar
Q(t)\dbE_t[X^\e+X^*],X^\e-X^*\ran ds\\
\ns\ds\qq\qq+\int_t^{t+\e}\lan R(t)(u+u^*)+\bar
R(t)\dbE_t[u+u^*],u-u^*\ran ds\\
\ns\ds\qq\qq\qq+\lan
G(t)\big[X^\e(T)+X^*(s)\big]+\bar G(t)\dbE_t\big[X^\e(T)+X^*(T)\big],X^\e(T)-X^*(T)\ran\Big\}.\ea$$
Let $(Y(\cd\,,t),Z(\cd\,,t))$ be the adapted solution to BSDE
(\ref{BSDE(t)}). Then one has (suppressing $(s,t)$)
$$\ba{ll}
\ns\ds d\big[\lan X^\e-X^*,Y\ran\big]=\[\lan A(X^\e-X^*)
+B(u^\e-u^*),Y\ran\\
\ns\ds\qq\qq\qq\qq\qq-\lan X^\e-X^*,\{A^TY+C^TZ+QX^*+\bar Q\dbE_t[X^*]\}\ran\\
\ns\ds\qq\qq\qq\qq\qq+\lan C(X^\e-X^*)+D(u^\e-u^*),Z\ran\]ds\\
\ns\ds\qq\qq\qq\qq\qq+\[\lan C(X^\e-X^*)+D(u^\e-u^*),Y\ran+\lan X^\e-X^*,Z\ran\]dW\\
\ns\ds\qq\qq\qq\qq=\[\lan u^\e-u^*,B^TY+D^TZ\ran-\lan
X^\e-X^*,QX^*+\bar Q\dbE_t[X^*]\ran\]ds\\
\ns\ds\qq\qq\qq\qq\qq+\[\lan C(X^\e-X^*)+D(u^\e-u^*),Y\ran+\lan
X^\e-X^*,Z\ran\]dW.\ea$$
Thus,
$$\ba{ll}
\ns\ds\dbE_t\big[\lan X^\e(T)-X^*(T),GX^*(T)+\bar
G\dbE_t[X^*(T)]\ran\big]\\
\ns\ds=\dbE_t\2n\int_t^T\2n\[\lan u^\e\1n-u^*,B^TY+D^T\1n Z\ran-\lan
X^\e\1n-\1n X^*,QX^*\1n+\1n\bar Q\dbE_t[X^*]\ran\]ds.\ea$$
Consequently,
$$\ba{ll}
\ns\ds J(t,X^*(t);u^\e(\cd))-J(t,X^*(t),u^*(\cd))\\
\ns\ds=\dbE_t\Big\{\int_t^T\[\lan Q(X^\e+X^*)+\bar
Q\dbE_t[X^\e+X^*]-2QX^*-2\bar Q\dbE_t[X^*],X^\e-X^*\ran\ea$$
$$\ba{ll}
\ns\ds\qq+\lan R(u^\e+u^*)+\bar R\dbE_t(u^\e+u^*)
+2B^TY+2D^TZ,u^\e-u^*\ran\]ds\\
\ns\ds\qq+\1n\lan\1n G\big[X^\e(T)\1n+\1n X^*(T)\big]\1n+\1n\bar
G\dbE_t\big[X^\e(T)\1n+\1n X^*(T)\big]\1n-\1n 2GX^*(T)\1n-\1n\bar
G\dbE_t[X^*(T)],X^\e(T)\1n-\1n X^*(T)\1n\ran\1n\Big\}\\
\ns\ds=\dbE_t\Big\{\int_t^T\[\lan Q(X^\e-X^*)+\bar
Q\dbE_t[X^\e-X^*],X^\e-X^*\ran\\
\ns\ds\qq+\lan R(u^\e+u^*)+\bar R\dbE_t(u^\e+u^*)
+2B^TY+2D^TZ,u^\e-u^*\ran\]ds\\
\ns\ds\qq+\1n\lan\1n G\big[X^\e(T)\1n-\1n X^*(T)\big]\1n+\1n\bar
G\dbE_t\big[X^\e(T)\1n-\1n X^*(T)\big],X^\e(T)\1n-\1n
X^*(T)\1n\ran\1n\Big\}\\
\ns\ds=\dbE_t\Big\{\int_t^{t+\e}\lan R(t)[u+u^*]+\bar
R(t)\dbE_t[u+u^*]+2B^TY+2D^TZ,u-u^*\ran ds\\
\ns\ds\qq+\int_t^T|
Q(t)^{1\over2}\big\{X^\e-X^*-\dbE_t[X^\e-X^*]\big\}|^2+|\h Q(t)^{1\over2}\dbE_t[X^\e-X^*]|^2ds\\
\ns\ds\qq+|G(t)^{1\over2}\big\{X^\e(T)-X^*(T)-\dbE_t[X^\e(T)-X^*(T)]\big\}|^2
+|\h G(t)^{1\over2}\dbE_t[X^\e(T)-X^*(T)]|^2\Big\}.\ea$$
Let us denote
$$\left\{\ba{ll}
\ns\ds\h u(s)=u(s)-u^*(s),\\
\ns\ds\th(s,t)=R(s,t)u^*(s)+\bar
R(s,t)\dbE_t[u^*(s)]+B(s)^TY(s,t)+D(s)^TZ(s,t).\ea\right.$$
Then
$$\ba{ll}
\ns\ds\dbE_t\[\lan R(t)(u+u^*)+\bar R(t)\dbE_t[u+u^*] +2B^TY+2\bar
B^T\dbE_t[Y]+2D^TZ+2\bar D^T\dbE_t[Z],u-u^*\ran\Big]\\
\ns\ds=\dbE_t\big\{\lan R(t)\h u+\bar R(t)\dbE_t[\h u]+2\th,\h u\ran\big\}\\
\ns\ds=\dbE_t\[\lan R(t)\big\{\h u-\dbE_t[\h u]\big\}+\h
R(t)\dbE_t[\h u]
+2\th,\h u-\dbE_t[\h u]+\dbE_t[\h u]\ran\]\\
\ns\ds=\dbE_t\Big\{|R(t)^{1\over2}\big\{\h u-\dbE_t[\h
u]\big\}|^2+2\lan\th,\h u-\dbE_t[\h u]\ran+|\h
R(t)^{1\over2}\dbE_t[\h u]|^2+2\lan\th,\dbE_t[\h
u]\ran\Big\}\\
\ns\ds=\dbE_t\Big\{\big|R(t)^{1\over2}\big\{\h u-\dbE_t[\h
u]\big\}+R(t)^{-{1\over2}}\th\big|^2+\big|\h R(t)^{1\over2}\dbE_t[\h
u]+\h R(t)^{-{1\over2}}\dbE_t[\th]\big|^2\\
\ns\ds\qq-\lan R(t)^{-1}\th,\th\ran-\lan\h
R(t)^{-1}\dbE_t[\th],\dbE_t[\th]\ran\Big\}.\ea$$
Hence,
$$\ba{ll}
\ns\ds J(t,X^*(t);u^\e(\cd))-J(t,X^*(t);u^*(\cd))\\
\ns\ds=\dbE_t\Big\{\int_t^{t+\e}\[\big|\h
R(t)^{1\over2}\dbE_t[u-u^*]+\h R(t)^{-{1\over2}}\dbE_t[\th]\big|^2
+\big|R(t)^{1\over2}\{\h u-\dbE_t[\h u]\}+R(t)^{-{1\over2}}\th\big|^2\]ds\\
\ns\ds\qq+\int_t^T\[|
Q(t)^{1\over2}\big\{X^\e-X^*-\dbE_t[X^\e-X^*]\big\}|^2+|\h Q(t)^{1\over2}\dbE_t[X^\e-X^*]|^2\]ds\\
\ns\ds\qq+|G(t)^{1\over2}\big\{X^\e(T)-X^*(T)-\dbE_t[X^\e(T)-X^*(T)]\big\}|^2
+|\h G(t)^{1\over2}\dbE_t[X^\e(T)-X^*(T)]|^2\\
\ns\ds\qq-\int_t^{t+\e}\(\lan R(t)^{-1}\th,\th\ran+\lan\h R(t)^{-1}\dbE_t[\th],\dbE_t[\th]\ran\)ds\Big\}\\
\ns\ds\ge-\dbE_t\int_t^{t+\e}\(\lan
R(s,t)^{-1}\th(s,t),\th(s,t)\ran+\lan\h
R(s,t)^{-1}\dbE_t[\th(s,t)],\dbE_t[\th(s,t)]\ran\)ds.\ea$$
Note that under our continuity condition (see (\ref{4.5})),
\bel{}\lim_{s\da t}\dbE_t|\th(s,t)|^2=|\th(t,t)|^2=0,\qq
t\in[0,T].\ee
Hence, (\ref{local}) holds and $(X^*(\cd),u^*(\cd))$ is an open-loop
equilibrium pair of Problem (MF-LQ) for the initial state $x$.
\endpf

\ms

The above result leads to the following FBSDE family (parameterized
by $t\in[0,T)$):
\bel{FBSDE}\left\{\2n\ba{ll}
\ns\ds dX^*(s)=\Big\{A(s)X^*(s)+B(s)u^*(s)\Big\}ds+\Big\{C(s)X^*(s)+D(s)u^*(s)\Big\}dW(s),\q s\in[0,T],\\
\ns\ds dY(s,t)=-\Big\{A(s)^TY(s,t)+C(s)^TZ(s,t)+Q(s,t)X^*(s)+\bar
Q(s,t)\dbE_t[X^*(s)]\Big\}ds+Z(s,t)dW(s),\\
\ns\ds\qq\qq\qq\qq\qq\qq\qq\qq\qq\qq\qq\qq\qq\qq\qq0\le t\le s\le T,\\
\ns\ds X^*(0)=x,\qq Y(T,t)=G(t)X^*(T)+\bar
G(t)\dbE_t[X^*(T)],\\
\ns\ds\h R(t,t)u^*(t)+B(t)^TY(t,t)+D(t)^TZ(t,t)=0,\q
t\in[0,T].\ea\right.\ee
Inspired by \cite{MY 1999}, we suppose
$$Y(s,t)=P(s,t)X^*(s)+\bar P(s,t)\dbE_t[X^*(s)],\qq s\in[t,T],$$
for some deterministic functions $P(\cd\,,\cd)$ and $\bar
P(\cd\,,\cd)$. Then
$$\ba{ll}
\ns\ds-\Big\{A(s)^TY(s,t)\1n+\1n C(s)^T\1n Z(s,t)\1n+\1n
Q(s,t)X^*(s)\1n+\1n\bar
Q(s,t)\dbE_t[X^*(s)]\Big\}ds\1n+\1n Z(s,t)dW(s)\1n=\1n dY(s,t)\\
\ns\ds=\Big\{P_s(s,t)X^*(s)+P(s,t)[A(s)X^*(s)+B(s)u^*(s)]+\bar
P_s(s,t)\dbE_t[X^*(s)]\\
\ns\ds\qq+\bar
P(s,t)\(A(s)\dbE_t[X^*(s)]+B(s)\dbE_t[u^*(s)]\)\Big\}ds+P(s,t)\(C(s)X^*(s)+D(s)u^*(s)\)dW(s).\ea$$
Hence, we need
$$Z(s,t)=P(s,t)\(C(s)X^*(s)+D(s)u^*(s)\).$$
Consequently,
$$\ba{ll}
\ns\ds\h R(t,t)u^*(t)+B(t)^TY(t,t)=-D(t)^TZ(t,t)=-D(t)^TP(t,t)\(C(t)X^*(t)+D(t)u^*(t)\).\ea$$
Hence, (with $\h P(s,t)=P(s,t)+\bar P(s,t)$)
\bel{open-u}u^*(t)=-\big[\h
R(t,t)+D(t)^TP(t,t)D(t)\big]^{-1}\[B(t)^T\h
P(t,t)+D(t)^TP(t,t)C(t)\]X^*(t),\ee
and
\bel{}\ba{ll}
\ns\ds Z(s,t)=P(s,t)\(C(s)X^*(s)+D(s)u^*(s)\)\\
\ns\ds\qq\q=P(s,t)C(s)X^*(s)-\1n P(s,t)D(s)\big[\h R(s,s)\1n+\1n
D(s)^T\1n
P(s,s)D(s)\big]^{-1}\\
\ns\ds\qq\qq\q\cd\1n\[B(s)^T\h
P(s,s)+D(s)^TP(s,s)C(s)\]X^*(s).\ea\ee
Therefore,
$$\ba{ll}
\ns\ds0=P_s(s,t)X^*(s)+P(s,t)[A(s)X^*(s)+B(s)u^*(s)]+\bar
P_s(s,t)\dbE_t[X^*(s)]\\
\ns\ds\qq+\bar
P(s,t)\(A(s)\dbE_t[X^*(s)]+B(s)\dbE_t[u^*(s)]\)\\
\ns\ds\qq+A(s)^TY(s,t)+C(s)^TZ(s,t)+Q(s,t)X^*(s)+\bar
Q(s,t)\dbE_t[X^*(s)]\\
\ns\ds\q=P_s(s,t)X^*(s)+P(s,t)A(s)X^*(s)\\
\ns\ds\qq-P(s,t)B(s)[\h R(s,s)+D(s)^TP(s,s)D(s)]^{-1}[B(s)^T\h P(s,s)
+D(s)^TP(s,s)C(s)]X^*(s)\\
\ns\ds\qq+\bar P_s(s,t)\dbE_t[X^*(s)]+\bar
P(s,t)A(s)\dbE_t[X^*(s)]\\
\ns\ds\qq-\bar P(s,t)B(s)[\h R(s,s)+D(s)^TP(s,s)D(s)]^{-1}[B(s)^T\h P(s,s)
+D(s)^TP(s,s)C(s)]\dbE_t[X^*(s)]\\
\ns\ds\qq+A(s)^TP(s,t)X^*(s)+A(s)^T\bar
P(s,t)\dbE_t[X^*(s)]+C(s)^TP(s,t)C(s)X^*(s)\\
\ns\ds\qq-C(s)^T\1n P(s,t)D(s)\big[\h R(s,s)\1n+\1n D(s)^T\1n
P(s,s)D(s)\big]^{-1}
\big[B(s)^T\1n\h P(s,s)\1n+\1n D(s)^T\1n P(s,s)C(s)\big]X^*(s)\\
\ns\ds\qq+Q(s,t)X^*(s)+\bar Q(s,t)\dbE_t[X^*(s)]\\
\ns\ds\q=\Big\{P_s(s,t)+P(s,t)A(s)+A(s)^TP(s,t)+C(s)^TP(s,t)C(s)+Q(s,t)\\
\ns\ds\qq-[P(s,t)B(s)+C(s)^TP(s,t)D(s)]\big[\h
R(s,s)+D(s)^TP(s,s)D(s)\big]^{-1}\\
\ns\ds\qq\qq\cd
[B(s)^T\h P(s,s)+D(s)^TP(s,s)C(s)]\Big\}X^*(s)\\
\ns\ds\qq+\Big\{\bar P_s(s,t)+\bar P(s,t)A(s)+A(s)^T\bar
P(s,t)+\bar Q(s,t)\\
\ns\ds\qq-\bar P(s,t)B(s)[\h R(s,s)\1n+\1n D(s)^T\1n
P(s,s)D(s)]^{-1}[B(s)^T\1n\h P(s,s)\1n+\1n D(s)^T\1n
P(s,s)C(s)]\Big\}\dbE_t[X^*(s)].\ea$$
The above will be true if we let $P(s,t)$ and $\h P(s,t)\equiv
P(s,t)+\bar P(s,t)$ be the solutions to the following coupled
Riccati equations:
\bel{open-Riccati P}\left\{\2n\ba{ll}
\ns\ds P_s(s,t)+P(s,t)A(s)+A(s)^TP(s,t)+C(s)^TP(s,t)C(s)+Q(s,t)\\
\ns\ds\qq-[P(s,t)B(s)\1n+\1n C(s)^T\1n P(s,t)D(s)][\h R(s,s)\1n+\1n
D(s)^TP(s,s)D(s)]^{-1}\\
\ns\ds\qq\qq\cd[B(s)^T\h P(s,s)+D(s)^T\1n P(s,s)C(s)]=0,\\
\ns\ds\h P_s(s,t)+\h P(s,t)A(s)+A(s)^T\h P(s,t)+C(s)^TP(s,t)C(s)+\h Q(s,t)\\
\ns\ds\qq-[\h P(s,t)B(s)\1n+\1n C(s)^T\1n P(s,t)D(s)][\h
R(s,s)\1n+\1n
D(s)^TP(s,s)D(s)]^{-1}\\
\ns\ds\qq\qq\cd[B(s)^T\h P(s,s)+D(s)^TP(s,s)C(s)]=0,\qq s\in[t,T],\\
\ns\ds P(T,t)=G(t),\qq\h P(T,t)=\h G(t).\ea\right.\ee
To summarize the above, we state the following result.

\ms

\bf Theorem 4.3. \sl Let {\rm(H1)--(H2)} hold. Suppose Riccati
equation system $(\ref{open-Riccati P})$ admits a unique solution
$(P(\cd\,,\cd),\h P(\cd\,,\cd))$ which is continuous in both
variables. Then $u^*(\cd)\in\sU[0,T]$ defined by $(\ref{open-u})$ is
an open-loop equilibrium control and the corresponding $X^*(\cd)$ is
an open-loop equilibrium state process.

\ms

\rm

On top of the above result, one may study the solvability of Riccati
equation system (\ref{open-Riccati P}). In principle, the above
gives a time-consistent solution to Problem (MF-LQ) of open-loop
type. Let us look at a couple special cases:

\ms

\it Case 1. The case of general discounting only. \rm Let
\bel{general discounting}\bar A(\cd)=\bar C(\cd)=0,\q\bar
B(\cd)=\bar D(\cd)=0,\q\bar Q(\cd\,,\cd)=0,\q\bar
R(\cd\,,\cd)=0,\q\bar G(\cd)=0.\ee
In this case, we have general discounting without conditional
expectations. Clearly,
$$\left\{\2n\ba{ll}
\ns\ds\h A(\cd)=A(\cd),\q\h B(\cd)=B(\cd),\q\h C(\cd)=C(\cd),\q\h
D(\cd)=D(\cd),\\
\ns\ds\h Q(\cd\,,\cd)=Q(\cd\,,\cd),\q\h
R(\cd\,,\cd)=R(\cd\,,\cd),\q\h G(\cd)=G(\cd),\ea\right.$$
and $P(\cd\,,\cd)=\h P(\cd\,,\cd)$ satisfying
\bel{discounting-Riccati P}\left\{\2n\ba{ll}
\ns\ds P_s(s,t)+P(s,t)A(s)+A(s)^TP(s,t)+C(s)^TP(s,t)C(s)+Q(s,t)\\
\ns\ds\qq-[P(s,t)B(s)\1n+\1n C(s)^T\1n P(s,t)D(s)][R(s,s)\1n+\1n
D(s)^TP(s,s)D(s)]^{-1}\\
\ns\ds\qq\qq\cd[B(s)^TP(s,s)+D(s)^T\1n P(s,s)C(s)]=0,\qq s\in[t,T],\\
\ns\ds P(T,t)=G(t).\ea\right.\ee
The open-loop equilibrium control $u^*(\cd)$ is given by
\bel{open-u-discounting}u^*(t)=-\big[R(t,t)+D(t)^TP(t,t)D(t)\big]^{-1}
\big[B(t)^TP(t,t)+D(t)^TP(t,t)C(t)\big]X^*(t),\q t\in[0,T].\ee
In addition, if $C(\cd)=0$ and $D(\cd)=0$, the problem becomes a
deterministic LQ problem with general discounting and we have
\bel{}\left\{\2n\ba{ll}
\ns\ds P_s(s,t)\1n+\1n P(s,t)A(s)\1n+\1n A(s)^T\1n P(s,t)\1n+\1n
Q(s,t)
\1n-\1n P(s,t)B(s)R(s,s)^{-1}\1n B(s)^T\1n P(s,s)=0,\q s\in[t,T],\\
\ns\ds P(T,t)=G(t),\ea\right.\ee
with the open-loop equilibrium control $u^*(\cd)$ given by
$$u^*(t)=-R(t,t)^{-1}B(t)^TP(t,t)X^*(t),\q t\in[0,T].$$
This result is comparable with those in \cite{Yong 2011} where
closed-loop equilibrium solution was obtained for deterministic LQ
problems with general discounting (without conditional expectation
terms).

\ms

\it Case 2. The case of conditional expectation at terminal only.
\rm Let
$$\left\{\ba{ll}
\ns\ds Q(s,t)=Q(s),\qq\bar Q(s,t)=0,\qq R(s,t)=R(s),\qq\bar R(s,t)=0,\\
\ns\ds G(t)=G,\qq\bar G(t)=\bar G.\ea\right.$$
In this case, the conditional expectation appears at the terminal
cost only without general discounting. This is a case studied in
\cite{Hu-Jin-Zhou 2011}. For such a case, both $P(s,t)$ and $\h
P(s,t)$ are independent of $t$ and (\ref{open-Riccati P}) becomes
\bel{open-Riccati P*}\left\{\2n\ba{ll}
\ns\ds\dot P(s)+P(s)A(s)+A(s)^TP(s)+C(s)^TP(s)C(s)+Q(s)\\
\ns\ds\q-[P(s)B(s)\1n+\1n C(s)^T\1n P(s)D(s)][R(s)\1n+\1n
D(s)^TP(s)D(s)]^{-1}[B(s)^T\h P(s)+D(s)^T\1n P(s)C(s)]=0,\\
\ns\ds\dot{\h P}(s)+\h P(s)A(s)+A(s)^T\h P(s)+C(s)^TP(s)C(s)+Q(s)\\
\ns\ds\q-[\h P(s)B(s)\1n+\1n C(s)^T\1n P(s)D(s)][R(s)\1n+\1n
D(s)^TP(s)D(s)]^{-1}[B(s)^T\h P(s)+D(s)^TP(s)C(s)]=0,\\
\ns\ds\qq\qq\qq\qq\qq\qq\qq\qq\qq\qq\qq\qq\qq\qq s\in[t,T],\\
\ns\ds P(T)=G,\qq\h P(T)=\h G,\ea\right.\ee
and the open-loop equilibrium control is given by
$$u^*(t)=-\big[R(t)+D(t)^TP(t)D(t)\big]^{-1}\big[B(t)^T\h P(t)+D(t)^TP(t)C(t)\big]
X^*(t),\q t\in[0,T].$$
This essentially covers a relevant result in \cite{Hu-Jin-Zhou
2011}. Note that as long as $\bar G\ne0$, $P(\cd)\ne\h P(\cd)$. Then
the equation for $P(\cd)$ contains $\h P(\cd)$ and it does not have
a desired symmetry. Therefore, $P(\cd)$ is not expected to be
symmetric. Consequently, $\h P(\cd)$ will not be symmetric either.
Finally, if, in addition, $\bar G=0$, then $P(\cd)=\h P(\cd)$ and it
satisfies
\bel{open-Riccati P**}\left\{\2n\ba{ll}
\ns\ds\dot P(s)+P(s)A(s)+A(s)^TP(s)+C(s)^TP(s)C(s)+Q(s)\\
\ns\ds\q-[P(s)B(s)\1n+\1n C(s)^T\1n P(s)D(s)][R(s)\1n+\1n
D(s)^T\1n P(s)D(s)]^{-1}[B(s)^T\1n P(s)\1n+\1n D(s)^T\1n P(s)C(s)]=0,\\
\ns\ds\qq\qq\qq\qq\qq\qq\qq\qq\qq\qq\qq\qq\qq\qq\qq\qq s\in[t,T],\\
\ns\ds P(T)=G.\ea\right.\ee
In this case, the problem is reduced to a classical stochastic LQ
problem and it is time-consistent. For this case, the open-loop
equilibrium control $u^*(\cd)$ is the optimal control and is given
by
$$u^*(t)=-\big[R(t)+D(t)^TP(t)D(t)\big]^{-1}\big[B(t)^TP(t)+D(t)^TP(t)C(t)\big]
X^*(t),\q t\in[0,T],$$
which recovers the result for classical stochastic LQ problem.

\ms

We now make a couple of comments on this.

\ms

\it The advantages: \rm The approach is direct and the derivation of
equilibrium pair is not very complicated. Moreover, the open-loop
equilibrium control $u^*(\cd)$ admits a closed-loop representation
(\ref{open-u}).

\ms

\it The disadvantages: \rm (i) The Riccati equations in
(\ref{open-Riccati P}) do not have symmetry structure. Therefore the
solutions $P(\cd\,,\cd)$ and $\h P(\cd\,,\cd)$ of the system are not
necessarily symmetric. This leads to some difficulties in establish
the well-posedness of the system. (ii) If the state equation
contains conditional expectation, even the definition of open-loop
equilibrium pair is not clear to us.

\ms

\subsection{Close-loop equilibrium strategy.}

In this subsection, we introduce closed-loop equilibrium strategies.
To this end, we first introduce the following: For any $t\in[0,T)$,
\bel{}\ba{ll}
\ns\ds\wt J(t;X(\cd),u(\cd))=\dbE_t\Big\{\int_t^T\[\lan
Q(s,t)X(s),X(s)\ran+\lan\bar
Q(s,t)\dbE_t[X(s)],\dbE_t[X(s)]\ran\\
\ns\ds\qq\qq\qq\qq\qq\qq+\lan R(s,t)u(s),u(s)\ran+\lan\bar
R(s,t)\dbE_t[u(s)],\dbE_t[u(s)]\ran\]ds\\
\ns\ds\qq\qq\qq\qq\qq\qq+\lan G(t)X(T),X(T)\ran+\lan\bar
G(t)\dbE_t[X(T)],\dbE_t[X(T)]\ran\Big\},\ea\ee
for any $(X(\cd),u(\cd))\in\sX[t,T]\times\sU[t,T]$. We point out
that in the above $(X(\cd),u(\cd))$ does not have to be a
state-control pair of the original control system. Thus, $\wt
J(t;X(\cd),u(\cd))$ is an extension of the cost functional
$J(t,x;u(\cd))$, and
$$\wt
J\big(t;X(\cd\,;t,x,u(\cd)),u(\cd)\big)=J(t,x;u(\cd)),\qq\forall(t,x)\in\sD,~
u(\cd)\in\sU[t,T].$$
Next, and hereafter, we denote any partition of $[0,T]$ by $\D$:
$$\D=\big\{t_k\bigm|0\le k\le N\big\}\equiv\big\{0=t_0<t_1<t_2<\cds<t_{N-1}
<t_N=T\big\},$$
with $N$ being some natural number, and define its {\it mesh size}
by the following:
$$\|\D\|=\max_{0\le k\le N-1}(t_{k+1}-t_k).$$
For the above $\D$, we define \bel{}\ba{ll}
\ns\ds J^\D_k(X(\cd),u(\cd))=\dbE_{t_k}\Big\{\int_{t_k}^T\[\lan
Q(s,t_k)X(s),X(s)\ran+\lan\bar
Q(s,t_k)\dbE_{t_k}[X(s)],\dbE_{t_k}[X(s)]\ran\\
\ns\ds\qq\qq\qq\qq\qq\qq+\lan R(s,t_k)u(s),u(s)\ran+\lan\bar
R(s,t_k)\dbE_{t_k}[u(s)],\dbE_{t_k}[u(s)]\ran\]ds\\
\ns\ds\qq\qq\qq\qq\qq\qq+\lan G(t_k)X(T),X(T)\ran+\lan\bar
G(t_k)\dbE_{t_k}[X(T)],\dbE_{t_k}[X(T)]\ran\Big\},\ea\ee
for any $(X(\cd),u(\cd))\in\sX[t_k,T]\times\sU[t_k,T]$,
$k=0,1,2,\cds,N-1$. Again, in the above, $(X(\cd),u(\cd))$ does not
have to be a state-control pair of the original control system.

\ms

Now, we introduce some notions.

\ms

\bf Definition 4.4. \rm Let
$\D=\{0=t_0<t_1<\cds<t_{N-1}<t_N=T\big\}$ be a partition of $[0,T]$,
and let $\Th^\D,\h\Th^\D:[0,T]\to\dbR^{m\times n}$ be two given
maps, possibly depending on $\D$.

\ms

(i) For any $x\in\dbR^n$ fixed, let $X^\D(\cd)\equiv X^\D(\cd\,;x)$
be the solution to the following linear MF-SDE:
\bel{closed-loop D}\left\{\2n\ba{ll}
\ns\ds dX^\D\1n(s)\1n=\2n\Big\{\1n\big[A(s)\1n-\1n
B(s)\Th^\D(s)\big]\1n
X^\D\1n(s)\\
\ns\ds\qq\qq+\1n\[\bar A(s)\1n+\1n
B(s)\1n\big[\Th^\D(s)\1n-\1n\h\Th^\D(s)\big]\1n-\1n\bar
B(s)\h\Th^\D(s)\]\dbE_{\rho^\D(s)}[X^\D\1n(s)]\Big\}ds\\
\ns\ds\qq\qq+\Big\{\1n\big[C(s)\1n-\1n
D(s)\Th^\D(s)\big]X^\D(s)\\
\ns\ds\qq\qq+\1n\[\bar C(s)\1n+\2n
D(s)\1n\big[\Th^\D(s)\1n-\1n\h\Th^\D(s)\big]\1n-\1n\bar
D(s)\h\Th^\D(s)\]\dbE_{\rho^\D(s)}[X^\D(s)]\Big\}dW(s),~s\1n\in\1n[0,\1n T],\\
\ns\ds X^\D(0)=x,\ea\right.\ee
where
$$\rho^\D(s)=\sum_{k=0}^{N-1}t_kI_{[t_k,t_{k+1})}(s),\qq s\in[0,T],$$
and let $u^\D(\cd)\equiv u^\D(\cd\,;x)$ be defined by
\bel{}u^\D(s)=-\Th^\D(s)X^\D(s)+\big[\Th^\D(s)-\h\Th^\D(\cd)\big]
\dbE_{\rho^\D(s)}[X^\D(s)],\qq s\in[0,T].\ee
The pair $(X^\D(\cd),u^\D(\cd))$ is called the {\it closed-loop
pair} associated with $\D$ and $(\Th^\D(\cd),\h\Th^\D(\cd))$,
starting from $x$.

\ms

(ii) For each $t_k\in\D$ and any $u_k(\cd)\in\sU[t_k,t_{k+1}]$, let
$X_k(\cd)$ be the solution to the following system:
\bel{closed-loop D}\left\{\2n\ba{ll}
\ns\ds dX_k(s)\1n=\1n\Big\{A(s)X_k(s)+\bar
A(s)\dbE_{t_k}[X_k(s)]+B(s)u_k(s)+\bar
B(s)\dbE_{t_k}[u_k(s)]\Big\}ds\\
\ns\ds\qq\qq+\Big\{C(s)X_k(s)\1n+\1n\bar
C(s)\dbE_{t_k}[X_k(s)]\1n+\1n D(s)u_k(s)\1n+\1n\bar
D(s)\dbE_{t_k}[u_k(s)]\Big\}dW(s),\q s\1n\in\1n[t_k,t_{k+1}],\\
\ns\ds X_k(t_k)=X^\D(t_k),\ea\right.\ee
and $X^\D_{k+1}(\cd)$ be the solution to the following:
\bel{closed-loop D(k+1)}\left\{\2n\ba{ll}
\ns\ds dX^\D_{k+1}\1n(s)\1n=\2n\Big\{\1n\big[A(s)\1n-\1n
B(s)\Th^\D(s)\big]\1n
X^\D_{k+1}\1n(s)\\
\ns\ds\qq\qq+\1n\[\bar A(s)\1n+\1n
B(s)\1n\big[\Th^\D(s)\1n-\1n\h\Th^\D(s)\big]\1n-\1n\bar
B(s)\h\Th^\D(s)\]\dbE_{\rho^\D(s)}[X^\D_{k+1}\1n(s)]\Big\}ds\\
\ns\ds\qq\qq+\Big\{\1n\big[C(s)\1n-\1n
D(s)\Th^\D(s)\big]X^\D_{k+1}(s)\\
\ns\ds\qq\qq+\1n\[\bar C(s)\1n+\2n
D(s)\1n\big[\Th^\D\1n(s)\2n-\1n\h\Th^\D\1n(s)\big]\2n-\2n\bar
D(s)\h\Th^\D(s)\]\dbE_{\rho^\D(s)}[X^\D_{k\1n+\1n1}(s)]\Big\}dW\1n(s),~s\1n\in
\1n[t_{k\1n+\1n1},\1n T],\\
\ns\ds X^\D_{k+1}(t_{k+1})=X_k(t_{k+1}).\ea\right.\ee
 Denote
\bel{variation}\left\{\ba{ll}
\ns\ds X_k(\cd)\oplus X^\D(\cd)\equiv
X_k(\cd)I_{[t_k,t_{k+1})}(\cd)+X^\D_{k+1}(\cd)I_{[t_{k+1},T]}(\cd),\\
\ns\ds u_k(\cd)\oplus
u^\D(\cd)=u_k(\cd)I_{[t_k,t_{k+1})}(\cd)-\big\{\Th^\D(\cd)
X_{k+1}^\D(\cd)\1n+\2n\big[\Th^\D(\cd)\1n-\1n\h\Th^\D(\cd)\big]\dbE_{\rho^\D(\cd)}[X_{k+1}^\D(\cd)]\big\}I_{[t_{k+1},T]}(\cd).
\ea\right.\ee
We call $(X_k(\cd)\oplus X^\D(\cd),u_k(\cd)\oplus u^\D(\cd))$ a {\it
local variation} of $(X^\D(\cd),u^\D(\cd))$ on $[t_k,t_{k+1}]$.
Suppose the following {\it local optimality condition} holds:
\bel{}\ba{ll}
\ns\ds J^\D_k\big(X_k^\D(\cd),u_k^\D(\cd)\big)\le J_k^\D\big(
X_k(\cd)\oplus X^\D(\cd),u_k(\cd)\oplus u^\D(\cd)\big),\qq\forall
u_k(\cd)\in\sU[t_k,t_{k+1}].\ea\ee
Then we call $(\Th^\D(\cd),\h\Th^\D(\cd))$ a {\it closed-loop
$\D$-equilibrium strategy} of Problem (MF-LQ), and call
$\big(X^\D(\cd\,;x),$\par\no$u^\D(\cd\,;x)\big)$ a {\it closed-loop
$\D$-equilibrium pair} of Problem (MF-LQ) for the initial state $x$.

\ms

(iii) If the following holds:
\bel{}\lim_{\|\D\|\to0}\[\|\Th^\D(\cd)-\Th(\cd)\|_{C([0,T];\dbR^{m\times
n})}+\|\h\Th^\D(\cd)-\h\Th(\cd)\|_{C([0,T];\dbR^{m\times
n})}\]=0,\ee
for some $\Th,\h\Th\in C([0,T];\dbR^{m\times n})$, then
$(\Th(\cd),\h\Th(\cd))$ is called a {\it closed-loop equilibrium
strategy} of Problem (MF-LQ). For any $(t,x)\in\sD$, let $\h
X^*(\cd)\equiv\h X^*(\cd\,;t,x)$ be the solution to the following
system:
\bel{closed-loop}\left\{\2n\ba{ll}
\ns\ds d\h X^*\1n(s)\1n=\2n\big[\h A(s)\1n-\1n \h
B(s)\h\Th(s)\big]\h X^*(s)ds+\big[\h C(s)\1n-\1n\h D(s)\h\Th(s)\big]
\h X^*(s)dW(s),\qq s\in[t,T],\\
\ns\ds\h X^*(t)=x,\ea\right.\ee
and define $\h u^*(\cd)\equiv\h u^*(\cd\,;t,x)$ as follows:
\bel{}\h u^*(s)=-\h\Th(s)\h X^*(s),\qq s\in[t,T].\ee
Then $(t,x)\mapsto(\h X^*(\cd\,;t,x)$, $\h u^*(\cd\,;t,x))$ is
called a {\it closed-loop equilibrium pair flow} of Problem (MF-LQ).
Further,
\bel{}\h V(t,x)=\wt J(t,x;\h X^*(\cd\,;t,x),\h
u^*(\cd\,;t,x)),\qq(t,x)\in\sD\ee
is called a {\it closed-loop equilibrium value function} of Problem
(MF-LQ).

\ms

We point out that $(\Th^\D(\cd),\h\th^\D(\cd))$ and
$(\Th(\cd),\h\Th(\cd))$ are independent of the initial state
$x\in\dbR^n$. Let us now state the main result of this paper.

\ms

\bf Theorem 4.5. \sl Let {\rm(H1)}--{\rm(H3)} hold. Then there
exists a unique pair $(\G(\cd\,,\cd),\h\G(\cd\,,\cd))$ of
$\dbS^n$-valued functions solving the following system of equations:
\bel{Riccati-G(s,t)}\left\{\2n\ba{ll}
\ns\ds\G_s(s,t)+\G(s,t)\big[\h A(s)-\h B(s)\h\Th(s)\big]
+\big[\h A(s)-\h B(s)\h\Th(s)\big]^T\G(s,t)+Q(s,t)\\
\ns\ds\q+\big[\h C(s)\1n-\1n\h D(s)\h\Th(s)\big]^T\1n\G(s,t)\big[\h
C(s)\1n-\1n\h D(s)\h\Th(s)\big]\1n+\1n\Th(s)^T\1n R(s,t)\Th(s)=0,\\
\ns\ds\h\G_s(s,t)+\h\G(s,t)\big[\h A(s)-\h
B(s)\h\Th(s)\big]+\big[\h A(s)-\h B(s)\h\Th(s)\big]^T\h\G(s,t)+\h Q(s,t)\\
\ns\ds\q+\big[\h C(s)\1n-\1n\h D(s)\h\Th(s)\big]^T\1n\G(s,t)\big[\h
C(s)\1n -\1n\h D(s)\h\Th(s)\big]\1n+\1n\h\Th(s)^T\1n
\h R(s,t)\h\Th(s)=0,\qq0\le t\le s\le T,\\
\ns\ds\G(T,t)=G(t),\qq\h\G(T,t)=\h G(t),\qq0\le t\le
T,\ea\right.\ee
where $\h\Th(\cd)$ is given by the following:
\bel{Th(N-1)}\h\Th(s)=\big[\h R(s,s)+\h D(s)^T\G(s,s)\h
D(s)\big]^{-1}\big[\h B(s)^T\h \G(s,s)+\h D(s)^T\G(s,s)\h
C(s)\big],\q s\in[0,T].\ee
The closed-loop equilibrium state process $X^*(\cd)$ is the solution to the
following system:
\bel{}\left\{\2n\ba{ll}
\ns\ds dX^*(s)=\big[\h A(s)-\h B(s)\h\Th(s)\big]X^*(s)ds+\big[\h C(s)-\h D(s)\h\Th(s)\big]X^*(s)dW(s),\q s\in[0,T],\\
\ns\ds X^*(0)=x,\ea\right.\ee
the closed-loop equilibrium control admits the following representation:
\bel{feedback(N-2)e} u^*(s)=-\h\Th(s)X^*(s),\q s\in[0,T],\ee
and the closed-loop equilibrium value function is given by the following:
\bel{V}\h V(t,x)=\lan\h\G(t,t)x,x\ran,\qq\forall(t,x)\in\sD.\ee

\ms

\rm

Note that in (\ref{Riccati-G(s,t)}), the equations for $\G(\cd\,,\cd)$ and
$\h\G(\cd\,,\cd)$ are different: $(Q(\cd\,,\cd),R(\cd\,,\cd),G(\cd))$ appears in
the former and $(\h Q(\cd\,,\cd),\h R(\cd\,,\cd),\h G(\cd))$ appears in the
later. Also, we see that the system is fully coupled.

\ms

Let us look at two special cases.

\ms

\it Case 1. The case of general discounting only. \rm As in the
previous subsection, let
$$\bar A(\cd)=\bar C(\cd)=0,\q\bar B(\cd)=\bar
D(\cd)=0,\q\bar Q(\cd\,,\cd)=0,\q\bar R(\cd\,,\cd)=0,\q\bar
G(\cd)=0.$$
In this case, we have
$$\left\{\2n\ba{ll}
\ns\ds\h A(\cd)=A(\cd),\q\h B(\cd)=B(\cd),\q\h C(\cd)=C(\cd),\q\h
D(\cd)=D(\cd),\\
\ns\ds\h Q(\cd\,,\cd)=Q(\cd\,,\cd),\q\h
R(\cd\,,\cd)=R(\cd\,,\cd),\q\h G(\cd)=G(\cd).\ea\right.$$
Thus,
\bel{Th(N-1)}\h\Th(s)=\big[R(s,s)+D(s)^T\G(s,s)D(s)\big]^{-1}\big[B(s)^T\h
\G(s,s)+D(s)^T\G(s,s)C(s)\big],\q s\in[0,T].\ee
Hence, (\ref{Riccati-G(s,t)}) becomes
\bel{}\left\{\2n\ba{ll}
\ns\ds\G_s(s,t)+\G(s,t)\big[A(s)-B(s)\h\Th(s)\big]
+\big[A(s)-B(s)\h\Th(s)\big]^T\G(s,t)+Q(s,t)\\
\ns\ds\q+\big[C(s)\1n-\1n
D(s)\h\Th(s)\big]^T\1n\G(s,t)\big[C(s)\1n-\1n
D(s)\h\Th(s)\big]\1n+\1n\h\Th(s)^T\1n
R(s,t)\h\Th(s)=0,\\
\ns\ds\h\G_s(s,t)+\h\G(s,t)\big[A(s)-B(s)\h\Th(s)\big]
+\big[A(s)-B(s)\h\Th(s)\big]^T\h\G(s,t)+Q(s,t)\\
\ns\ds\q+\big[C(s)\1n-\1n
D(s)\h\Th(s)\big]^T\1n\G(s,t)\big[C(s)\1n-\1n
D(s)\h\Th(s)\big]\1n+\1n\h\Th(s)^T\1n
R(s,t)\h\Th(s)=0,\q s\in[t,T],\\
\ns\ds\G(T,t)=\h\G(T,t)=G(t).\ea\right.\ee
Then, we see that
$$\h\G(s,t)=\G(s,t),\qq0\le t\le s\le T.$$
Consequently, $\h\Th(\cd)=\Th(\cd)$ and the equation for
$\G(\cd\,,\cd)$ becomes
\bel{G(N-3)}\left\{\2n\ba{ll}
\ns\ds\G_s(s,t)+\G(s,t)\big[A(s)-B(s)\Th(s)\big]+\big[A(s)-B(s)\Th(s)\big]^T
\G(s,t)+Q(s,t)\\
\ns\ds\q+\big[C(s)\1n-\1n D(s)\Th(s)\big]^T\1n\G(s,t)\big[C(s)\1n-\1n D(s)
\Th(s)\big]+\Th(s)^T\1n R(s,t)\Th(s)=0,\q s\in[t,T],\\
\ns\ds\G(T,t)=G(t),\ea\right.\ee
with
$$\Th(s)=\big[R(s,s)+D(s)^T\G(s,s)D(s)\big]^{-1}\big[B(s)^T\G(s,s)+D(s)^T\G(s,s)C(s)\big].$$
This coincides with a special case of the general result found in
\cite{Yong 2012b}.

\ms

\it Case 2. The case of conditional expectation only. \rm Let
$$\left\{\2n\ba{ll}
\ns\ds\bar A(\cd)=\bar C(\cd)=0,\q\bar B(\cd)=\bar D(\cd)=0,\\
\ns\ds Q(s,t)=Q(s),\q\bar Q(s,t)=\bar Q(s),\q R(s,t)=R(s),\q\bar
R(s,t)=\bar R(s),\q G(t)=G,\q\bar G(t)=\bar G.\ea\right.$$
Note that unlike Case 2 in the previous subsection, we allow conditional
expectation terms to appear in the running costs. In this case, one has
$$\left\{\ba{ll}
\ns\ds\h A(\cd)=A(\cd),\q\h B(\cd)=B(\cd),\q\h C(\cd)=C(\cd),\q\h
D(\cd)=D(\cd),\\
\ns\ds\h Q(s,t)=\h Q(s),\q\h R(s,t)=\h R(s),\q\h G(t)=\h
G.\ea\right.$$
Thus, both $\G(s,t)$ and $\h\G(s,t)$ are independent of $t$ and
system (\ref{Riccati-G(s,t)}) becomes
\bel{Riccati-G(s,t)*}\left\{\2n\ba{ll}
\ns\ds\dot\G(s)+\G(s)\big[A(s)-B(s)\h\Th(s)\big]+\big[A(s)
-B(s)\h\Th(s)\big]^T\G(s)+Q(s)\\
\ns\ds\q+\big[C(s)\1n-\1n D(s)\h\Th(s)\big]^T\1n\G(s)
\big[C(s)\1n-\1n D(s)\h\Th(s)\big]\1n+\1n\h\Th(s)^T\1n R(s)\h\Th(s)=0,\\
\ns\ds\dot{\h\G}(s)+\h\G(s)\big[A(s)-
B(s)\h\Th(s)\big]+\big[A(s)-B(s)\h\Th(s)\big]^T\h\G(s)+\h Q(s)\\
\ns\ds\q+\big[C(s)\1n-\1n D(s)\h\Th(s)\big]^T\1n\G(s)\big[
C(s)\1n-\1n D(s)\h\Th(s)\big]\1n+\1n\h\Th(s)^T\1n
\h R(s)\h\Th(s)=0,\qq0\le s\le T,\\
\ns\ds\G(T)=G,\qq\h\G(T)=\h G,\ea\right.\ee
where $\h\Th(\cd)$ is given by the following:
\bel{Th(N-1)}\h\Th(s)=\big[\h R(s)+D(s)^T\G(s)D(s)\big]^{-1}\big[
B(s)^T\h\G(s)+D(s)^T\G(s)C(s)\big],\q s\in[0,T].\ee
Further, if
$$\bar Q(\cd)=0,\qq\bar R(\cd)=0,$$
then we are in the case of conditional expectation at terminal cost
only, and (\ref{Riccati-G(s,t)*}) becomes
\bel{G(s)*}\left\{\2n\ba{ll}
\ns\ds\dot\G(s)+\G(s)\big[A(s)-B(s)\h\Th(s)\big]+\big[A(s)-B(s)\h\Th(s)\big]^T
\G(s)+Q(s)\\
\ns\ds\qq+\big[C(s)\1n-\1n D(s)\h\Th(s)\big]^T\1n\G(s)\big[
C(s)\1n-\1n D(s)\h\Th(s)\big]\1n+\1n\h\Th(s)^T\1n R(s)\h\Th(s)=0,\\
\ns\ds\dot{\h\G}(s)+\h\G(s)\big[A(s)-B(s)\h\Th(s)\big]
+\big[A(s)-B(s)\h\Th(s)\big]^T\h\G(s)+Q(s)\\
\ns\ds\qq+\big[C(s)\1n-\1n
D(s)\h\Th(s)\big]^T\1n\G(s)\big[C(s)\1n-\1n
D(s)\h\Th(s)\big]\1n+\1n\h\Th(s)^T\1n
R(s)\h\Th(s)=0,\q0\le s\le T,\\
\ns\ds\G(T)=G,\qq\h\G(T)=\h G.\ea\right.\ee
The closed-loop equilibrium state process $X^*(\cd)$ satisfies
\bel{}\left\{\2n\ba{ll}
\ns\ds dX^*(s)=\big[A(s)-B(s)\h\Th(s)\big]X^*(s)ds+\big[C(s)-D(s)\h\Th(s)\big]X^*(s)dW(s),\q s\in[0,T],\\
\ns\ds X^*(0)=x,\ea\right.\ee
the closed-loop equilibrium control admits the following representation:
\bel{feedback(N-2)e} u^*(s)=-\h\Th(s)X^*(s),\q s\in[0,T],\ee
and the closed-loop equilibrium value function is given by the following:
\bel{V}\h V(t,x)=\lan\h\G(t)x,x\ran,\qq\forall(t,x)\in\sD.\ee
Clearly, in the case that $\bar G=0$, the problem is reduced to a
classical LQ problem. In this case, $\G(\cd)=\h\G(\cd)$, and our
result recovers the classical one.

\ms

Directly comparing the results of this subsection with those in the
previous subsection, we see that the open-loop and closed-loop
equilibrium solutions are different for Problem (MF-LQ), even for
the above two special cases. The results coincide when the probelm
is reduced to classical LQ problems.

\section{Multi-Person Differential Games}

In this section, we fix a partition
$\D:0=t_0<t_1<\cds<t_{N-1}<t_N=T$ of $[0,T]$, and construct a
closed-loop $\D$-equilibrium strategy for Problem (MF-LQ). To this end, we
introduce an associated $N$-person differential game. For notational
convenience, we label the players by $0,1,2,\cds,N-1$. The $k$-th
player controls the system on $[t_k,t_{k+1})$, with the cost
functional being ``sophisticatedly'' constructed (see below). The
main rule are the following:

\ms

(i) Each player will play optimally based on the assumption that the
later players will play optimally.

\ms

(ii) The $k$-th player will affect the $(k+1)$-th player's action
through her terminal state (which is the initial state of the
$(k+1)$-th player).

\ms

(iii) Although the $k$-th player will not be able to control the
system from $t_{k+1}$ on, she will still ``discount'' the cost
functional in her own way on the interval $[t_{k+1},T]$.

\ms

Let us now look at the $N$-person differential games in details.
In what follows, we denote
$$\left\{\2n\ba{ll}
\ns\ds Q_k(s)=Q(s,t_k),\q R_k(s)=R(s,t_k),\q G_k=G(t_k),\\
\ns\ds\bar Q_k(s)=\bar Q(s,t_k),\q\bar R_k(s)=R(s,t_k),\q\bar
G_k=\bar G(t_k),\\
\ns\ds\h Q_k(s)\1n=\1n Q_k(s)\1n+\1n\bar Q_k(s),\q\h R_k(s)\1n=\1n
R_k(s)\1n+\1n\bar R_k(s),\q\h G_k\1n=\1n G_k\1n+\1n\bar
G_k,\ea\right.\qq k=0,1,\cds,N-1.$$
We will carefully look at Players $(N-1)$, $(N-2)$, and $(N-3)$, who
will have different features. Once we have done that, the situations
for the rest players will be clear and can be treated inductively. We now
begin with Player $(N-1)$. The state equation for this player is the following:
\bel{state(N-1)}\left\{\2n\ba{ll}
\ns\ds dX_{N-1}(s)\1n=\1n\Big\{A(s)X_{N-1}(s)+\bar
A(s)\dbE_{t_{N-1}}[X_{N-1}(s)]+B(s)u_{N-1}(s)+\bar
B(s)\dbE_{t_{N-1}}[u_{N-1}(s)]\Big\}ds\\
\ns\ds\qq\qq\q+\Big\{C(s)X_{N-1}(s)\1n+\1n\bar
C(s)\dbE_{t_{N-1}}[X_{N-1}(s)]\1n+\2n D(s)u_{N-1}(s)\1n+\2n\bar D(s)
\dbE_{t_{N-1}}[u_{N-1}(s)]\Big\}dW(s),\\
\ns\ds\qq\qq\qq\qq\qq\qq\qq\qq\qq\qq\qq\qq\qq\qq\qq\qq s\in[t_{N-1},t_N],\\
\ns\ds X_{N-1}(t_{N-1})=x_{N-1}\in\sX_{t_{N-1}},\ea\right.\ee
and the cost functional is given by
\bel{cost(N-1)}\ba{ll}
\ns\ds J_{N-1}^\D(x_{N-1};u_{N-1}(\cd))\\
\ns\ds=\1n\dbE_{t_{N\1n-\1n1}}\1n\Big\{\1n\int_{t_{N\1n-\1n1}}^{t_N}\2n
\[\1n\lan Q_{N\1n-\1n1}(s)X_{N\1n-\1n1}(s),X_{N\1n-\1n1}(s)\ran\1n
+\1n\lan\bar Q_{N\1n-\1n1}(s)\dbE_{t_{N\1n-\1n1}}[X_{N\1n-\1n1}(s)],
\dbE_{t_{N\1n-\1n1}}[X_{N\1n-\1n1}(s)]\ran\\
\ns\ds\qq\qq+\lan R_{N-1}(s)u_{N-1}(s),u_{N-1}(s)\ran\1n+\1n\lan\bar
R_{N-1}(s)\dbE_{t_{N-1}}[u_{N-1}(s)],\dbE_{t_{N-1}}[u_{N-1}(s)]\ran\1n\]ds\\
\ns\ds\qq\qq+\1n\lan\1n
G_{N-1}X_{N-1}(t_N),X_{N-1}(t_N)\ran\1n+\1n\lan\bar
G_{N-1}\dbE_{t_{N-1}}[X_{N-1}(t_N)],\dbE_{t_{N-1}}[X_{N-1}(t_N)]\ran\1n\Big\}.\ea\ee
Player $(N-1)$ wants to solve the following problem.

\ms

\bf Problem (MF-LQ)$_{N-1}$. \rm For any $x_{N-1}\in\sX_{t_{N-1}}$,
find a $u_{N-1}^*(\cd)\in\sU[t_{N-1},t_N]$ such that
\bel{Problem(N-1)}J^\D_{N-1}(x_{N-1};u_{N-1}^*(\cd))=\inf_{u_{N-1}(\cd)\in\sU[t_{N-1},t_N]}
J^\D_{N-1}(x_{N-1};u_{N-1}(\cd)).\ee

This is a standard LQ problem for MF-SDEs. According to Proposition
3.1, under (H1)--(H2), one has a pair of $\dbS^n$-valued functions
$(P_{N-1}(\cd),\h P_{N-1}(\cd))$ uniquely solve the following
Riccati equations ($s$ is suppressed):
\bel{Riccati P(N-1)}\left\{\2n\ba{ll}
\ns\ds\dot P_{N-1}+P_{N-1}A+A^TP_{N-1}+C^TP_{N-1}C+Q_{N-1}\\
\ns\ds~-\1n(P_{N-1}B\1n+\1n C^T\1n P_{N-1}D)(R_{N-1}\2n+\2n D^T\2n
P_{N-1}D)^{-1}(B^T\1n P_{N-1}\2n+\1n D^T\1n P_{N-1}C)\1n=\1n0,\\
\dot{\h P}_{N-1}+\h P_{N-1}\h A+\h A^T\h P_{N-1}+\h C^TP_{N-1}\h C
+\h Q_{N-1}\\
\ns\ds~-\1n(\h P_{N-1}\h B\1n+\1n\h C^T\1n P_{N-1}\h D)(\h
R_{N-1}\2n+\1n\h D^T\1n P_{N-1}\h D)^{-1}(\h B^T\1n\h
P_{N-1}\2n+\1n\h D^T\1n P_{N-1}
\h C)\1n=\1n0,\q s\1n\in\1n[t_{N-1},t_N),\\
\ns\ds P_{N-1}(t_N)=G_{N-1},\qq\h P_{N-1}(t_N)=\h G_{N-1}.\ea\right.\ee
Define
\bel{Th(N-1)}\left\{\2n\ba{ll}
\ns\ds\Th_{N-1}=(R_{N-1}+D^TP_{N-1}D)^{-1}(B^TP_{N-1}+D^TP_{N-1}C),\\
\ns\ds\h\Th_{N-1}=(\h R_{N-1}+\h D^TP_{N-1}\h D)^{-1}(\h B^T\h
P_{N-1}+\h D^TP_{N-1}\h C),\ea\right.\q s\in[t_{N-1},t_N].\ee
Then the optimal state process $X^*_{N-1}(\cd)\equiv
X^*_{N-1}(\cd\,;t_{N-1},x_{N-1})$ solves the following closed-loop
system:
\bel{closed-loop(N-1)}\left\{\2n\ba{ll}
\ns\ds dX_{N-1}^*=\Big\{\big[A-B\Th_{N-1}\big]X_{N-1}^*+\big[\bar
A+B(\Th_{N-1}-\h\Th_{N-1})-\bar
B\h\Th_{N-1}\big]\dbE_{t_{N-1}}[X_{N-1}^*]\Big\}ds\\
\ns\ds\qq\qq~+\1n\Big\{\1n\big[C\1n-\1n
D\Th_{N-1}\big]X_{N-1}^*\1n+\2n\big[\bar C\1n+\1n
D(\Th_{N-1}\1n-\1n\h\Th_{N-1})\1n-\1n\bar
D\h\Th_{N-1}\big]\dbE_{t_{N-1}}[X_{N-1}^*]\Big\}dW(s),\\
\ns\ds\qq\qq\qq\qq\qq\qq\qq\qq\qq\qq\qq\qq\qq s\in[t_{N-1},t_N],\\
\ns\ds X_{N-1}^*(t_{N-1})=x_{N-1},\ea\right.\ee
and the optimal control $u^*_{N-1}(\cd)\equiv
u^*_{N-1}(\cd\,;t_{N-1},x_{N-1})$ admits the following state
feedback representation:
\bel{feedback(N-1)}
u_{N-1}^*(s)=-\Th_{N-1}(s)X_{N-1}^*(s)+\big[\Th_{N-1}(s)-\h\Th_{N-1}(s)\big]\dbE_{t_{N-1}}
[X_{N-1}^*(s)],\q s\in[t_{N-1},t_N].\ee
Finally,
\bel{inf(N-1)}\ba{ll}
\ns\ds\inf_{u_{N-1}(\cd)\in\sU[t_{N-1},t_N]}J^\D_{N-1}(x_{N-1};u_{N-1}(\cd))
=J^\D_{N-1}(x_{N-1};u^*_{N-1}(\cd))\\
\ns\ds\qq\qq\qq\qq\qq\qq\qq\qq~=\lan\h
P_{N-1}(t_{N-1})x_{N-1},x_{N-1}\ran,\q\forall
x_{N-1}\in\sX_{t_{N-1}}.\ea\ee
This way, Play $(N-1)$ has solved Problem (MF-LQ)$_{N-1}$.

\ms

Next, we consider Player $(N-2)$ whose state equation is the
following:
\bel{state(N-2)}\left\{\2n\ba{ll}
\ns\ds dX_{N-2}(s)\1n=\1n\Big\{A(s)X_{N-2}(s)+\bar
A(s)\dbE_{t_{N-2}}[X_{N-2}(s)]+B(s)u_{N-2}(s)+\bar
B(s)\dbE_{t_{N-2}}[u_{N-2}(s)]\Big\}ds\\
\ns\ds\qq\qq\q+\Big\{\1n C(s)X_{N-2}(s)\1n+\1n\bar
C(s)\dbE_{t_{N-2}}[X_{N-2}(s)]\1n+\2n D(s)u_{N-2}(s)\1n+\2n\bar D(s)
\dbE_{t_{N-2}}[u_{N-2}(s)]\1n\Big\}dW(s),\\
\ns\ds\qq\qq\qq\qq\qq\qq\qq\qq\qq\qq\qq\qq\qq\qq s\in[t_{N-2},t_{N-1}],\\
\ns\ds X(t_{N-2})=x_{N-2}\in\sX_{t_{N-2}}.\ea\right.\ee
Let $X_{N-2}(\cd)=X(\cd\,;t_{N-2},x_{N-2},u_{N-2}(\cd))$ be the corresponding
solution. At $s=t_{N-1}$, Player $(N-1)$ takes over the system, and will use
her optimal control $u_{N-1}^*(\cd)$ of state feedback form
(\ref{feedback(N-1)}) on $[t_{N-1},t_N]$, where $X_{N-1}^*(\cd)$
satisfies closed-loop system (\ref{closed-loop(N-1)}) with initial state
$x_{N-1}=X_{N-2}(t_{N-1})$. Because of this, Player $(N-2)$
considers the following (sophisticated) cost functional (suppressing
$s$)
\bel{cost(N-2)a}\ba{ll}
\ns\ds
J^\D_{N-2}(x_{N-2};u_{N-2}(\cd))\\
\ns\ds=\dbE_{t_{N-2}}\Big\{\int_{t_{N-2}}^{t_{N-1}}\[\lan
Q_{N-2}X_{N-2},X_{N-2}\ran+\lan\bar
Q_{N-2}\dbE_{t_{N-2}}[X_{N-2}],\dbE_{t_{N-2}}[X_{N-2}]\ran\\
\ns\ds\qq\qq\qq\q+\lan R_{N-2}u_{N-2},u_{N-2}\ran+\lan\bar
R_{N-2}\dbE_{t_{N-2}}[u_{N-2}],\dbE_{t_{N-2}}[u_{N-2}]\ran\]ds\\
\ns\ds\qq\qq+\int_{t_{N-1}}^{t_N}\1n\[\lan
Q_{N-2}X^*_{N-1},X^*_{N-1}\ran+\lan\bar
Q_{N-2}\dbE_{t_{N-2}}[X^*_{N-1}],\dbE_{t_{N-2}}[X^*_{N-1}]\ran\\
\ns\ds\qq\qq\qq\q+\lan R_{N-2}u^*_{N-1},u^*_{N-1}\ran+\lan\bar
R_{N-2}\dbE_{t_{N-2}}[u^*_{N-1}],\dbE_{t_{N-2}}[u^*_{N-1}]\ran\]ds\\
\ns\ds\qq\qq\qq\q+\lan
G_{N-2}X^*_{N-1}(t_N),X^*_{N-1}(t_N)\ran\1n+\1n\lan\bar
G_{N-2}\dbE_{t_{N-2}}[X_{N-1}^*(t_N)],\dbE_{t_{N-2}}[X^*_{N-1}(t_N)]\ran
\1n\Big\}.\ea\ee
Note that the appearance of $Q_{N-2},\bar Q_{N-2},R_{N-2},\bar
R_{N-2},G_{N-2},\bar G_{N-2}$ and $\dbE_{t_{N-2}}$ in the running
cost over $[t_{N-1},t_N]$ and in the terminal cost exactly explains
the meaning of ``discounting in her own way'' for Player $(N-2)$
mentioned earlier (see Rule (iii) at the beginning of the section).
We now want to rewrite the above cost functional so
that Player $(N-2)$ will face a standard LQ problem for an MF-SDE on
$[t_{N-2},t_{N-1}]$. To this end, we observe the following: (noting
(\ref{feedback(N-1)}))
\bel{Eu(N-1)}\dbE_{t_{N-2}}[u_{N-1}^*(s)]=-\h\Th_{N-1}(s)\dbE_{t_{N-2}}[X_{N-1}^*(s)],\qq
s\in[t_{N-1},t_N],\ee
\bel{closed-loop(N-1)a}\left\{\2n\ba{ll}
\ns\ds d\big(\dbE_{t_{N-1}}[X_{N-1}^*]\big)=(\h A-\h B\h\Th_{N-1})
\dbE_{t_{N-1}}[X_{N-1}^*]ds,\qq s\in[t_{N-1},t_N],\\
\ns\ds\dbE_{t_{N-1}}[X_{N-1}^*(t_{N-1})]=x_{N-1},\ea\right.\ee
and
\bel{closed-loop(N-1)b}\left\{\2n\ba{ll}
\ns\ds d\big(\dbE_{t_{N-2}}[X_{N-1}^*]\big)=(\h A-\h B\h\Th_{N-1})
\dbE_{t_{N-2}}[X_{N-1}^*]ds,\qq s\in[t_{N-1},t_N],\\
\ns\ds\dbE_{t_{N-2}}[X_{N-1}^*(t_{N-1})]=\dbE_{t_{N-2}}[x_{N-1}].\ea\right.\ee
Thus
$$\ba{ll}
\ns\ds\dbI_{N-2}\equiv\dbE_{t_{N-2}}\Big\{\int_{t_{N-1}}^{t_N}\1n\[\lan
Q_{N-2}X^*_{N-1},X^*_{N-1}\ran+\lan\bar
Q_{N-2}\dbE_{t_{N-2}}[X^*_{N-1}],\dbE_{t_{N-2}}[X^*_{N-1}]\ran\\
\ns\ds\qq\qq\qq\q+\lan R_{N-2}u^*_{N-1},u^*_{N-1}\ran+\lan\bar
R_{N-2}\dbE_{t_{N-2}}[u^*_{N-1}],\dbE_{t_{N-2}}[u^*_{N-1}]\ran\]ds\\
\ns\ds\qq\qq\qq\q+\1n\lan
G_{N-2}X^*_{N-1}(t_N),X^*_{N-1}(t_N)\ran\1n+\1n\lan\bar
G_{N-2}\dbE_{t_{N-2}}[X_{N-1}^*(t_N)],\dbE_{t_{N-2}}[X^*_{N-1}(t_N)]\ran\1n\Big\}\\
\ns\ds\qq=\dbE_{t_{N-2}}\Big\{\int_{t_{N-1}}^{t_N}\1n\[\lan
Q_{N-2}X^*_{N-1},X^*_{N-1}\ran+\lan\bar
Q_{N-2}\dbE_{t_{N-2}}[X^*_{N-1}],\dbE_{t_{N-2}}[X^*_{N-1}]\ran\\
\ns\ds\qq\qq\qq\q+\lan
R_{N-2}\big\{\Th_{N-1}X_{N-1}^*+(\h\Th_{N-1}-\Th_{N-1})\dbE_{t_{N-1}}
[X_{N-1}^*]\big\},\\
\ns\ds\qq\qq\qq\qq\qq\qq\qq\qq
\Th_{N-1}X_{N-1}^*+(\h\Th_{N-1}-\Th_{N-1})\dbE_{t_{N-1}}[X_{N-1}^*]\ran\\
\ns\ds\qq\qq\qq\q+\lan\bar
R_{N-2}\h\Th_{N-1}\dbE_{t_{N-2}}[X^*_{N-1}],\h\Th_{N-1}\dbE_{t_{N-2}}
[X^*_{N-1}]\ran\]ds\\
\ns\ds\qq\qq\qq\q+\1n\lan
G_{N-2}X^*_{N-1}(t_N),X^*_{N-1}(t_N)\ran\1n+\1n\lan\bar
G_{N-2}\dbE_{t_{N-2}}[X_{N-1}^*(t_N)],\dbE_{t_{N-2}}[X^*_{N-1}(t_N)]\ran
\1n\Big\}\\
\ns\ds\qq=\dbE_{t_{N-2}}\Big\{\int_{t_{N-1}}^{t_N}\1n\[\lan
\big(Q_{N-2}+\Th_{N-1}^TR_{N-2}\Th_{N-1}\big)X^*_{N-1},X^*_{N-1}\ran\\
\ns\ds\qq\q+\lan
\big(\h\Th_{N-1}^TR_{N-2}\h\Th_{N-1}-\Th_{N-1}^TR_{N-2}\Th_{N-1}\big)
\dbE_{t_{N-1}}[X_{N-1}^*],\dbE_{t_{N-1}} [X_{N-1}^*]\ran\\
\ns\ds\qq\q+\lan\big(\bar Q_{N-2}\1n+\1n\h\Th_{N-1}^T\bar
R_{N-2}\h\Th_{N-1}\big)
\dbE_{t_{N-2}}[X^*_{N-1}],\dbE_{t_{N-2}}[X^*_{N-1}]\ran\]ds\\
\ns\ds\qq\q+\lan G_{N-2}X^*_{N-1}(t_N),X^*_{N-1}(t_N)\ran+\lan\bar
G_{N-2}\dbE_{t_{N-2}}[X_{N-1}^*(t_N)],\dbE_{t_{N-2}}[X^*_{N-1}N(t_N)]\ran
\1n\Big\}.\ea$$
Note that on the right hand side of the above, both $\dbE_{t_{N-1}}$
and $\dbE_{t_{N-2}}$ appear quadratically. Suggested by Lemma 2.4,
we introduce the following three Lyapunov equations on
$[t_{N-1},t_N]$:
\bel{wt G(N-2)}\left\{\2n\ba{ll}
\ns\ds\dot{\wt\G}_{N-2}\1n+\1n\wt\G_{N-2}(A\1n-\1n
B\Th_{N-1})\1n+\1n(A-B\Th_{N-1})^T\wt\G_{N-2}
\1n+Q_{N-2}\\
\ns\ds\qq+(C-D\Th_{N-1})^T\wt\G_{N-2}(C-D\Th_{N-1})+\Th_{N-1}^TR_{N-2}\Th_{N-1}=0,
\qq s\in[t_{N-1},t_N),\\
\ns\ds\wt\G_{N-2}(t_N)=G_{N-2},\ea\right.\ee
\bel{G(N-2)}\left\{\2n\ba{ll}
\ns\ds\dot\G_{N-2}+\G_{N-2}(\h A-\h B\h\Th_{N-1})+(\h A-\h
B\h\Th_{N-1})^T\G_{N-2}+Q_{N-2}\\
\ns\ds\qq+(\h C-\h D\h\Th_{N-1})^T\wt\G_{N-2}(\h C-\h D\h\Th_{N-1})
+\h\Th_{N-1}^TR_{N-2}\h\Th_{N-1}=0,\qq s\in[t_{N-1},t_N),\\
\ns\ds\G_{N-2}(t_N)=G_{N-2},\ea\right.\ee
and
\bel{bG(N-2)}\left\{\2n\ba{ll}
\ns\ds\dot{\bar\G}_{N-2}+\bar\G_{N-2}(\h A-\h B\h\Th_{N-1})+(\h A
-\h B\h\Th_{N-1})^T\bar\G_{N-2}+\bar Q_{N-2}+\h\Th_{N-1}^T\bar R_{N-2}
\h\Th_{N-1}=0,\\
\ns\ds\qq\qq\qq\qq\qq\qq\qq\qq\qq\qq\qq\qq s\in[t_{N-1},t_N),\\
\ns\ds\bar\G_{N-2}(t_N)=\bar G_{N-2}.\ea\right.\ee
Then, using Lemma 2.4, the sophisticated cost functional of Player
$(N-2)$ can be written as
\bel{cost(N-2)b}\ba{ll}
\ns\ds
J_{N-2}^\D(x_{N-2};u_{N-2}(\cd))\\
\ns\ds=\dbE_{t_{N-2}}\Big\{\int_{t_{N-2}}^{t_{N-1}}\[\lan
Q_{N-2}X_{N-2},X_{N-2}\ran+\lan\bar
Q_{N-2}\dbE_{t_{N-2}}[X_{N-2}],\dbE_{t_{N-2}}[X_{N-2}]\ran\\
\ns\ds\qq\qq+\lan R_{N-2}u_{N-2},u_{N-2}\ran+\lan\bar
R_{N-2}\dbE_{t_{N-2}}[u_{N-2}],\dbE_{t_{N-2}}[u_{N-2}]\ran\]ds\\
\ns\ds\qq\qq+\lan\G_{N-2}(t_{N-1})X_{N-2}(t_{N-1}),X_{N-2}(t_{N-1})\ran\\
\ns\ds\qq\qq+\lan\bar\G_{N-2}(t_{N-1})\dbE_{t_{N-2}}[X_{N-2}(t_{N-1})],
\dbE_{t_{N-2}}[X_{N-2}(t_{N-1})]\ran\Big\}.\ea \ee
This is a standard cost functional of an LQ problem for MF-SDEs on
the interval $[t_{N-2},t_{N-1}]$ now. We see that Lemma 2.4 plays an
interesting role here. \ms

Let us make two observations. First, we note that, in general,
$$\ba{ll}
\ns\ds\lan\h P_{N-1}(t_{N-1})X_{N-2}(t_{N-1}),X_{N-2}(t_{N-1})\ran
=J^\D_{N-1}(X_{N-2}(t_{N-1});u_{N-1}^*(\cd))\\
\ns\ds=\dbE_{t_{N-1}}\Big\{\int_{t_{N-1}}^{t_N}\1n\[\lan
Q_{N-1}X^*_{N-1},X^*_{N-1}\ran+\lan\bar
Q_{N-1}\dbE_{t_{N-1}}[X^*_{N-1}],\dbE_{t_{N-1}}[X^*_{N-1}]\ran\\
\ns\ds\qq\qq\qq\q+\lan R_{N-1}u^*_{N-1},u^*_{N-1}\ran+\lan\bar
R_{N-1}\dbE_{t_{N-1}}[u^*_{N-1}],\dbE_{t_{N-1}}[u^*_{N-1}]\ran\]ds\\
\ns\ds\qq\qq\qq\q+\1n\lan
G_{N-1}X^*_{N-1}(t_N),X^*_{N-1}(t_N)\ran\1n+\1n\lan\bar
G_{N-1}\dbE_{t_{N-1}}[X_{N-1}^*(t_N)],
\dbE_{t_{N-1}}[X^*_{N-1}(t_N)]\ran\1n\Big\}\\
\ns\ds\ne\dbI_{N-2}=\dbE_{t_{N-2}}\Big\{\int_{t_{N-1}}^{t_N}\1n\[\lan
Q_{N-2}X^*_{N-1},X^*_{N-1}\ran+\lan\bar
Q_{N-2}\dbE_{t_{N-2}}[X^*_{N-1}],\dbE_{t_{N-2}}[X^*_{N-1}]\ran\\
\ns\ds\qq\qq\qq\q+\lan R_{N-2}u^*_{N-1},u^*_{N-1}\ran+\lan\bar
R_{N-2}\dbE_{t_{N-2}}[u^*_{N-1}],\dbE_{t_{N-2}}[u^*_{N-1}]\ran\]ds\\
\ns\ds\qq\qq\qq\q+\1n\lan
G_{N-2}X^*_{N-1}(t_N),X^*_{N-1}(t_N)\ran\1n+\1n\lan\bar
G_{N-2}\dbE_{t_{N-2}}[X_{N-1}^*(t_N)],
\dbE_{t_{N-2}}[X^*_{N-1}(t_N)]\ran\1n\Big\}\\
\ns\ds=\lan\G_{N-2}(t_{N-1})X_{N-2}(t_{N-1}),X_{N-2}(t_{N-1})\ran\\
\ns\ds\qq\qq+\lan\bar\G_{N-2}(t_{N-1})\dbE_{t_{N-2}}[X_{N-2}(t_{N-1})],
\dbE_{t_{N-2}}[X_{N-2}(t_{N-1})]\ran\Big\},\ea$$
since the weighting
matrices are possibly different, and conditional expectation terms
are different.

\ms

Second, if we let
$$\left\{\2n\ba{ll}
\ns\ds\cA(s)=A(s)-B(s)\Th_{N-1}(s),\qq\cC(s)=C(s)-D(s)\Th_{N-1}(s),\\
\ns\ds\bar\cA(s)=\bar
A(s)+B(s)\big[\Th_{N-1}(s)-\h\Th_{N-1}(s)\big]-\bar
B(s)\h\Th_{N-1}(s),\\
\ns\ds\bar\cC(s)=\bar
C(s)+D(s)\big[\Th_{N-1}(s)-\h\Th_{N-1}(s)\big]-\bar
D(s)\h\Th_{N-1}(s),\\
\ns\ds\cQ(s)=Q_{N-2}(s)+\Th_{N-1}(s)^TR_{N-2}(s)\Th_{N-1}(s),\\
\ns\ds\wt\cQ(s)=\h\Th_{N-1}^TR_{N-2}(s)\h\Th_{N-1}(s)-\Th_{N-1}(s)^TR_{N-2}(s)
\Th_{N-1}(s),\\
\ns\ds\bar\cQ(s)=\bar Q_{N-2}(s)+\h\Th_{N-1}(s)^T\bar
R_{N-2}(s)\h\Th_{N-1}(s),\ea\right.$$
then, under (H1)--(H2), one has
$$\left\{\2n\ba{ll}
\ns\ds\cQ(s)\ge0,\\
\ns\ds\cQ(s)+\wt\cQ(s)=Q_{N-2}(s)+\h\Th_{N-1}(s)^TR_{N-2}(s)\h\Th_{N-1}(s)\ge0,\\
\ns\ds\cQ(s)+\wt\cQ(s)+\bar\cQ(s)=\h Q_{N-2}(s)+\h\Th_{N-1}(s)^T\h
R_{N-2}(s)\h\Th_{N-2}(s)\1n\ge\1n0.\ea\right.$$
Thus, by Lemma 2.4, we have
\bel{G(N-2)>0}\G_{N-2}(t_{N-1}),~\G_{N-2}(t_{N-1})+\bar\G_{N-2}(t_{N-1})\ge0.\ee
Having the above, we now pose the following problem for Player
$(N-2)$.

\ms

\bf Problem (MF-LQ)$_{N-2}$. \rm For $x_{N-2}\in\sX_{t_{N-2}}$, find
a $u^*_{N-2}(\cd)\in\sU[t_{N-2},t_{N-1}]$ such that
\bel{Problem(N-2)}J^\D_{N-2}(x_{N-2};u_{N-2}^*(\cd))=\inf_{u_{N-2}(\cd)\in\sU[t_{N-2},
t_{N-1}]}J^\D_{N-2}(x_{N-2};u_{N-2}(\cd)).\ee

\ms

Similar to the case of Player $(N-1)$ above, under (H1)--(H2), there
exists a pair $(P_{N-2}(\cd),\h P_{N-2}(\cd))$ uniquely solve the
following Riccati equation system:
\bel{Riccati P(N-2)}\left\{\2n\ba{ll}
\ns\ds\dot P_{N-2}+P_{N-2}A+A^TP_{N-2}+C^TP_{N-2}C+Q_{N-2}\\
\ns\ds\q-\1n(P_{N-2}B\1n+\1n C^T\1n
P_{N-2}D)(R_{N-2}\1n+\1n D^T\1n P_{N-2}D)^{-1}\1n
(B^T\1n P_{N-2}\1n+\1n D^T\1n P_{N-2}C)\1n=\1n0,\\
\ns\ds\dot{\h P}_{N-2}+\h P_{N-2}\h A+\h A^T\h P_{N-2}+\h
C^TP_{N-2}\h C+\h Q_{N-2}\\
\ns\ds\q-\1n(\h P_{N-2}\h B\1n+\1n\h C^T\1n P_{N-2}\h D)(\h
R_{N-2}\1n+\2n\h D^T\1n P_{N-2}\h D)^{-1}\1n
(\h B^T\1n\h P_{N-2}\1n+\2n\h D^T\1n P_{N-2}\h C)\1n=\1n0,\q s\1n\in
\1n[t_{N-2},t_{N-1}),\\
\ns\ds P_{N-2}(t_{N-1})=\G_{N-2}(t_{N-1}),\qq
\h P_{N-2}(t_{N-1})=\G_{N-2}(t_{N-1})+\bar\G_{N-2}(t_{N-1}).\ea\right.\ee
Define
\bel{Th_(N-2)}\left\{\2n\ba{ll}
\ns\ds\Th_{N-2}=(R_{N-2}+D^TP_{N-2}D)^{-1}(B^TP_{N-2}+D^TP_{N-2}C),\\
\ns\ds\h\Th_{N-2}=(\h R_{N-2}+\h D^TP_{N-2}\h D)^{-1}(\h B^T\h
P_{N-2}+\h D^TP_{N-2}\h C),\ea\right.\q s\in[t_{N-2},t_{N-1}].\ee
Then the optimal state process $X_{N-2}^*(\cd)$ solves the following
closed-loop system:
\bel{closed-loop(N-2)}\left\{\2n\ba{ll}
\ns\ds dX_{N-2}^*=\Big\{(A-B\Th_{N-2})X_{N-2}^*+\big[\bar
A+B(\Th_{N-2}-\h\Th_{N-2})-\bar
B\h\Th_{N-2}\big]\dbE_{t_{N-2}}[X_{N-2}^*]\Big\}ds\\
\ns\ds\qq\qq\q+\Big\{(C\1n-\1n D\Th_{N-2})X_{N-2}^*\2n+\2n\big[\bar
C\1n+\1n D(\Th_{N-2}\1n-\1n\h\Th_{N-2})\1n-\1n\bar
D\h\Th_{N-2}\big]\dbE_{t_{N-2}}[X_{N-2}^*]\Big\}dW(s),\\
\ns\ds\qq\qq\qq\qq\qq\qq\qq\qq\qq\qq\qq\qq\qq s\in[t_{N-2},t_{N-1}],\\
\ns\ds X_{N-2}^*(t_{N-2})=x_{N-2},\ea\right.\ee
and the optimal control $u^*_{N-2}(\cd)\equiv
u^*_{N-2}(\cd\,;t_{N-2},x_{N-2})$ admits the following state
feedback representation:
\bel{feedback(N-2)}
u_{N-2}^*(s)\1n=\1n-\Th_{N-2}(s)X_{N-2}^*(s)\1n+\1n\big[\Th_{N-2}(s)
\1n-\1n\h\Th_{N-2}(s)\big]\dbE_{t_{N-2}}[X_{N-2}^*(s)],\q
s\in[t_{N-2},t_{N-1}].\ee
Finally, the following holds:
\bel{inf(N-2)}\ba{ll}
\ns\ds\inf_{u_{N-2}(\cd)\in\sU[t_{N-2},\,t_{N-1}]}J^\D_{N-2}(x_{N-2};
u_{N-2}(\cd))=J^\D_{N-2}(x_{N-2};u^*_{N-2}(\cd))\\
\ns\ds\qq\qq\qq\qq\qq\qq\qq\qq=\lan\h
P_{N-2}(t_{N-2})x_{N-2},x_{N-2}\ran,\qq\forall
x_{N-2}\in\sX_{t_{N-2}}.\ea\ee
Thus, Player $(N-2)$ has solved Problem (FM-LQ)$_{N-2}$. Different
from Player $(N-1)$, we need to solve three Lyapunov equations in
order to transform the complicated-looking sophisticated cost
functional into a standard cost functional for a controlled linear
MF-SDE. Thus, the situation of Player $(N-2)$ is interestingly
different from that of Player $(N-1)$.

\ms

Next, we consider Player $(N-3)$ on $[t_{N-3},t_{N-2}]$. For any
$x_{N-3}\in\sX_{t_{N-3}}$, consider
\bel{state(N-3)}\left\{\2n\ba{ll}
\ns\ds dX_{N-3}(s)\1n=\1n\Big\{A(s)X_{N-3}(s)+\bar
A(s)\dbE_{t_{N-3}}[X_{N-3}(s)]+B(s)u_{N-3}(s)+\bar
B(s)\dbE_{t_{N-3}}[u_{N-3}(s)]\Big\}ds\\
\ns\ds\qq\qq\q+\Big\{\1n C(s)X_{N\1n-\1n3}(s)\1n+\1n\bar
C(s)\dbE_{t_{N\1n-\1n3}}[X_{N\1n-\1n3}(s)]\1n+\2n
D(s)u_{N-3}(s)\1n+\2n\bar D(s)
\dbE_{t_{N-3}}[u_{N-3}(s)]\Big\}dW(s),\\
\ns\ds\qq\qq\qq\qq\qq\qq\qq\qq\qq\qq\qq\qq\qq s\in[t_{N-3},t_{N-2}],\\
\ns\ds X_{N-3}(t_{N-3})=x_{N-3}\in\sX_{t_{N-3}}.\ea\right.\ee
Due to the fact that Players $(N-2)$ and $(N-1)$ will play optimally
on $[t_{N-2},t_{N-1})$ and $[t_{N-1},t_N]$, respectively, the
resulting state process, denoted by $X^\D(\cd)$, on $[t_{N-2},t_N]$
will satisfy the following:
\bel{closed-loop(N-2)a}\left\{\2n\ba{ll}
\ns\ds dX^\D=\Big\{(A-B\Th^\D)X^\D+\big[\bar
A+B(\Th^\D-\h\Th^\D)-\bar
B\h\Th^\D\big]\dbE_{\rho^\D(s)}[X^\D]\Big\}ds\\
\ns\ds\qq\qq+\Big\{(C-D\Th^\D)X^\D+\big[\bar
C+D(\Th^\D-\h\Th^\D)-\bar
D\h\Th^\D\big]\dbE_{\rho^\D(s)}[X^\D]\Big\}dW(s),\q s\in[t_{N-2},t_N],\\
\ns\ds X^\D(t_{N-2})=X_{N-3}(t_{N-2}),\ea\right.\ee
where
$$\left\{\2n\ba{ll}
\ns\ds\Th^\D(s)=\Th_{N-2}(s)I_{[t_{N-2},t_{N-1})}(s)+\Th_{N-1}(s)I_{[t_{N-1},
t_N]}(s),\\
\ns\ds\h\Th^\D(s)=\h\Th_{N-2}(s)I_{[t_{N-2},t_{N-1})}(s)+\h\Th_{N-1}(s)
I_{[t_{N-1},t_N]}(s),\ea\right.\q s\in(t_{N-2},t_N],$$
and
$$\rho^\D(s)=t_{N-2}I_{[t_{N-2},t_{N-1})}(s)+t_{N-1}I_{[t_{N-1},t_N]}(s),
\qq s\in[t_{N-2},t_N].$$
The corresponding control on $[t_{N-2},t_N]$ takes the following
form:
\bel{feedback(N-2)e}\ba{ll}
\ns\ds u^\D(s)=-\Th^\D(s)X^\D(s)+\big[\Th^\D(s)-\h\Th^\D(s)\big]
\dbE_{\rho^\D(s)}[X^\D(s)]\\
\ns\ds\qq~=\left\{\2n\ba{ll}
\ns\ds-\Th_{N-1}(s)X^\D(s)+\big[\Th_{N-1}(s)-\h\Th_{N-1}(s)\big]
\dbE_{t_{N-1}}[X^\D(s)],\q s\in[t_{N-1},t_N],\\
\ns\ds-\Th_{N-2}(s)X^\D(s)+\big[\Th_{N-2}(s)-\h\Th_{N-2}(s)\big]
\dbE_{t_{N-2}}[X^\D(s)],\q s\in[t_{N-2},t_{N-1}).\ea\right.\ea\ee
Similar to Player $(N-2)$, Player $(N-3)$ considers the following
(sophisticated) cost functional
\bel{cost(N-3)}\ba{ll}
\ns\ds
J^\D_{N-3}(x_{N-3};u_{N-3}(\cd))\\
\ns\ds=\dbE_{t_{N-3}}\Big\{\int_{t_{N-3}}^{t_{N-2}}\[\lan
Q_{N-3}X_{N-3},X_{N-3}\ran+\lan\bar
Q_{N-3}\dbE_{t_{N-3}}[X_{N-3}],\dbE_{t_{N-3}}[X_{N-3}]\ran\\
\ns\ds\qq\qq\qq\q+\lan R_{N-3}u_{N-3},u_{N-3}\ran+\lan\bar
R_{N-3}\dbE_{t_{N-3}}[u_{N-3}],\dbE_{t_{N-3}}[u_{N-3}]\ran\]ds\\
\ns\ds\qq\q+\int_{t_{N-2}}^{t_N}\1n\[\lan
Q_{N-3}X^\D,X^\D\ran+\lan\bar
Q_{N-3}\dbE_{t_{N-3}}[X^\D],\dbE_{t_{N-3}}[X^\D]\ran\\
\ns\ds\qq\qq\qq\q+\lan R_{N-3}u^\D,u^\D\ran+\lan\bar
R_{N-3}\dbE_{t_{N-3}}[u^\D],\dbE_{t_{N-3}}[u^\D]\ran\]ds\\
\ns\ds\qq\qq\qq\q+\lan G_{N-3}X^\D(t_N),X^\D(t_N)\ran+\lan\bar
G_{N-3}\dbE_{t_{N-3}}[X^\D(t_N)],\dbE_{t_{N-3}}[X^\D(t_N)]\ran\1n
\Big\}\\
\ns\ds\equiv\1n\dbE_{t_{N-3}}\Big\{\1n\int_{t_{N-3}}^{t_{N-2}}\[\lan
Q_{N-3}X_{N-3},X_{N-3}\ran+\lan\bar
Q_{N-3}\dbE_{t_{N-3}}[X_{N-3}],\dbE_{t_{N-3}}[X_{N-3}]\ran\\
\ns\ds\qq\qq\qq+\1n\lan R_{N-3}u_{N-3},u_{N-3}\ran\1n+\1n\lan\bar
R_{N-3}\dbE_{t_{N-3}}[u_{N-3}],\dbE_{t_{N-3}}[u_{N-3}]\ran\]ds\1n+\1n\dbI_{N-3},\ea\ee
where
$$\ba{ll}
\ns\ds\dbI_{N-3}\equiv\dbE_{t_{N-3}}\Big\{\int_{t_{N-2}}^{t_N}\1n\[\lan
Q_{N-3}X^\D,X^\D\ran+\lan\bar
Q_{N-3}\dbE_{t_{N-3}}[X^\D],\dbE_{t_{N-3}}[X^\D]\ran\\
\ns\ds\qq\qq\qq\q+\lan R_{N-3}u^\D,u^\D\ran+\lan\bar
R_{N-3}\dbE_{t_{N-3}}[u^\D],\dbE_{t_{N-3}}[u^\D]\ran\]ds\\
\ns\ds\qq\qq\qq\q+\1n\lan
G_{N-3}X^\D(t_N),X^\D(t_N)\ran\1n+\1n\lan\bar
G_{N-3}\dbE_{t_{N-3}}[X^\D(t_N)],\dbE_{t_{N-3}}[X^\D(t_N)]\ran\1n\Big\}\\
\ns\ds\qq=\dbE_{t_{N-3}}\Big\{\int_{t_{N-2}}^{t_N}\1n\[\lan
Q_{N-3}X^\D,X^\D\ran+\lan\bar
Q_{N-3}\dbE_{t_{N-3}}[X^\D],\dbE_{t_{N-3}}[X^\D]\ran\\
\ns\ds\qq\qq\qq\q+\lan\1n R_{N-3}\big\{\Th^\D
X^\D\2n+\1n(\h\Th^\D\2n-\1n\Th^\D)\dbE_{\rho^\D(s)}[X^\D]\big\},
\Th^\D X^\D\2n+\1n(\h\Th^\D\2n-\1n\Th^\D)\dbE_{\rho^\D(s)}[X^\D]\ran\\
\ns\ds\qq\qq\qq\q+\lan\bar
R_{N-3}\h\Th^\D\dbE_{t_{N-3}}[X^\D],\h\Th^\D\dbE_{t_{N-3}}
[X^\D]\ran\]ds\\
\ns\ds\qq\qq\qq\q+\lan G_{N-3}X^\D(t_N),X^\D(t_N)\ran+\lan\bar
G_{N-3}\dbE_{t_{N-3}}[X^\D(t_N)],\dbE_{t_{N-3}}[X^\D(t_N)]\ran\1n\Big\}\\
%
%
%
\ns\ds\qq\qq\qq\q+\lan\big[\bar Q_{N-3}\1n+\1n(\h\Th^\D)^T\bar
R_{N-3}\h\Th^\D\big]
\dbE_{t_{N-3}}[X^\D],\dbE_{t_{N-3}}[X^\D]\ran\]ds\\
\ns\ds\qq\qq\qq\q+\lan G_{N-3}X^\D(t_N),X^\D(t_N)\ran+\lan\bar
G_{N-3}\dbE_{t_{N-3}}[X^\D(t_N)],\dbE_{t_{N-3}}[X^\D(t_N)]\ran
\1n\Big\}\\
\ns\ds\qq=\dbE_{t_{N-3}}\Big\{\int_{t_{N-2}}^{t_{N-1}}\1n\[\lan
\big[Q_{N-3}+(\Th^\D)^TR_{N-3}\Th^\D\big]X^\D,X^\D\ran\\
\ns\ds\qq\qq\qq\q+\lan\big[(\h\Th^\D)^TR_{N-3}\h\Th^\D\1n
-\1n(\Th^\D)^TR_{N-3}\Th^\D\big]
\dbE_{t_{N-2}}[X^\D],\dbE_{t_{N-2}}[X^\D]\ran\\
\ns\ds\qq\qq\qq\q+\lan\big[\bar Q_{N-3}\1n+\1n(\h\Th^\D)^T\bar
R_{N-3}\h\Th^\D\big]
\dbE_{t_{N-3}}[X^\D],\dbE_{t_{N-3}}[X^\D]\ran\]ds\Big\}\\
\ns\ds\qq\qq\qq\q+\dbE_{t_{N-3}}\Big\{\int_{t_{N-1}}^{t_N}\1n\[\lan
\big[Q_{N-3}+(\Th^\D)^TR_{N-3}\Th^\D\big]X^\D,X^\D\ran\\
\ns\ds\qq\qq\qq\q+\lan\big[(\h\Th^\D)^TR_{N-3}\h\Th^\D\1n
-\1n(\Th^\D)^TR_{N-3}\Th^\D\big]
\dbE_{t_{N-1}}[X^\D],\dbE_{t_{N-1}}[X^\D]\ran\\
\ns\ds\qq\qq\qq\q+\lan\big[\bar Q_{N-3}\1n+\1n(\h\Th^\D)^T\bar
R_{N-3}\h\Th^\D\big]
\dbE_{t_{N-3}}[X^\D],\dbE_{t_{N-3}}[X^\D]\ran\]ds\\
\ns\ds\qq\qq\qq\q+\lan G_{N-3}X^\D(t_N),X^\D(t_N)\ran+\lan\bar
G_{N-3}\dbE_{t_{N-3}}[X^\D(t_N)],\dbE_{t_{N-3}}[X^\D(t_N)]\ran
\1n\Big\}.\ea$$
To take care of the above, we apply Lemma 2.4, by introducing the
following three Lyapunov equations on $[t_{N-2},t_N]$:
\bel{wt G(N-3)}\left\{\2n\ba{ll}
\ns\ds\dot{\wt\G}_{N-3}+\wt\G_{N-3}(A-B\Th^\D)+(A-B\Th^\D)^T\wt\G_{N-3}+Q_{N-3}\\
\ns\ds\qq+(C-D\Th^\D)^T\wt\G_{N-3}(C-D\Th^\D)+(\Th^\D)^TR_{N-3}\Th^\D=0,\\
\ns\ds\qq\qq\qq\qq\qq\qq\qq\qq
s\in(t_{N-2},t_{N-1})\cup(t_{N-1},t_N),\\
\ns\ds\wt\G_{N-3}(t_N)=G_{N-3},\q\wt\G_{N-3}(t_{N-1})=\G_{N-3}(t_{N-1}),\ea\right.\ee
\bel{G(N-3)}\left\{\2n\ba{ll}
\ns\ds\dot\G_{N-3}+\G_{N-3}(\h A-\h B\h\Th^\D)+(\h A-\h B\h\Th^\D)^T
\G_{N-3}+Q_{N-3}\\
\ns\ds\qq+(\h C-\h D\h\Th^\D)^T\wt\G_{N-3}(\h C-\h D\h\Th^\D)
+(\h\Th^\D)^TR_{N-3}\h\Th^\D=0,
\qq s\in[t_{N-2},t_N),\\
\ns\ds\G_{N-3}(t_N)=G_{N-3},\ea\right.\ee
and
\bel{bG(N-3)}\left\{\2n\ba{ll}
\ns\ds\dot{\bar\G}_{N-3}+\bar\G_{N-3}(\h A-\h B\h\Th^\D)+(\h A-\h B)
\h\Th^\D)^T\bar\G_{N-3}+\bar Q_{N-3}+(\h\Th^\D)^T\h R_{N-3}\h\Th^\D=0,\\
\ns\ds\qq\qq\qq\qq\qq\qq\qq\qq\qq\qq\qq\qq\qq s\in[t_{N-2},t_N),\\
\ns\ds\bar\G_{N-3}(t_N)=\bar G_{N-3}.\ea\right.\ee
Note that different from the situation of Player $(N-2)$, the
Lyapunov equation for $\wt\G_{N-3}(\cd)$ is formulated on two
intervals and there might have a jump at $t_{N-1}$ for
$\wt\G_{N-3}(\cd)$. Then
$$\ba{ll}
\ns\ds\dbE_{t_{N-3}}\Big\{\int_{t_{N-1}}^{t_N}\1n\[\lan
\big[Q_{N-3}+(\Th^\D)^TR_{N-3}\Th^\D\big]X^\D,X^\D\ran\\
\ns\ds\qq+\lan\big[(\h\Th^\D)^TR_{N-3}\h\Th^\D\1n
-\1n(\Th^\D)^TR_{N-3}\Th^\D\big]
\dbE_{t_{N-1}}[X^\D],\dbE_{t_{N-1}}[X^\D]\ran\\
\ns\ds\qq+\lan\big[\bar Q_{N-3}\1n+\1n(\h\Th^\D)^T\bar
R_{N-3}\h\Th^\D\big]
\dbE_{t_{N-3}}[X^\D],\dbE_{t_{N-3}}[X^\D]\ran\]ds\\
\ns\ds\qq+\lan G_{N-3}X^\D(t_N),X^\D(t_N)\ran+\lan\bar
G_{N-3}\dbE_{t_{N-3}}[X^\D(t_N)],\dbE_{t_{N-3}}[X^\D(t_N)]\ran
\1n\Big\}\\
\ns\ds=\dbE_{t_{N-3}}\Big\{\lan\G_{N-3}(t_{N-1})X^\D(t_{N-1}),X^\D(t_{N-1})\ran\\
\ns\ds\qq\qq+\lan\bar\G_{N-3}(t_{N-1})\dbE_{t_{N-3}}[X^\D(t_{N-1})],
\dbE_{t_{N-3}}[X^\D(t_{N-1})]\ran\Big\}.\ea$$
In the same manner, we have
$$\ba{ll}
\ns\ds\dbI_{N-3}=\dbE_{t_{N-3}}\Big\{\int_{t_{N-2}}^{t_{N-1}}\1n\[\lan
\big[Q_{N-3}+(\Th^\D)^TR_{N-3}\Th^\D\big]X^\D,X^\D\ran\\
\ns\ds\qq\qq+\lan\big[(\h\Th^\D)^TR_{N-3}\h\Th^\D\1n
-\1n(\Th^\D)^TR_{N-3}\Th^\D\big]
\dbE_{t_{N-2}}[X^\D],\dbE_{t_{N-2}}[X^\D]\ran\\
\ns\ds\qq\qq+\lan\big[\bar Q_{N-3}\1n+\1n(\h\Th^\D)^T\bar
R_{N-3}\h\Th^\D\big]
\dbE_{t_{N-3}}[X^\D],\dbE_{t_{N-3}}[X^\D]\ran\]ds\\
\ns\ds\qq\qq+\lan\G_{N-3}(t_{N-1})X^\D(t_{N-1}),X^\D(t_{N-1})\ran+\lan\bar\G_{N-3}(t_{N-1})\dbE_{t_{N-3}}[X^\D(t_{N-1})],
\dbE_{t_{N-3}}[X^\D(t_{N-1})]\ran\Big\}\\
\ns\ds\qq\q=\1n\dbE_{t_{N-3}}\1n\Big\{\1n\lan\G_{N-3}(t_{N-2})X^\D(t_{N-2}),X^\D(t_{N-2})\ran
\1n+\1n\lan\bar\G_{N-3}(t_{N-2})\dbE_{t_{N-3}}[X^\D(t_{N-2})],
\dbE_{t_{N-3}}[X^\D(t_{N-2})]\ran\1n\Big\}.\ea$$
Hence, the sophisticated cost functional for Player $(N-3)$ can be
written as
\bel{cost(N-3)b}\ba{ll}
\ns\ds
J^\D_{N-3}(x_{N-3};u_{N-3}(\cd))\\
\ns\ds=\dbE_{t_{N-3}}\Big\{\int_{t_{N-3}}^{t_{N-2}}\[\lan
Q_{N-3}X_{N-3},X_{N-3}\ran+\lan\bar
Q_{N-3}\dbE_{t_{N-3}}[X_{N-3}],\dbE_{t_{N-3}}[X_{N-3}]\ran\\
\ns\ds\qq\qq+\lan R_{N-3}u_{N-3},u_{N-3}\ran+\lan\bar
R_{N-3}\dbE_{t_{N-3}}[u_{N-3}],\dbE_{t_{N-3}}[u_{N-3}]\ran\]ds\\
\ns\ds\qq\qq+\lan\G_{N-3}(t_{N-2})X_{N-3}(t_{N-2}),X_{N-3}(t_{N-2})\ran\\
\ns\ds\qq\qq+\lan\bar\G_{N-3}(t_{N-2})\dbE_{t_{N-3}}[X_{N-3}(t_{N-2})],
\dbE_{t_{N-3}}[X_{N-3}(t_{N-2})]\ran\Big\}.\ea \ee
Note that in solving (\ref{wt G(N-3)})--(\ref{bG(N-3)}), we first
solve (\ref{wt G(N-3)}) on $[t_{N-1},t_N]$, then solve
(\ref{G(N-3)}) on $[t_{N-1},t_N]$. After that, solve (\ref{wt
G(N-3)}) on $[t_{N-2},t_{N-1}]$, and (\ref{G(N-3)}) on
$[t_{N-2},t_{N-1}]$. Solving (\ref{bG(N-3)}) is more directly.

Similar to the situation for Player $(N-2)$ (see the derivation of
(\ref{G(N-2)>0})), we first have
$$\G_{N-3}(t_{N-1}),~\G_{N-3}(t_{N-1})+\bar\G_{N-3}(t_{N-1})\ge0.$$
Then applying the same argument, we will have
$$\G_{N-3}(t_{N-2}),~\G_{N-3}(t_{N-2})+\bar\G_{N-3}(t_{N-2})\ge0.$$
Therefore, we can pose the following problem for Player $(N-3)$.

\ms

\bf Problem (MF-LQ)$_{N-3}$. \rm For $x_{N-3}\in\sX_{t_{N-3}}$, find
a $u^*_{N-3}(\cd)\in\sU[t_{N-3},t_{N-2}]$ such that
\bel{inf(N-3)}J^\D_{N-3}(x_{N-3};u_{N-3}^*(\cd))=\inf_{u_{N-3}(\cd)\in\sU[t_{N-3},
t_{N-2}]}J^\D_{N-3}(x_{N-3};u_{N-3}(\cd)).\ee

\ms

This is a standard LQ problem for MF-SDEs now. Under (H1)--(H2),
there exist a pair $(P_{N-3}(\cd),\h P_{N-3}(\cd))$ uniquely
solve the following Riccati equation system:
\bel{Riccati P(N-3)}\left\{\2n\ba{ll}
\ns\ds\dot P_{N-3}+P_{N-3}A+A^TP_{N-3}+C^TP_{N-3}C+Q_{N-3}\\
\ns\ds\q-(P_{N-3}B\1n+\1n C^T\1n
P_{N-3}D)(R_{N-3}\1n+\1n D^T\1n P_{N-3}D)^{-1}\1n
(B^T\1n P_{N-3}\1n+\1n D^T\1n P_{N-3}C)\1n=\1n0,\q s\1n\in\1n[t_{N-3},
t_{N-2}),\\
\dot{\h P}_{N-3}+\h P_{N-3}\h A+\h A^T\h P_{N-3}+\h C^TP_{N-3}\h C
+\h Q_{N-3}\\
\ns\ds\q-(\h P_{N-3}\h B\1n+\1n\h C^T\1n P_{N-3}\h D)(\h
R_{N-3}\1n+\1n\h D^T\1n P_{N-3}\h D)^{-1}\1n
(\h B^T\1n\h P_{N-3}\1n+\1n\h D^T\1n P_{N-3}\h C)\1n=\1n0,\q s\1n\in\1n
[t_{N-3},t_{N-2}),\\
\ns\ds P_{N-3}(t_{N-2})=\G_{N-3}(t_{N-2}+0),\qq\h
P_{N-3}(t_{N-2})=\G_{N-3}(t_{N-2})+\bar\G_{N-3}(t_{N-2}).\ea\right.\ee
Define
\bel{Th_(N-3)}\left\{\2n\ba{ll}
\ns\ds\Th_{N-3}=(R_{N-3}+D^TP_{N-3}D)^{-1}(B^TP_{N-3}+D^TP_{N-3}C),\\
\ns\ds\h\Th_{N-3}=(\h R_{N-3}+\h D^TP_{N-3}\h D)^{-1}(\h B^T\h
P_{N-3}+\h D^TP_{N-3}\h C),\ea\right.\q s\in[t_{N-3},t_{N-2}].\ee
Then the optimal state process $X_{N-3}^*(\cd)$ solves the following
closed-loop system:
\bel{closed-loop(N-3)}\left\{\2n\ba{ll}
\ns\ds dX_{N-3}^*=\Big\{(A-B\Th_{N-3})X_{N-3}^*+\big[\bar
A+B(\Th_{N-3}-\h\Th_{N-3})-\bar
B\h\Th_{N-3}\big]\dbE_{t_{N-3}}[X_{N-3}^*]\Big\}ds\\
\ns\ds\qq\qq\q+\Big\{(C\1n-\1n D\Th_{N-3})X_{N-3}^*\2n+\2n\big[\bar
C\1n+\1n D(\Th_{N-3}\1n-\1n\h\Th_{N-3})\1n-\1n\bar
D\h\Th_{N-3}\big]\dbE_{t_{N-3}}[X_{N-3}^*]\Big\}dW(s),\\
\ns\ds\qq\qq\qq\qq\qq\qq\qq\qq\qq\qq\qq\qq\qq s\in[t_{N-3},t_{N-2}],\\
\ns\ds X_{N-3}^*(t_{N-3})=x_{N-3},\ea\right.\ee
and the optimal control $u^*_{N-3}(\cd)\equiv
u^*_{N-3}(\cd\,;t_{N-3},x_{N-3})$ admits the following
representation:
\bel{feedback(N-3)}
u_{N-3}^*(s)\1n=\1n-\Th_{N-3}(s)X_{N-3}^*(s)\1n+\1n\big[\Th_{N-3}(s)
\1n-\1n\h\Th_{N-3}(s)\big]\dbE_{t_{N-3}}[X_{N-3}^*(s)],\q
s\in[t_{N-3},t_{N-2}].\ee
Finally, the following holds:
\bel{inf(N-3)}\ba{ll}
\ns\ds\inf_{u_{N-3}(\cd)\in\sU[t_{N-3},\,t_{N-2}]}J^\D_{N-3}(x_{N-3};
u_{N-3}(\cd))=J^\D_{N-3}(x_{N-3};u^*_{N-3}(\cd))\\
\ns\ds\qq\qq\qq\qq\qq\qq\qq\qq=\lan\h
P_{N-3}(t_{N-3})x_{N-3},x_{N-3}\ran, \qq\forall
x_{N-3}\in\sX_{t_{N-3}}.\ea\ee
As a result, Player $(N-3)$ solves Problem (MF-LQ)$_{N-3}$. The
major difference between the situations of Player $(N-3)$ and Player
$(N-2)$ is that the Lyapunov equation for $\wt\G_{N-3}(\cd)$ has a
possible {\it impulse} at $t_{N-1}$. We actually have used Lemma 2.4 on
$[t_{N-1},t_N]$ and $[t_{N-2},t_{N-1}]$ separately.

\ms

By induction, we can construct the following finite sequences
$$\left\{\2n\ba{ll}
\ns\ds(P_k(\cd),\h P_k(\cd)),\qq\qq~0\le k\le N-1,\\
\ns\ds(\wt\G_k(\cd),\G_k(\cd),\bar\G_k(\cd)),\qq0\le k\le
N-2,\ea\right.$$
where $(P_k(\cd),\h P_k(\cd))$ satisfies the following Riccati
equation system: For $0\le k\le N-1$,
\bel{5Ric P(k)}\left\{\2n\ba{ll}
\ns\ds\dot P_k+P_kA+A^TP_k+C^TP_kC+Q_k\\
\ns\ds\q-(P_kB+C^TP_kD)(R_k+D^TP_kD)^{-1}(B^TP_k+D^TP_kC)=0,\\
\ns\ds\,\,\dot{\2n\h P}_k+\h P_k\h A+\h A^T\h P_k+\h C^TP_k\h C
+\h Q_k\\
\ns\ds\q-(\h P_k\h B+\h C^TP_k\h D)(\h R_k+\h D^TP_k\h D)^{-1}(\h
B^T\h P_k+\h D^TP_k
\h C)=0,\qq s\in[t_k,t_{k+1}),\\
\ns\ds P_k(t_{k+1})=\G_k(t_{k+1}),\q\h
P_k(t_{k+1})=\G_k(t_{k+1})+\bar\G_k(t_{k+1}),\ea\right.\ee
with the convention that
$$\G_{N-1}(t_N)=G_{N-1},\q\G_{N-1}(t_N)+\h\G_{N-1}(t_N)=\h G_{N-1}.$$
and $(\wt\G_k(\cd),\G_k(\cd),\bar\G_k(\cd))$ satisfies the following
Lyapunov equation systems: For $0\le k\le N-2$,

\bel{wt G(k)}\left\{\2n\ba{ll}
\ns\ds\dot{\wt\G}_k+\wt\G_k(A-B\Th^\D)+(A-B\Th^\D)^T\wt\G_k+Q_k\\
\ns\ds\qq+(C-D\Th^\D)^T\wt\G_k(C-D\Th^\D\1n)+(\Th^\D)^TR_k\Th^\D=0,
\qq s\in\bigcup_{i=k+1}^{N-1}(t_i,t_{i+1}),\\
\ns\ds\wt\G_k(t_N)=G_k,\q\wt\G_k(t_{N-1})=\G_k(t_{N-1}),\q\cds,\q
\wt\G_k(t_{k+2})=\G_k(t_{k+2}),\ea\right.\ee
\bel{G(k)}\left\{\2n\ba{ll}
\ns\ds\dot\G_k+\G_k(\h A-\h B\h\Th^\D)+(\h A-\h B\h\Th^\D)^T\G_k+Q_k\\
\ns\ds\qq+(\h C-\h D\h\Th^\D)^T\wt\G_k(\h C-\h D\h\Th^\D)+(\h\Th^\D)^T
R_k\h\Th^\D=0,\qq s\in[t_{k+1},t_N),\\
\ns\ds\G_k(t_N)=G_k,\ea\right.\ee
and
\bel{bG(k)}\left\{\2n\ba{ll}
\ns\ds\dot{\bar\G}_k\1n+\1n\bar\G_k(\h A\1n-\1n\h
B\h\Th^\D)\1n+\1n(\h A\1n-\1n\h B \h\Th^\D)^T\bar\G_k\1n+\1n\bar
Q_k\1n+\1n(\h\Th^\D)^T\1n\bar R_k\h\Th^\D
\1n=0,\q s\in[t_{k+1},t_N),\\
\ns\ds\bar\G_k(t_N)=\bar G_k.\ea\right.\ee
where, for $0\le k\le N-1$,
$$\left\{\2n\ba{ll}
\ns\ds\Th^\D(s)=\big[R_k(s)+D(s)^TP_k(s)D(s)\big]^{-1}\big[B(s)^TP_k(s)+D(s)^T
P_k(s)C(s)\big],\\
\ns\ds\h\Th^\D(s)=\big[\h R_k(s)+\h D(s)^TP_k(s)\h
D(s)\big]^{-1}\big[\h B(s)^T\h P_k(s)+\h D(s)^T P_k(s)\h
C(s)\big],\ea\right.\q s\in[t_k,t_{k+1}].$$
Next, we define
\bel{feedback*}
u^\D(s)=-\Th^\D(s)X^\D(s)+\big[\Th^\D(s)-\h\Th^\D(s)\big]
\dbE_{\rho^\D(s)}[X^\D(s)],\q s\in[0,t_N],\ee
where
\bel{rho}\rho^\D(s)=\sum_{k=0}^{N-2}t_kI_{[t_k,t_{k+1})}(s)+t_{N-1}I_{[t_{N-1},t_N]}(s),
\qq s\in[0,T],\ee
and define $X^\D(\cd)$ by the solution of the following closed-loop
system:
\bel{closed-loop(N-2)a}\left\{\2n\ba{ll}
\ns\ds dX^\D=\Big\{(A-B\Th^\D)X^\D+\big[\bar
A+B(\Th^\D-\h\Th^\D)-\bar
B\h\Th^\D\big]\dbE_{\rho^\D(s)}[X^\D]\Big\}ds\\
\ns\ds\qq\qq+\Big\{\1n(C\1n-\1n D\Th^\D\1n)X^\D\2n+\1n\big[\bar
C\1n+\1n D(\Th^\D\2n-\1n\h\Th^\D\1n)\1n-\1n\bar
D\h\Th^\D\big]\dbE_{\rho^\D(s)}[X^\D]\1n\Big\}dW\1n(s),\q s\1n\in\1n[0,t_N],\\
\ns\ds X^\D(0)=x.\ea\right.\ee
Then one has that
\bel{inf(k)}\ba{ll}
\ns\ds\inf_{u_k(\cd)\in\sU[t_k,\,t_{k+1}]}J^\D_k\big(X^\D(t_k);
u_k(\cd)\big)=J^\D_k\big(X^\D(t_k);u^\D(\cd)\big|_{[t_k,t_N]}\big)\\
\ns\ds\qq\qq\qq\qq\qq\qq\qq=\lan\h
P^\D(t_k)X^\D(t_k),X^\D(t_k)\ran,\qq0\le k\le N-1.\ea\ee
From the above, we see that $(\Th^\D(\cd),\h\Th^\D(\cd))$ is a closed-loop
$\D$-equilibrium strategy of Problem (MF-LQ), and
$(X^\D(\cd),u^\D(\cd))$ is the associated closed-loop $\D$-equilibrium pair.

\ms

To conclude this section, we note that
$\h\G_k(\cd)\equiv\G_k(\cd)+\bar\G_k(\cd)$ satisfies the following
Lyapunov equation:
\bel{hG(k)}\left\{\2n\ba{ll}
\ns\ds\dot{\h\G}_k+\h\G_k(\h A-\h B\h\Th^\D)+(\h A-\h B\h\Th^\D)^T\h\G_k+\h Q_k
+(\h C-\h D\h\Th^\D)^T\wt\G_k(\h C-\h D\h\Th^\D)+(\h\Th^\D)^T\1n\h R_k\h\Th^\D=0,\\
\ns\ds\qq\qq\qq\qq\qq\qq\qq\qq\qq\qq\qq\qq s\in[t_{k+1},t_N),\\
\ns\ds\h\G_k(t_N)=\h G_k.\ea\right.\ee

\section{Convergence.}

In the previous section, for any given partition $\D$ of $[0,T]$, a closed-loop
$\D$-equilibrium strategy
$(\Th^\D(\cd),\h\Th^\D(\cd))$ is constructed. The goal of this section is to
establish the convergence of $(\Th^\D(\cd),\h\Th^\D(\cd))$ as $\|\D\|\to0$. To achieve
this, we first present a technical result in which a partition
$\D:0=t_0<\cds<t_N=T$ of $[0,T]$ is given.

\ms

\bf Proposition 6.1. \sl Let {\rm(H1)--(H3)} hold. Let $Q_*,\h
Q_*,G_*,\h G_*\in\dbS^n$, and $R_*,\h R_*\in\dbS^m$ such that
\bel{}\left\{\2n\ba{ll}
\ns\ds Q(s,T)\le Q_*,\q\h Q(s,T)\le\h Q_*,\q R(s,T)\le R_*,\q\h
R(s,t)\le\h R_*,\qq s\in[0,T],\\
\ns\ds G(T)\le G_*,\q\h G(T)\le\h G_*.\ea\right.\ee
For $0\le k\le N-1$, let
\bel{6Ric P(ell)*}\left\{\2n\ba{ll}
\ns\ds\dot P^*_k\1n+\1n P^*_kA\1n+\1n A^T\1n P^*_k\1n+\1n C^T\1n P^*_kC\1n+\1n
Q_*\1n-\1n(P^*_kB\1n+\1n C^TP^*_kD)(R_*\1n+\1n D^T\1n P^*_kD)^{-1}(B^T\1n P^*_k
\1n+\1n D^T\1n P^*_kC)=0,\\
\ns\ds\dot{\h P^*_k}\1n+\1n\h P^*_k\h A\1n+\1n\h A^T\1n\h P^*_k\1n+\1n\h
C^T\1n P^*_k\h C\1n+\1n\h Q_*\1n-\1n(\h P^*_k\h B\1n+\1n\h C^T\1n P^*_k\h D)
(\h R_*\1n+\1n\h D^T\1n P^*_k\h
D)^{-1}(\h B^T\1n\h P^*_k\1n+\1n\h D^T\1n P^*_k\h C)=0,\\
\ns\ds\qq\qq\qq\qq\qq\qq\qq\qq\qq\qq\qq\qq\qq s\in[t_k,t_{k+1}),\\
\ns\ds P^*_k(t_{k+1})=\G^*_k(t_{k+1}),\qq\h P^*_k(t_{k+1})=\h
P^*_{k+1}(t_{k+1}),\ea\right.\ee
where for $1\le k\le N-1$,
\bel{6G(ell-1)*}\left\{\2n\ba{ll}
\ns\ds\dot\G_{k-1}^*+\G_{k-1}^*(\h A-\h B\h\Th^*_k)+(\h
A-\h B\h\Th^*_k)^T\G_{k-1}^*+Q_*\\
\ns\ds\qq+(\h C-\h D\h\Th^*_k)^TP_k^*(\h C-\h D\h\Th^*_k)
+(\h\Th^*_k)^TR_*\h\Th^*_k=0,\qq s\in[t_k,t_{k+1}),\\
\ns\ds\G_{k-1}^*(t_{k+1})=\G_k^*(t_{k+1}),\ea\right.\ee
with
$$\h\Th^*_k=(\h R_*+\h D^TP_k^*\h D)^{-1}(\h B^T\h P^*_k+\h
D^TP^*_k\h C),\q s\in[t_k,t_{k+1}],\q0\le k\le N-1.$$
and with the convention that
$$\G^*_{N-1}(t_N)=G_*,\qq\h P^*_N(t_N)=\h G_*.$$
Further, for $0\le k\le N-1$, let
\bel{6Pi(ell)*}\left\{\2n\ba{ll}
\ns\ds\dot\Pi^*_k+\Pi^*_kA+A^T\Pi^*_k+C^T\Pi^*_kC+Q_*=0,\qq s\in[t_k,t_{k+1}),\\
\ns\ds\Pi^*_k(t_{k+1})=\G^*_k(t_{k+1}),\ea\right.\ee
and
\bel{}\left\{\2n\ba{ll}
\ns\ds\1n\dot{\;\h\Pi_*}+\h\Pi_*\h A+\h A^T\h\Pi_*+\h
C^T\Pi_*\h C+\h Q_*=0,\qq s\in[0,T),\\
\ns\ds\h\Pi_*(T)=\h G_*.\ea\right.\ee
Then
\bel{6.6*}\left\{\2n\ba{ll}
\ns\ds0\le\wt\G_\ell(s)\le P_k(s)\le P^*_k(s)\le\Pi^*_k(s),\\
\ns\ds0\le\h\G_\ell(s)\le\h P_k(s)\le\h
P^*_k(s)\le\h\Pi_*(s),\ea\right.\qq s\in[t_k,t_{k+1}],\q0\le\ell<k\le
N-1,\ee
with the convention that $\wt\G_{-1}(s)=\h\G_{-1}(s)=0$.

\ms

\rm

\it Proof. \rm We look at the case $k=N-1$. Recall that
\bel{6Ric P(N-1)}\left\{\2n\ba{ll}
\ns\ds\dot P_{N-1}+P_{N-1}A+A^TP_{N-1}+C^TP_{N-1}C+Q_{N-1}\\
\ns\ds\q-(P_{N-1}B+C^TP_{N-1}D)(R_{N-1}+D^TP_{N-1}D)^{-1}(B^TP_{N-1}+D^TP_{N-1}C)=0,\\
\ns\ds\;\dot{\2n\h P}_{N-1}+\h P_{N-1}\h A+\h A^T\h P_{N-1}+\h
C^TP_{N-1}\h C
+\h Q_{N-1}\\
\ns\ds\q-(\h P_{N-1}\h B\1n+\1n\h C^T\1n P_{N-1}\h D)(\h R_{N-1}+\h
D^TP_{N-1}\h D)^{-1}(\h B^T\h P_{N-1}\1n+\1n\h D^T\1n P_{N-1}
\h C)=0,\\
\ns\ds\qq\qq\qq\qq\qq\qq\qq\qq s\in[t_{N-1},t_N),\\
\ns\ds P_{N-1}(t_N)=G_{N-1},\qq\h P_{N-1}(t_N)=\h
G_{N-1}.\ea\right.\ee
According to Proposition 3.2, together with the usual comparison of solutions to
Riccati equations, we have
\bel{6.5}\left\{\2n\ba{ll}
\ns\ds0\le P_{N-1}(s)\le P_{N-1}^*(s)\le\Pi_{N-1}^*(s),\\
\ns\ds0\le\h P_{N-1}(s)\le\h
P_{N-1}^*(s)\le\h\Pi_{N-1}^*(s),\ea\right.\qq s\in[t_{N-1},t_N].\ee
Define
$$\left\{\2n\ba{ll}
\ns\ds\Th_{N-1}=(R_{N-1}+D^TP_{N-1}D)^{-1}(B^TP_{N-1}+D^TP_{N-1}C),\\
\ns\ds\h\Th_{N-1}=(\h R_{N-1}+\h D^TP_{N-1}\h D)^{-1}(\h B^T\h
P_{N-1}+\h D^TP_{N-1}\h C),\ea\right.\q s\in[t_{N-1},t_N].$$
Then Riccati equation system (\ref{6Ric P(N-1)}) can be written as
\bel{6Ric P(N-1)e}\left\{\2n\ba{ll}
\ns\ds\dot P_{N-1}+P_{N-1}(A-B\Th_{N-1})+(A-B\Th_{N-1})^TP_{N-1}+Q_{N-1}\\
\ns\ds\q+(C-D\Th_{N-1})^TP_{N-1}(C-D\Th_{N-1})+\Th_{N-1}^TR_{N-1}\Th_{N-1}=0,\\
\ns\ds\;\dot{\2n\h P}_{N-1}+\h P_{N-1}(\h A-\h B\h\Th_{N-1})+(\h
A-\h
B\h\Th_{N-1})^T\h P_{N-1}+\h Q_{N-1}\\
\ns\ds\q+(\h C-\h D\h\Th_{N-1})^TP_{N-1}(\h C-\h
B\h\Th_{N-1})+\h\Th_{N-1}^T\h R_{N-1}\h\Th_{N-1}=0,\qq s\in[t_{N-1},t_N),\\
\ns\ds P_{N-1}(t_N)=G_{N-1},\q\h P_{N-1}(t_N)=\h
G_{N-1}.\ea\right.\ee
Next, for $\ell=N-2,N-3,\cds,2,1,0$, we recall
\bel{6wtG(N-1)}\left\{\2n\ba{ll}
\ns\ds\2n\dot{~\wt\G_\ell}+\wt\G_\ell(A-B\Th_{N-1})+(A-B\Th_{N-1})^T\wt\G_\ell
+Q_\ell\\
\ns\ds\qq+(C-D\Th_{N-1})^T\wt\G_\ell(C-D\Th_{N-1}\1n)+(\Th_{N-1})^TR_\ell
\Th_{N-1}=0,
\qq s\in[t_{N-1},t_N),\\
\ns\ds\wt\G_\ell(t_N)=G_\ell,\ea\right.\ee
\bel{6G(N-1)}\left\{\2n\ba{ll}
\ns\ds\dot\G_\ell+\G_\ell(\h A-\h B\h\Th_{N-1})+(\h A-\h
B\h\Th_{N-1})^T\G_\ell+Q_\ell\\
\ns\ds\qq+(\h C-\h D\h\Th_{N-1})^T\wt\G_\ell(\h C-\h D\h\Th_{N-1})
+(\h\Th_{N-1})^T\1n R_\ell\h\Th_{N-1}=0,\qq s\in[t_{N-1},t_N),\\
\ns\ds\G_\ell(t_N)=G_\ell,\ea\right.\ee
and
\bel{6hG(N-1)}\left\{\2n\ba{ll}
\ns\ds\dot{\h\G}_\ell+\h\G_\ell(\h A-\h B\h\Th_{N-1})+(\h A-\h
B\h\Th_{N-1})^T\h\G_\ell+\h Q_\ell\\
\ns\ds\qq+(\h C-\h D\h\Th_{N-1})^T\wt\G_\ell(\h C-\h
D\h\Th_{N-1})+(\h\Th_{N-1})^T\h R_\ell\h\Th_{N-1}=0,\qq s\in[t_{N-1},t_N),\\
\ns\ds\h\G_\ell(t_N)=\h G_\ell.\ea\right.\ee
Comparing (\ref{6Ric P(N-1)}) with (\ref{6wtG(N-1)}) and
(\ref{6hG(N-1)}), making use of (H3), we see that
\bel{}\left\{\2n\ba{ll}
\ns\ds0\le\wt\G_\ell(s)\le P_{N-1}(s),\\
\ns\ds0\le\h\G_\ell(s)\le\h P_{N-1}(s),\ea\right.\qq s\in[t_{N-1},t_N],\q0\le\ell\le N-2.
\ee
This proves (\ref{6.6*}) for $k=N-1$. The proofs for general
$k=N-2,N-3,\cds,2,1,0$ are essentially the same. \endpf

\ms

Since $\h\Pi_*(\cd)$ is bounded uniform in $\D$, we see from
(\ref{6.6*}) that $\h\G_k(\cd)$ and $\h P_k(\cd)$ are
uniformly bounded. We now establish the uniform boundedness of
$\wt\G_k(\cd)$, $\G_k(\cd)$ and $P_k(\cd)$.

\ms

\bf Proposition 6.2. \sl Let {\rm (H1)--(H3)} hold. Then
$\wt\G_k(\cd)$, $\G_k(\cd)$ and $P_k(\cd)$ are bounded uniformly in
partition $\D$ of $[0,T]$.

\ms

\it Proof. \rm Let $P_k^*(\cd)$, $\Pi_k^*(\cd)$ and $\G_k^*(\cd)$ be as in
Proposition 6.1. Choose $\b>\|\h R(\cd\,,\cd)\|_\infty$ large enough so that
\bel{b}\b>\|\h C(\cd)\|_\infty\|\h D(\cd)\|_\infty|P_k^*(s)|,\qq s\in[t_k,t_{k+1}],~
0\le k\le N-1.\ee
Then
\bel{BPB}\h D(s)^TP_k^*(s)\h C(s)\h C(s)^TP_k^*(s)\h D(s)\le\big[\b I+\h D(s)^TP_k^*(s)\h D(s)\big]^2,\q
s\in[t_k,t_{k+1}],~0\le k\le N-1.\ee
We claim that
\bel{|Th|}|\h\Th_k^*(s)|\le\b^{-1}\|\h
B(\cd)\|_\infty\|\h\Pi_*(\cd)\|_\infty+1,\qq
s\in[t_k,t_{k+1}],~0\le k\le N-1.\ee
In fact, by taking $\h R_*=\b I$, we have
$$\ba{ll}
\ns\ds0\le(\h R_*+\h D^TP_k^*\h D)^{-1}\h D^TP_k^*\h C\h
C^TP_k^*\h D(\h R_*+\h D^TP_k^*\h D)^{-1}\\
\ns\ds\q=(\h R_*+\h D^TP_k^*\h D)^{-1}\[\h D^TP_k^*\h C\h C^TP_k^*\h
D-(\h R_*+\h D^TP_k^*\h D)^2\](\h R_*+\h D^TP_k^*\h D)^{-1}+I\le
I.\ea$$
Then
$$\ba{ll}
\ns\ds|\h\Th_k^*|\le|(\h R_*+\h D^TP_k^*\h D)^{-1}||\h B||\h
P_k^*|+|(\h R_*+\h D^TP_k^*\h D)^{-1}\h D^TP_k^*\h C|\le\b^{-1}|\h B||\h\Pi_*|+1,\ea$$
from which one obtains (\ref{|Th|}).

\ms

Next, we recall the equations (\ref{6G(ell-1)*}) and (\ref{6Pi(ell)*}) for $\G^*_k(\cd)$ and $\Pi^*_k(\cd)$. Define
$$\G^*_\D(s)=\sum_{k=1}^{N-1}\G^*_{k-1}(s)I_{(t_k,t_{k+1}]}(s),\q
\Pi^*_\D(s)=\sum_{k=1}^{N-1}\Pi^*_k(s)I_{(t_k,t_{k+1}]}(s),\qq
s\in[0,T].$$
Let us begin on $(t_{N-1},t_N]$. We have
$$\ba{ll}
\ns\ds|\Pi^*_\D(s)|=|\Pi_{N-1}^*(s)|\le|G_*|+\int_s^{t_N}\[(2|A(r)|+|C(r)|^2)|\Pi^*_{N-1}(r)|+|Q_*|\]dr,\\
\ns\ds\qq\qq\;\le|G_*|+\int_s^{t_N}\[\(2\|A(\cd)\|_\infty+\|C(\cd)\|^2_\infty\)|\Pi_{N-1}^*(r)|
+|Q_*|\]dr\\
\ns\ds\qq\qq\;\equiv|G_*|+|Q_*|(t_N-s)+K_\Pi\int_s^{t_N}|\Pi_\D^*(r)|dr,\qq
s\in[t_{N-1},t_N],\ea$$
with
$$K_\Pi=2\|A(\cd)\|_\infty+\|C(\cd)\|^2_\infty,$$
and
$$\ba{ll}
\ns\ds|\G^*_\D(s)|=|\G_{N-2}^*(s)|\le|G_*|+\int_s^{t_N}\[2|\h
A(r)-\h
B(r)\h\Th_{N-1}^*(s)|\,|\G_{N-2}^*(r)|+|Q_*|\\
\ns\ds\qq\qq\qq\qq\q+|\h C(r)-\h
D(r)\h\Th^*_{N-1}(r)|^2|P^*_{N-1}(r)|+|R_*||\h\Th^*_{N-1}(r)|^2\]dr\\
\ns\ds\qq\le|G_*|+\int_s^{t_N}\[2\(\|\h A(\cd)\|_\infty+\|\h
B(\cd)\|_\infty(\b^{-1}\|\h
B(\cd)\|_\infty\|\h\Pi_*(\cd)\|_\infty+1)\)|\G_{N-2}^*(r)|+|Q_*|\\
\ns\ds\qq\qq+\(\|\h C(\cd)\|_\infty+\|\h D(\cd)\|_\infty
\big(\b^{-1}\|\h B(\cd)\|_\infty\|\h\Pi_*(\cd)\|_\infty+1\big)\)^2|\Pi_{N-1}^*(r)|\\
\ns\ds\qq\qq+|R_*|\big(\b^{-1}\|\h B(\cd)\|_\infty\|\h\Pi_*(\cd)\|_\infty+1
\big)^2\]dr\\
\ns\ds\qq\equiv|G_*|+K_Q(t_N-s)+\int_s^{t_N}\[K_\G|\G_\D^*(r)|
+\wt K_\Pi|\Pi^*_\D(r)|\]dr,\qq s\in[t_{N-1},t_N],\ea$$
with
$$\left\{\2n\ba{ll}
\ns\ds K_\G=2\(\|\h A(\cd)\|_\infty+\|\h B(\cd)\|_\infty(\b^{-1}\|\h B(\cd)\|_\infty\|\h\Pi_*(\cd)\|_\infty+1)\),\\
\ns\ds\wt K_\Pi=\(\|\h C(\cd)\|_\infty+\|\h D(\cd)\|_\infty(\b^{-1}\|\h B(\cd)\|_\infty\|\h\Pi_*(\cd)\|_\infty+1)\)^2,\\
\ns\ds K_Q=|Q_*|+|R_*|\big(\b^{-1}\|\h B(\cd)\|_\infty\|\h\Pi_*(\cd)\|_\infty+1\big)^2.\ea\right.$$
Thus,
\bel{}\ba{ll}
\ns\ds|\Pi^*_\D(s)|+|\G^*_\D(s)|\\
\ns\ds\le2|G_*|+(|Q_*|+K_Q)(t_N-s)+\int_s^{t_N}\[(K_\Pi\vee\wt
K_\Pi)|\Pi^*_\D(r)|+K_\G|\G^*_\D(r)|\]dr\\
\ns\ds\le2|G_*|+2K_Q(t_N-s)+\int_s^{t_N}\[(K_\Pi+\wt
K_\Pi)|\Pi^*_\D(r)|+K_\G|\G^*_\D(r)|\]dr,\qq
s\in[t_{N-1},t_N].\ea\ee
Next, on $(t_{N-2},t_{N-1}]$,
$$\ba{ll}
\ns\ds|\Pi^*_\D(s)|=|\Pi_{N-2}^*(s)|\le|\G^*_{N-2}(t_{N-1})|+\int_s^{t_{N-1}}\[(2|A(r)|+|C(r)|^2)|\Pi^*_{N-2}(r)|+|Q_*|\]dr\\
\ns\ds\qq=|G_*|\1n+\1n K_Q(t_N\1n-\1n
t_{N-1})\1n+\3n\int_{t_{N-1}}^{t_N}\2n\(K_\G|\G_{N-2}^*(r)|\1n+\1n
\wt
K_\Pi|\Pi^*_{N-1}(r)|\)dr\1n+\3n\int_s^{t_{N-1}}\2n\(K_\Pi|\Pi^*_{N-2}(r)|\1n
+\1n|Q_*|\)dr\\
\ns\ds\qq\le|G_*|+K_Q(t_N-s)+\int_{t_{N-1}}^{t_N}\(K_\G|\G_{N-2}^*(r)|+
\wt
K_\Pi|\Pi^*_{N-1}(r)|\)dr+K_\Pi\int_s^{t_{N-1}}|\Pi^*_{N-2}(r)|dr,\\
\ns\ds\qq\qq\qq\qq\qq\qq\qq\qq\qq\qq\qq\qq\qq
s\in[t_{N-2},t_{N-1}],\ea$$
and
$$\ba{ll}
\ns\ds|\G^*_\D(s)|=|\G_{N-3}^*(s)|\le|\G^*_{N-2}(t_{N-1})|+\int_s^{t_{N-1}}\[2|\h
A(r)-\h
B(r)\h\Th_{N-2}^*(s)|\,|\G_{N-3}^*(r)|+|Q_*|\\
\ns\ds\qq\qq\qq\qq\q+|\h C(r)-\h
D(r)\h\Th^*_{N-2}(r)|^2|P^*_{N-2}(r)|+|R_*||\Th^*_{N-2}(r)|^2\]dr\\
\ns\ds\le|G_*|+K_Q(t_N-t_{N-1})+\int_{t_{N-1}}^{t_N}\(K_\G|\G^*_{N-2}(r)|+
\wt K_\Pi|\Pi^*_{N-1}(r)|\)dr\\
\ns\ds\qq\qq\qq+\int_s^{t_{N-1}}\(K_\G|\G^*_{N-3}(r)|+\wt
K_\Pi|\Pi^*_{N-2}(r)|+|Q_*|\)dr\\
\ns\ds=|G_*|+K_Q(t_N-s)+\int_s^{t_N}\(K_\G|\G^*_\D(r)|+ \wt
K_\Pi|\Pi^*_\D(r)|\)dr,\qq s\in[t_{N-2},t_{N-1}].\ea$$
Then
\bel{}\ba{ll}
\ns\ds|\Pi^*_\D(s)|+|\G^*_\D(s)|\le2|G_*|+2K_Q(t_N-s)
+2\int_{t_{N-1}}^{t_N}\(K_\G|\G^*_\D(r)|+\wt
K_\Pi|\Pi^*_\D(r)|\)dr\\
\ns\ds\qq\qq\qq\qq\qq\q+\int_s^{t_{N-1}}\[(K_\Pi+\wt
K_\Pi)|\Pi^*_\D(r)|+K_\G|\G^*_\D(r)|\]dr,\qq
s\in[t_{N-2},t_{N-1}].\ea\ee
By induction, we can show that
$$\ba{ll}
\ns\ds|\Pi^*_\D(s)|+|\G^*_\D(s)|\le2|G_*|+2K_Q(T-s)+\int_s^T\[(K_\Pi\vee\wt
K_\Pi+\wt K_\Pi)|\Pi^*_\D(r)|+2K_\G|\G^*_\D(r)|\]dr\\
\ns\ds\qq\qq\qq\qq\le2|G_*|+2K_Q(T-s)+K_0\int_s^T\[|\Pi^*_\D(r)|+|\G^*_\D(r)|\]dr,\qq
s\in[0,T],\ea$$
with
$$K_0=\max\{K_\Pi\vee\wt K_\Pi+\wt K_\Pi,2K_\G\}.$$
Hence, by Gronwall's inequality, one obtains
$$|\Pi^*_\D(s)|+|\G_\D^*(s)|\le2(|G_*|+K_QT)e^{K_0T},\qq
s\in[0,T],$$
which implies
$$|P_k^*(s)|\le2\big(|G_*|+K_QT\big)e^{K_0T},\qq s\in[t_k,t_{k+1}],\q0\le k\le N-1.$$
Note that
$$\ba{ll}
\ns\ds K_0\le K_\Pi+2(\wt K_\Pi+K_\G)\\
\ns\ds=2\|A(\cd)\|_\infty+\|C(\cd)\|_\infty^2+2\(\|\h
C(\cd)\|_\infty+\|\h D(\cd)\|_\infty(\b^{-1}\|\h
B(\cd)\|_\infty\|\h\Pi_*(\cd)\|_\infty+1)\)^2\\
\ns\ds\q+2\(\|\h A(\cd)\|_\infty+\|\h
B(\cd)\|_\infty(\b^{-1}\|\h
B(\cd)\|_\infty\|\h\Pi_*(\cd)\|_\infty+1)\)\\
\ns\ds\le2\|A(\cd)\|_\infty+\|C(\cd)\|_\infty^2+2\(\|\h
C(\cd)\|_\infty+\|\h D(\cd)\|_\infty(\|\h
B(\cd)\|_\infty\|\h\Pi_*(\cd)\|_\infty+1)\)^2\\
\ns\ds\q+2\(\|\h A(\cd)\|_\infty+\|\h
B(\cd)\|_\infty(\|\h
B(\cd)\|_\infty\|\h\Pi_*(\cd)\|_\infty+1)\)\equiv\bar K_0,\ea$$
and
$$K_Q\le|Q_*|+|R_*|\big(\|\h B(\cd)\|_\infty\|\h\Pi_*(\cd)\|_\infty+1\big)^2
\equiv\bar K_Q,$$
with $\bar K_0$ and $\bar K_Q$ being absolute constants. Now, we may take
$$\b=1+2\big(|G_*|+\bar K_QT\big)e^{\bar K_0T}\|\h C(\cd)\|_\infty
\|\h D(\cd)\|_\infty.$$
Then (\ref{b}) holds and we obtain uniform boundedness of $P_k^*(\cd)$,
$\G^*_k(\cd)$, and $\Pi^*_k(\cd)$. Consequently, from (\ref{6.6*}), we obtain
the uniform boundedness of $P_k(\cd)$ and $\wt\G_k(\cd)$. Finally, by the
equation (\ref{G(k)}) for $\G_k(\cd)$, we obtain the uniform boundedness of
$\G_k(\cd)$. \endpf

\ms

The following result essentially gives some estimates on the jumps.

\ms

\bf Proposition 6.3. \sl Let {\rm(H1)--(H3)} hold. Then
\bel{6.5}\ba{ll}
\ns\ds|\G_k(s)-\G_{k-1}(s)|+|\wt\G_k(s)-\wt\G_{k-1}(s)|\\
\ns\ds\le
K\[|G_k-G_{k-1}|+\int_s^{t_N}\(|Q_k(r)-Q_{k-1}(r)|+|R_k(r)-R_{k-1}(r)|\)dr\],\q
s\in[t_{k+1},t_N],\ea\ee
and
\bel{6.6}\ba{ll}
\ns\ds|\h\G_k(s)-\h\G_{k-1}(s)|\le K\[|G_k-G_{k-1}|+|\h G_k-\h
G_{k-1}|+\int_s^{t_N}\(|Q_k(r)-Q_{k-1}(r)|\\
\ns\ds\qq\qq\qq\qq\qq+|R_k(r)-R_{k-1}(r)|\1n+\1n|\h Q_k(r)-\h
Q_{k-1}(r)|\1n+\1n|\h
R_k(r)-\h R_{k-1}(r)|\)dr\],\\
\ns\ds\qq\qq\qq\qq\qq\qq\qq\qq\qq\qq\qq\qq\qq\qq
s\in[t_{k+1},t_N].\ea\ee

\ms

\it Proof. \rm Under (H1)--(H3), we know that $P_k(\cd),\h P_k(\cd)$
are uniformly bounded. Thus, all the coefficients of the Lyapunov
equations (\ref{wt G(k)})--(\ref{bG(k)}) and (\ref{hG(k)}) are
bounded. Now, comparing the equations for $\G_k(\cd)$ and
$\G_{k-1}(\cd)$, we have
\bel{6.7}\ba{ll}
\ns\ds|\G_k(s)-\G_{k-1}(s)|\le|G_k-G_{k-1}|+K\int_s^{t_N}\(|\G_k(r)-
\G_{k-1}(r)|+|\wt\G_k(r)-\wt\G_{k-1}(r)|\)dr\\
\ns\ds\qq\qq\qq\qq\q+K\2n\int_s^{t_N}\2n\(|Q_k(r)-Q_{k-1}(r)|+|R_k(r)-R_{k-1}(r)|
\)dr,\q s\in[t_{k+1},t_N].\ea\ee
Next, comparing the equations for $\wt\G_k(\cd)$ and
$\wt\G_{k-1}(\cd)$, we see that for $s\in(t_{N-1},t_N]$,
$$\ba{ll}
\ns\ds|\wt\G_k(s)-\wt\G_{k-1}(s)|\le|G_k-G_{k-1}|+K\int_s^{t_N}
|\wt\G_k(r)-\wt\G_{k-1}(r)|dr\\
\ns\ds\qq\qq\qq\qq\qq+K\int_s^{t_N}\(|Q_k(r)-Q_{k-1}(r)|+|R_k(r)-R_{k-1}(r)|
\)dr.\ea$$
For $s\in(t_{N-2},t_{N-1}]$,
$$\ba{ll}
\ns\ds|\wt\G_k(s)-\wt\G_{k-1}(s)|\le
|\G_k(t_{N-1})-\G_{k-1}(t_{N-1})|+K\int_s^{t_{N-1}}|\wt\G_k(s)-
\wt\G_{k-1}(r)|dr\\
\ns\ds\qq\qq\qq\qq\qq+K\int_s^{t_{N-1}}\(|Q_k(r)-Q_{k-1}(r)|+|R_k(r)-R_{k-1}(r)|
\)dr\\
\ns\ds\le|G_k-G_{k-1}|+K\int_{t_{N-1}}^{t_N}\(|\G_k(r)-\G_{k-1}(r)|
+|\wt\G_k(r)-\wt\G_{k-1}(r)|\)dr\\
\ns\ds\qq\qq+K\int_{t_{N-1}}^{t_N}\(|Q_k(r)-Q_{k-1}(r)|+|R_k(r)-R_{k-1}(r)|\)dr\\
\ns\ds\qq\qq+K\int_s^{t_{N-1}}\(|\wt\G_k(s)-
\wt\G_{k-1}(r)|+|Q_k(r)-Q_{k-1}(r)|+|R_k(r)-R_{k-1}(r)|
\)dr\\
\ns\ds\le|G_k-G_{k-1}|+K\int_s^{t_N}\(|\G_k(r)-\G_{k-1}(r)|
+|\wt\G_k(r)-\wt\G_{k-1}(r)|\)dr\\
\ns\ds\qq\qq+K\int_s^{t_N}\(|Q_k(r)-Q_{k-1}(r)|+|R_k(r)-R_{k-1}(r)|\)dr.\ea$$
For $s\in(t_{N-3},t_{N-2}]$,
$$\ba{ll}
\ns\ds|\wt\G_k(s)-\wt\G_{k-1}(s)|\le
|\G_k(t_{N-2})-\G_{k-1}(t_{N-2})|+K\int_s^{t_{N-2}}|\wt\G^\D_k(s)-
\wt\G^\D_{k-1}(r)|dr\\
\ns\ds\qq\qq\qq\qq\qq+K\int_s^{t_{N-2}}\big[|Q_k(r)-Q_{k-1}(r)|+|R_k(r)-R_{k-1}(r)|
\big]dr\\
\ns\ds\le|G_k-G_{k-1}|+K\int_{t_{N-2}}^{t_N}\(|\G_k(r)-\G_{k-1}(r)|
+|\wt\G_k(r)-\wt\G_{k-1}(r)|\)dr\\
\ns\ds\qq\qq+K\int_{t_{N-2}}^{t_N}\(|Q_k(r)-Q_{k-1}(r)|+|R_k(r)-R_{k-1}(r)|
\)dr\\
\ns\ds\qq\qq+K\int_s^{t_{N-2}}\(|\wt\G_k(s)-
\wt\G_{k-1}(r)|+|Q_k(r)-Q_{k-1}(r)|+|R_k(r)-R_{k-1}(r)|\)dr\\
\ns\ds\le|G_k-G_{k-1}|+K\int_s^{t_N}\(|\G_k(r)-\G_{k-1}(r)|+
|\wt\G_k(r)-\wt\G_{k-1}(r)|\)dr\\
\ns\ds\qq\qq+K\int_s^{t_N}\(|Q_k(r)-Q_{k-1}(r)|+|R_k(r)-R_{k-1}(r)|
\)dr.\ea$$
By induction, we see that
$$\ba{ll}
\ns\ds|\wt\G_k(s)-\wt\G_{k-1}(s)|\le|G_k-G_{k-1}|+K\int_s^{t_N}\(
|\G_k(r)-\G_{k-1}(r)|+|\wt\G_k(r)-\wt\G_{k-1}(r)|\)dr\\
\ns\ds\qq\qq\qq\qq+K\int_s^{t_N}\(|Q_k(r)-Q_{k-1}(r)|+|R_k(r)-R_{k-1}(r)|\)dr,\q
s\in[t_{k+1},t_N].\ea$$
Combining (\ref{6.7}), we obtain
$$\ba{ll}
\ns\ds|\G_k(s)-\G_{k-1}(s)|+|\wt\G_k(s)-\wt\G_{k-1}(s)|\\
\ns\ds\le2|G_k-G_{k-1}|+K\int_s^{t_N}\( |\G_k(r)-\G_{k-1}(r)|
+|\wt\G_k(r)-\wt\G_{k-1}(r)|\)dr\\
\ns\ds\qq\qq\qq\qq+K\int_s^{t_N}\(|Q_k(r)-Q_{k-1}(r)|+|R_k(r)-R_{k-1}(r)|\)dr,\q
s\in[t_{k+1},t_N].\ea$$
By Gronwall's inequality, we obtain (\ref{6.5}). Similarly,
comparing equations for $\h\G_k(\cd)$ and $\h\G_{k-1}(\cd)$, we
obtain
$$\ba{ll}
\ns\ds|\h\G_k(s)-\h\G_{k-1}(s)|\le|\h G_k-\h G_{k-1}|+K\int_s^{t_N}
|\h\G_k(r)-\h\G_{k-1}(r)|dr\\
\ns\ds\qq\qq\qq\qq+\int_s^{t_N}\(|\h Q_k(r)-\h Q_{k-1}(r)|+|\h
R_k(r)-\h R_{k-1}(r)|+|\wt\G_k(r)-\wt\G_{k-1}(r)|\)dr.\ea$$
By Gronwall's inequality, we obtain (\ref{6.6}). \endpf

\ms

Now, for given partition $\D$ of $[0,T]$, let us denote
\bel{}\left\{\2n\ba{ll}
\ns\ds Q^\D(s)=\sum_{k=0}^{N-1} Q(s,t_k)I_{(t_k,t_{k+1}]}(s),\q \h
Q^\D(s)=\sum_{k=0}^{N-1}\2n\bar
Q(s,t_k)I_{(t_k,t_{k+1}]}(s),\\
\ns\ds R^\D(s)=\sum_{k=0}^{N-1} R(s,t_k)I_{(t_k,t_{k+1}]}(s),\q\h
R^\D(s)=\sum_{k=0}^{N-1}\h R(s,t_k)I_{(t_k,t_{k+1}]}(s),\ea\right.\q
s\in[0,t_N].\ee
The dependence on the partition $\D$ is indicated. Denote
\bel{}P^\D(s)=\sum_{k=0}^{N-1}P_k(s)I_{(t_k,t_{k+1}]}(s),\q \h
P^\D(s)=\sum_{k=0}^{N-1}\h P_k(s)I_{(t_k,t_{k+1}]}(s),\qq
s\in[0,t_N],\ee
and
\bel{G(s,t)}\left\{\ba{ll}
\ns\ds\wt\G^\D(s,\t)=\sum_{k=0}^{N-1}\wt\G_k^\D(s)I_{(t_{k+1},t_{k+2}]}(\t),\\
\ns\ds\G^\D(s,\t)=\sum_{k=0}^{N-1}\G_k^\D(s)I_{(t_{k+1},t_{k+2}]}(\t),\\
\ns\ds\bar\G^\D(s,\t)=\sum_{k=0}^{N-1}\bar\G_k^\D(s)I_{(t_{k+1},t_{k+2}]}(\t),
\ea\right.\qq0\le\t\le s\le t_N,\ee
with the following extension on $[0,t_N]\times[0,t_N]$ as follows:
\bel{}\wt\G^\D(s,\t)=\wt\G^\D(\t,\t),\q\G^\D(s,\t)=\G^\D(\t,\t),\q\bar\G^\D(s,\t)
=\bar\G^\D(\t,\t),\qq s\in[0,\t].\ee
Note that $s\mapsto(\G^\D(s,\t),\bar\G^\D(s,\t))$ is continuous on
$[0,T]$, whereas $s\mapsto\wt\G^\D(s,\t)$ might have jumps at
$s=t_k$, with $t_k\in[\t,T]$.

\ms

From Propositions 6.1 and 6.2, we have the uniform boundedness of
$P^\D(\cd),\h P^\D(\cd),\wt\G^\D(\cd\,,\cd),\G^\D(\cd\,,\cd)$, and
$\h\G^\D(\cd\,,\cd)$. This implies the equicontinuity of $s\mapsto
(P^\D(s),\h P^\D(s),\wt\G^\D(s,\t),\G^\D(s,\t),\h\G^\D(s,\t))$ on
each $(t_k,t_{k+1})$, for fixed $\t\in[0,T]$. Also, from Proposition
6.3, we see the jump sizes of the above functions are controlled by
the mesh size $\|\D\|$ of the partition $\D$, thanks to the
continuity of functions $(s,t)\mapsto(Q(s,t),\h Q(s,t),R(s,t),\h
R(s,t),G(t),\h G(t))$.

\ms

Finally, we will need the following result which is a kind of
extension of the well-known Arzela-Ascoli Theorem.

\ms

\bf Lemma 6.4. \sl Let $f_k:[0,T]\to\dbR^n$, $k\ge1$, be a sequence
of uniformly bounded piecewise continuous functions. Suppose that
for any $\e>0$, there exists a $\d(\e)>0$ and a $k_0(\e)\ge1$ such
that
\bel{c}|f_k(t)-f_k(s)|<\e,\qq\forall|t-s|<\d(\e),~k>k_0(\e).\ee
Then there exists a subsequence $f_{k_j}(\cd)$ and a continuous
function $f:[0,T]\to\dbR^n$ such that
$$\lim_{j\to\infty}\|f_{k_j}(\cd)-f(\cd)\|_\infty=0.$$

\ms

\it Sketch of the Proof. \rm Let $\dbQ[0,T]\equiv\{r_1,r_2,\cds\}$
be the set of all rational numbers in $[0,T]$. Without loss of
generality, let us assume that all the functions $f_k(\cd)$ are
right-continuous. First consider the sequence $\{f_k(r_1)\}_{k\ge1}$
which is bounded. Thus, we may have a convergent subsequence,
denoted by, $\{f_{\si_1(j)}(r_1)\}_{j\ge1}\subseteq\{f_k\}_{k\ge1}$.
Next, consider sequence $\{f_{\si_1(j)}(r_2)\}_{j\ge1}$. By the
boundedness, we again have a convergent subsequence
$\{f_{\si_2(j)}(r_2)\}_{j\ge1}\subseteq\{f_{\si_1(j)}(r_2)\}_{j\ge1}$.
Continuing the procedure, and using a usual diagonal argument, we
obtain a subsequence $\{f_{k_j}(\cd)\}_{j\ge1}$ and a function
$f(r)$ defined for all $r\in\dbQ[0,T]$ such that $f_{k_j}(r)\to
f(r)$ for all $r\in\dbQ[0,T]$. By (\ref{c}), we see that $f(\cd)$ is
bounded and continuous on $\dbQ[0,T]$ which is dense in $[0,T]$.
Then one can extend it naturally to a continuous function on
$[0,T]$. Further, we see that our convergence holds. \endpf

\ms

With the same argument, we see that the above result holds if
$[0,T]$ is replaced by $[0,T]\times[0,T]$ (actually, $[0,T]$ can be
replaced by any compact metric space). The point here is that the
functions $f_k(\cd)$ are not necessarily continuous, which will be
the case in our problem.

\ms

Having the above preparation, we are now at the position of stating
and proving the following convergence theorem.

\ms

\bf Theorem 6.5. \sl Let {\rm(H1)--(H3)} hold. Then
\bel{}\ba{ll}
\ns\ds\lim_{\|\D\|\to0}\(|\wt\G^\D(s,\t)-\G(s,\t)|+|\G^\D(s,\t)-\G(s,\t)|
+|\h\G^\D(s,\t)-\h\G(s,t)|\\
\ns\ds\qq\qq+|P^\D(s)-\G(s,s)|+|\h P^\D(s)-\h\G(s,s)|\)=0,\ea\ee
uniformly in $(s,\t)\in[0,T]\times[0,T]$ with $\G(\cd\,,\cd)$ and
$\h\G(\cd\,,\cd)$ satisfying $(\ref{Riccati-G(s,t)})$.

\ms

\it Proof. \rm Let a partition $\D:0=t_0<\cds<t_N=T$ of $[0,T]$ be
given. By our condition and Proposition 6.2, we have the boundedness
of $\wt\G^\D(s,\t)$, $\G^\D(s,\t)$ and $\h\G^\D(s,\t)$. Then noting
that $\G^\D(\cd\,,\cd)$ etc. satisfy proper differential equations,
we have
\bel{le K}\left\{\ba{ll}
\ns\ds|\G^\D_s(s,\t)|+|\h\G^\D_s(s,\t)|\le K,\qq s\ne\t,\q\t\ne
t_k,\q0\le k\le N,\\
\ns\ds|\wt\G^\D_s(s,\t)|\le K,\qq s\ne\t,t_k,\q0\le k\le N,\\
\ns\ds|\dot P^\D(s)|+|\dot{\h P^\D}(s)|\le K,\qq s\ne t_k,\q0\le
k\le N,\ea\right.\ee
with some uniform constant $K>0$ (independent of partition $\D$).
The above implies that
$$|\G^\D(s,\t)-\G^\D(\bar s,\t)|+|\h\G^\D(s,\t)-\h\G^\D(\bar s,\t)|
\le K|s-\bar s|,\qq s,\bar s\in[0,T],\q\t\in[0,T].$$
On the other hand, by Proposition 6.3, for any $\t,\bar\t\in[0,T]$
with
$$|\t-\bar\t|<\min_{0\le k\le N-1}|t_{k+1}-t_k|,$$
we have
$$|\wt\G^\D(s,\t)-\wt\G^\D(s,\bar\t)|+|\G^\D(s,\t)-\G^\D(s,\bar\t)|
+|\h\G^\D(s,\t)-\h\G^\D(s,\bar\t)|\le\bar\o(\|\D|\|),$$
with $\bar\o(\cd)$ being a modulus of continuity. Then by Lemma 6.4,
we have
$$\lim_{\|\D\|\to0}\(|\wt\G^\D(s,\t)-\wt\G(s,\t)|+|\G^\D(s,\t)-\G(s,\t)|
+|\h\G^\D(s,\t)-\h\G(s,\t)|\)=0,$$
for some continuous functions $\wt\G(\cd\,,\cd)$, $\G(\cd\,,\cd)$,
and $\h\G(\cd\,,\cd)$.

\ms

Next, let $s\in[t_{k+1},t_{k+2})$. Then
$$\ba{ll}
\ns\ds|P^\D(s)-\G^\D(s,s)|\le|P_k(s)-P_k(t_{k+1})|+|\G_k(s,s)-\G_k(s,t_{k+1})|\le
K\|\D\|,\\
\ns\ds|\h P^\D(s)-\h \G^\D(s,s)|\le|\h P_k(s)-\h
P_k(t_{k+1})|+|\h\G_k(s)-\h\G_k(t_{k+1})|\le K\|\D\|,\ea$$
and for $0\le\t\le s\le T$, let $\t\in[t_{k+1},t_{k+2})$, and
$s\in[t_{\ell+1},t_{\ell+2})$, with $\ell\ge k$, then
$$\ba{ll}
\ns\ds|\wt\G^\D(s,\t)-\G^\D(s,\t)|=|\wt\G^\D_k(s)-\G^\D_k(s)|\le|\wt\G^\D_k(s)-\wt\G^\D_k(t_{\ell+2})|+|\G^\D_k(t_{\ell+2})-\G^\D_k(s)|\le
K\|\D\|.\ea$$
Hence, we have
$$\lim_{\|\D\|\to0}\(|\wt\G^\D(s,\t)-\G(s,\t)|+|P^\D(s)-\G(s,s)|+|\h
P^\D(s)-\h\G(s,s)|\)=0.$$
Consequently,
$$\lim_{\|\D\|\to0}\(\Th^\D(s)-\Th(s)|+|\h\Th^\D(s)-\h\Th(s)|\)=0,$$
with
$$\left\{\2n\ba{ll}
\ns\ds\Th(s)=\big[R(s,s)+D(s)\G(s,s)D(s)\big]^{-1}\big[B(s)^T\G(s,s)+D(s)^T\G(s,s)C(s)\big],\\
\ns\ds\h\Th(s)=\big[\h R(s,s)+\h D(s)\G(s,s)\h D(s)\big]^{-1}\big[\h B(s)^T\h\G(s,s)+\h D(s)^T\G(s,s)\h C(s)\big].\\
\ea\right.$$
Then we see that $\G(\cd\,,\cd)$ and $\h\G(\cd\,,\cd)$ satisfy
(\ref{Riccati-G(s,t)}). \endpf

\ms

Note that the above convergence result gives an existence of a
solution $(\G(\cd\,,\cd),\h\G(\cd\,,\cd))$ to the Riccati equation
system (\ref{Riccati-G(s,t)}). To be complete, we have the following
result.

\ms

\bf Theorem 6.6. \sl Let {\rm(H1)}--{\rm(H2)} hold. Then
$(\ref{Riccati-G(s,t)})$ admits a unique solution.

\ms

\it Proof. \rm We need only to prove the uniqueness. We denote
$$\L(s,t)=\left(\2n\ba{c}\G(s,t)\\ \h\G(s,t)\ea\2n\right),\qq0\le t\le s\le T,\qq
\L_0(\t)=\left(\2n\ba{ll}G(\t)\\ \h G(\t)\ea\2n\right).$$
and rewrite (\ref{Riccati-G(s,t)}) as follows:
$$\left\{\2n\ba{ll}
\ns\ds\L_s(s,\t)=F(s,\t,\L(s,\t),\L(s,s)),\qq0\le\t\le s\le T,\\
\ns\ds\L(T,\t)=\L_0(\t),\ea\right.$$
where $F(s,\t,\g,\bar\g)$ is some continuous map which is
differentiable in $(\g,\bar\g)$. Now, suppose there are two
solutions $\L^1(\cd\,,\cd)$ and $\L^2(\cd\,,\cd)$. Then
$$\h\L(\cd\,,\cd)=\L^1(\cd\,,\cd)-\L^2(\cd\,,\cd)$$
satisfies
$$\left\{\2n\ba{ll}
\ns\ds\h\L_s(s,\t)=F_\g(s,\t)\h
L(s,\t)+F_{\bar\g}(s,\t)\h\L(s,s),\qq0\le\t\le s\le T,\\
\ns\ds\h\L(T,\t)=0,\ea\right.$$
for some continuous (matrix-valued) functions $F_\g(s,\t)$ and
$F_{\bar\g}(s,\t)$. If we let $\h\F(\cd\,,\cd)$ be the fundamental
matrix of $F_\g(\cd\,,\cd)$, then
$$\h\L(s,\t)=\int_s^T\h\F(r,s)\F_{\bar\g}(r,\t)\h\L(r,r)dr,\qq
0\le\t\le s\le T.$$
Set $s=\t$ in the above, we end up with the following
$$\h\L(\t,\t)=\int_\t^T\h\F(r,s)\F_{\bar\g}(r,\t)\h\L(r,r)dr,\qq
0\le\t\le T.$$
This is a linear homogeneous Volterra integral equation with a continuous
kernel. Hence, it is necessary that
$$\h\L(s,s)=0,\qq s\in[0,T].$$
Then,
$$\h\L(s,\t)=0,\qq0\le\t\le s\le T,$$
which leads to the uniqueness of the solution of
(\ref{Riccati-G(s,t)}). \endpf

\ms

Combining the above, we can easily obtain a proof of Theorem 4.5.

\rm

\end{document}